\documentclass[10pt]{amsart}
\usepackage{amssymb,amsmath}
\usepackage{stmaryrd}
\usepackage{algorithm}
\usepackage{algorithmic}
\usepackage{mathbbol,bm,tensor}
\usepackage{color}
\usepackage{epsfig}
\usepackage{subfigure}
\usepackage{overpic}
\usepackage[margin=30mm]{geometry}

\usepackage{booktabs}
\usepackage{amsmath}
\usepackage{setspace}
\usepackage{array,caption}
\usepackage[labelfont=bf]{caption}
\usepackage[numbers,sort&compress]{natbib}

\newcommand{\ep}{\epsilon}

\newcommand{\avg}[1]{\{\!\!\{{#1}\}\!\!\}}

\newcommand{\vertiii}[1]{{\left\vert\kern-0.25ex\left\vert\kern-0.25ex\left\vert #1
    \right\vert\kern-0.25ex\right\vert\kern-0.25ex\right\vert}}
	
\makeatletter \@addtoreset{equation}{section} \makeatother

\theoremstyle{plain}
\newtheorem{theorem}{Theorem}[section]
\newtheorem{lemma}{Lemma}[section]
\newtheorem{proposition}{Proposition}[section]

\theoremstyle{definition}
\newtheorem{remark}{Remark}[section]

\title[A HDG method for the Peterlin problem]{A Linear HDG scheme for the diffusion type Peterlin viscoelastic problem}
\author[Sibang Gou]{Sibang Gou}
\address{Graduated School of Informatics, Kyoto University, Kyoto, 6068501, Japan.}
\email{sibang.gou@outlook.com}
\author[Jingyan Hu]{Jingyan Hu}
\address{Institute of Fundamental and Frontier Sciences, University of Electronic Science and Technology of China, Chengdu, 610051, China}
\email{jingyan\_hu@std.uestc.edu.cn}
\author[Qi Wang]{Qi Wang}
\address{Institute of Fundamental and Frontier Sciences, University of Electronic Science and Technology of China, Chengdu, 610051, China}
\email{qi\_wang@std.uestc.edu.cn}
\author[Feifei Jing]{Feifei Jing}
\address{School of Mathematics and Statistics, Northwestern Polytehcnical University, Xi'an, 710129, China}
\email{ffjing@nwpu.edu.cn}
\author[Guanyu Zhou]{Guanyu Zhou}
\address{Institute of Fundamental and Frontier Sciences, School of Mathematical Sciences, University of Electronic Science and Technology of China, Chengdu, 610051, China}
\email{wind\_geno@live.com}
\date{\today}
\subjclass[2010]{}
\thanks{The research of Guanyu Zhou was partially supported by NSFC General Project (No. 12171071) and the Natural Science Foundation of Sichuan Province (No. 2023NSFSC0055).
\\
The research of Feifei Jing was supported by NSFC (No.12371407).}

\begin{document}
% abstract
\begin{abstract}
A linear semi-implicit hybridizable discontinuous Galerkin (HDG) scheme is proposed to solve the diffusive Peterlin viscoelastic model, allowing the diffusion coefficient $\ep$ of the conformation tensor to be arbitrarily small. 
We investigate the well-posedness, stability, and error estimates of the scheme. In particular, we demonstrate that the $L^2$-norm error of the conformation tensor  is independent of the reciprocal of $\ep$. 
Numerical experiments are conducted to validate the theoretical convergence rates. Our numerical examples show that the HDG scheme performs better in preserving the positive definiteness of the conformation tensor compared to the ordinary finite element method (FEM).  
\end{abstract}
% title
\maketitle
%%
%%\tableofcontents
%%%%%%%%%%%%%%%%%%%%%%%%%%%%%%%%%%%%
%%%%%%%%%%%%%%%%%%%%%%%%%%%%%%%%%%%%
%%%%%%%%%%%%%%%%%%%%%%%%%%%%%%%%%%%%
%%%%%%%%%%%%%%%%%%%%%%%%%%%%%%%%%%%%
\section{Introduction}\label{sec:1}

Unlike Newtonian fluids, non-Newtonian fluids do not have a linear relationship between the stress tensor and the deformation tensor. In these fluids, the total stress tensor is typically represented as the sum of the Newtonian viscous stress and an additional stress component that arises from the presence of polymers \cite{M1}. 
Non-Newtonian fluids can be categorized into several classes, including shear-thinning, shear-thickening, yield-stress fluids, and viscoelastic fluids \cite{H2016}. Among these, viscoelastic fluids exhibit both elastic and viscous behaviors. When subjected to shear or stretching, they deform elastically and are capable of recovering partially or fully once the applied force is removed. This dual property makes viscoelastic fluids of significant interest in both engineering applications and scientific research.
Viscoelastic fluids are used to model a wide range of materials, including molten plastics, biological fluids such as blood and egg whites, and lubricants with polymer additives \cite{Lions2000, H2016}. Various mathematical models have been developed to describe the behavior of these fluids, including the Oldroyd model \cite{Guillope90}, the Finite Extensible Nonlinear Elastic-Peterlin (FENE-P) model \cite{Zhang2013}, and the generalized Peterlin model \cite{Medivdova2017}. 

In this study, we focus on the diffusive Peterlin viscoelastic system with a constant viscosity $\nu$ and a tiny elastic diffusion coefficient $\epsilon$ for the conformation tensor:
\begin{subequations}\label{eq:P}
	\begin{align}
		\bm{u}_t +\nabla\cdot(\bm{u}\otimes \bm{u})-\nu \Delta \bm{u}+\nabla p -\nabla \cdot [({\rm tr}\bm{C})\bm{C}] &=\bm{f} && \quad {\rm in}\quad\Omega\times(0,T),\label{eq:P-a} \\
		\nabla \cdot \bm{u}&=0\ && \quad {\rm in}\quad\Omega\times(0,T),\label{eq:P-b} \\
		\bm{C}_t + \nabla\cdot(\bm{C}\otimes\bm{u})-\epsilon \Delta \bm{C}-(\nabla \bm{u})\bm{C}-\bm{C}(\nabla\bm{u})^\top -{\rm tr}\bm{C}{\bm I}+ ({\rm tr}\bm{C})^2 \bm{C} &=0  &&\quad {\rm in}\quad\Omega\times(0,T),\label{eq:P-c} \\
		\bm{u}=0,\quad \frac{\partial\bm{C}}{\partial\bm{n}}& =0 &&\quad {\rm on}\quad\Gamma\times(0,T),\label{eq:P-d} \\
		\bm{u}(\cdot,0)=\bm{u}^0,\quad \bm{C}(\cdot,0) &= \bm{C}^0 &&\quad {\rm in}\quad \Omega,\label{eq:P-e} 
	\end{align}
\end{subequations}
where $\Omega$ is a smooth, bounded domain in $\mathbb{R}^2$ with boundary $\Gamma$. $\bm{f} \in L^2(\Omega)$ is the given force, $(\bm{u}^0,\bm{C}^0)$ is the initial value, ${\bm I}$ is the identity matrix in $\mathbb{R}^{2 \times 2}$ and the diffusive coefficient $\epsilon$ can be arbitrarily small. 
The unknown $(\bm{u},p,\bm{C}):\Omega\times[0,T]\rightarrow\mathbb{R}^2\times\mathbb{R}\times\mathbb{R}^{2\times 2}_{\rm sym}$ represents the velocity, pressure and conformation tensor, respectively. 
Here, the operator $\otimes$ denotes
$(\bm{u}\otimes\bm{v})_{ij}=u_i v_j, ~(\bm{C}\otimes\bm{u})_{ijk}= C_{ij}u_k$. If $\nabla\cdot\bm{u}=0$, we have
$\nabla\cdot(\bm{u}\otimes \bm{u})=(\bm{u}\cdot \nabla)\bm{u}+(\nabla\cdot\bm{u})\bm{u}=(\bm{u}\cdot \nabla)\bm{u}$ and $\nabla\cdot(\bm{C}\otimes \bm{u})=(\bm{u}\cdot \nabla)\bm{C}$.
Therefore, $\nabla\cdot(\bm{u}\otimes \bm{u})$ and $\nabla\cdot(\bm{C}\otimes \bm{u})$ in \eqref{eq:P-a} and \eqref{eq:P-c} can be replaced by the convective terms $(\bm{u}\cdot \nabla)\bm{u}$ and $(\bm{u}\cdot \nabla)\bm{C}$.

The diffusive Peterlin model is discussed in \cite{Lukacova15}, where the global existence of weak solutions for large data is proved. The existence of weak solutions in 3D is established in \cite{Brunk2022}, as well as a  conditional weak-strong uniqueness results. On the other hand, various discrete schemes have been proposed for the numerical simulation of the diffusion Peterlin model. 
The linear and nonlinear schemes of stabilized Lagrange-Galerkin FEM for the Oseen-type Peterlin viscoelastic model are proposed in \cite{M1,M2}. 
Three fully discrete schemes by using Newton's iterative, Picard's iterative and implicit-explicit time-stepping are presented in \cite{Jiang2018}.
\cite{Xia2023} introduces a linear decoupled scheme with stabilizing terms, where  $(\bm u, p)$ and each components of $\bm C$ can be computed in parallel. 
Numerical studies of other viscoelastic models (such as the Oldroyd-B model) also exist \cite{Ervin2004, Najib1995, Ravindran2020, Zheng2017}. Moreover, some log-conformation formulations \cite{Becker2023, Boyaval2009, Fattal2004, Hu2007, Wittschieber2022} are used to guarantee symmetric positive-definiteness of discrete conformation tensor.

%todo
When $\ep$ is very small, \eqref{eq:P-c} becomes convection-dominated.
However, the numerical schemes used in most existing works \cite{Jiang2018, Xia2023, Han23, Zhangyz2024} are not specifically designed to address the convection-dominated issue. To tackle this problem, \cite{M1, M2, Hana2015} propose the Lagrange-Galerkin method for the Oseen-type model. 
By carefully examining the existing results on error estimates, we find that the constant in the error bound often depends on the reciprocal of $\epsilon$.
As $\epsilon$ approaches $0$, the parabolic equation \eqref{eq:P-c} transitions to a hyperbolic form, which brings difficulty to numerical approximation.
To address the convection-dominated issue and to derive an error bound that does not depend on the reciprocal of $\epsilon$, we propose a linear HDG scheme for the diffusion Peterlin model and study the error analysis. 

%% todo
The HDG can reduce the computational cost of the DG method while preserving its attractive properties of conservation and stability \cite{Cockburn20092, Labeur2007}. 
The method transfers the computation from the global domain to each element, resulting in a smaller stiffness matrix dimension that facilitates efficient inversion operation\cite{HDG2019}. 
Currently the HDG method has been widely applied to elliptic problem \cite{Cockburn2008, Oikawa2015}, convection-diffusion equation \cite{Oikawa2014}, linear elasticity system \cite{Qiuwf2018}, Stokes equation \cite{Cockburn20122, Cockburn2011, Sander2017} and Navier-Stokes equations \cite{Nguyen2011, Sander2018}.
Especially, for the Navier-Stokes equations, the HDG scheme is mass and momentum conserving, energy stable, and pressure-robust. 
In addition, the discrete velocity is pointwise divergence-free, and the error estimates are optimal and independent of the viscosity $\nu$ \cite{Kirk2019}.
The motivation of this work is to employ the HDG method to address the convection-dominated issue of \eqref{eq:P}. In particular, we aim to obtain error estimates that do not depend on $\epsilon^{-1}$.

In this work, we design a linear HDG scheme for the diffusion Peterlin problem. We prove the unique existences of the discrete problem. 
Furthermore, by utilizing the Ritz projections for Stokes and Poisson equations, we establish the error estimates.
The numerical experiments are carried out to verify the theoretical convergence rates. 
The HDG scheme remains stable even for the cases with $\epsilon = 10^{-3}$ and $\epsilon = 0$. 
A comparison with the ordinary FEM in simulation shows that the HDG scheme performs better in preserving the positive definiteness of the conformation tensor. 

The remainder of this paper is organized as follows. 
Section \ref{sec:2} introduces the fully discrete semi-implicit linearizing HDG scheme which satisfies the mass conservation. 
Section \ref{sec:3} focuses on stability analysis.
The theoretical convergence rate of the error estimates is rigorously proved under proper assumptions in Section \ref{sec:4}. 
We provide several numerical examples in Section \ref{sec:ex} to study the experimental convergence rates and illustrate the stability.
%%%%%%%%%%%%%%%%%%%%%%%%%%%%%%%%%%%%%
%%%%%%%%%%%%%%%%%%%%%%%%%%%%%%%%%%%%%
\section{The Linear HDG scheme}\label{sec:2}
Let $\mathcal{T}:=\{K\}$ be a regular quasi-uniform triangulation of the domain $\Omega$. The boundary of the element $K$ is denoted by $\partial K$ and the outward unit normal vector on $\partial K$ by $\bm{n}$. Two adjacent cells $K^+$ and $K^-$ share an interior facet $F:=\partial K^+\bigcap \partial K^-$. A facet of $\partial K$ that lies on the boundary $\partial\Omega$ is called a boundary facet. The interior and boundary facet sets are denoted by $\mathcal{F}_I$ and $\mathcal{F}_B$, respectively. The set of all facets is denoted by $\mathcal{F}:=\mathcal{F}_I\bigcup\mathcal{F}_B$. We introduce the function spaces:  
\[
	\begin{aligned}\label{eq:space}
    \bm{V}_h:=&\{\bm{v}_h\in[L^2(\mathcal{T})]^2:\bm{v}_h|_{K}\in[P_k(K)]^2(\forall K\in\mathcal{T})\},\\
    \hat{\bm{V}}_h:=&\{\hat{\bm{v}}_h\in[L^2(\mathcal{F})]^2:\hat{\bm{v}}_h|_{F}\in[P_k(F)]^2(\forall F\in\mathcal{F}),\hat{\bm{v}}_h=0\text{ on }\mathcal{F}_B\},\\
    Q_h:=&\{q_h\in L^2(\mathcal{T}):q_h|_{K}\in P_{k-1}(K)(\forall K\in\mathcal{T})\},\quad
    \hat{Q}_h:=\{\hat{q}_h\in L^2(\mathcal{F}):\hat{q}_h|_{F}\in P_k(F)(\forall F\in\mathcal{F})\},\\
    \bm{W}_h:=&\{\bm{D}_h\in[L^2(\mathcal{T})]^{2\times 2}_{\rm sym}:\bm{D}_h|_{K}\in[P_k(K)]^{2\times 2}_{\rm sym}(\forall K\in\mathcal{T})\},\\
    \hat{\bm{W}}_h:=&\{\hat{\bm{D}}_h|_{F}\in[L^2(\mathcal{F})]^{2\times 2}_{\rm sym}:\hat{\bm{D}}_h\in[P_k(F)]^{2\times 2}_{\rm sym}(\forall F\in\mathcal{F})\},
    \end{aligned}
\]
where $P_{k}(K)$ denotes the space of polynomials of degree $k$ on $K$. Note that the spaces $\bm{V}_h,Q_h,\bm{W}_h$ are defined on the triangulation $\mathcal{T}$, whereas the spaces $\hat{\bm{V}}_h,\hat{Q}_h,\hat{\bm{W}}_h$ are only defined on the facets $\mathcal{F}$. We write $h_K$ for the characteristic length of $K$, and $h:=\max_{K\in\mathcal{T}}h_K$ is the mesh size.

The spaces $\bm{V}_h,Q_h,\bm{W}_h$ are discontinuous across cell boundaries. Hence, the trace of a function $\bm{a}\in\bm{V}_h$ may be double-valued on cell boundaries. At an interior facet, $F$, we denote the traces of $\bm{v}_h\in\bm{V}_h$ by $\bm{v}_h^+:=\bm{v}_h|_{K^+}$ and $\bm{v}_h^-:=\bm{v}_h|_{K^-}$. We introduce the jump operator $[\![\bm{v}_h]\!]:=\bm{v}_h^+\cdot\bm{n}^++\bm{v}_h^-\cdot\bm{n}^-$, where $\bm{n}^{\pm}$ represents the outward unit normal vector  on $\partial K^{\pm}$. 

The HDG scheme is designed by using the continuity of the flux and adding the following terms on the edge of each element \cite{Sander2018}. 
Before introducing the semi-implicit fully discretization scheme, we present additional symmetry, penalty and convection balance terms: 
\[
\begin{aligned}
	&\text{$\bullet$ \bf symmetry term: }  	\int_{\partial K}\nu(\hat{\bm{u}} - \bm{u})\otimes\bm{n}:\nabla\bm{v}~ds, \quad \int_{\partial K}\epsilon(\hat{\bm{C}} - \bm{C} )\otimes\bm{n}:\nabla\bm{D}~ds.\\
	& \text{$\bullet$ \bf penalty term: }   \int_{\partial K}\frac{\nu\alpha}{h_{K}}(\bm{u}-\hat{\bm{u}})\otimes\bm{n}:(\bm{v}-\hat{\bm{v}})\otimes\bm{n}~ds, \quad \int_{\partial K}\frac{\epsilon\beta}{h_{K}}(\bm{C}-\hat{\bm{C}})\otimes\bm{n}:(\bm{D}-\hat{\bm{D}})\otimes\bm{n}~ds.\\
	&\begin{aligned}
		\text{$\bullet$ \bf convection balance term: } &\int_{\partial K}\frac{\bm{u}\cdot\bm{n}}{2}(\bm{u}+\hat{\bm{u}})(\bm{v}-\hat{\bm{v}})~ds +\int_{\partial K}\frac{|\bm{u}\cdot\bm{n}|}{2}(\bm{u}-\hat{\bm{u}})(\bm{v}-\hat{\bm{v}})~ds, \\
		&\int_{\partial K}\frac{\bm{u}\cdot\bm{n}}{2}(\bm{C}+\hat{\bm{C}})(\bm{D}-\hat{\bm{D}})~ds +\int_{\partial K}\frac{|\bm{u}\cdot\bm{n}|}{2}(\bm{C}-\hat{\bm{C}})(\bm{D}-\hat{\bm{D}})~ds.
	\end{aligned}
\end{aligned}
\]
%

% For all $\bm{u},\bm{v}, \bm{w}\in \bm{V}_h$, $\hat{\bm{u}},\hat{\bm{v}}\in \hat{\bm{V}}_h$, $p,q \in Q_h$, $\hat p,\hat q \in \hat Q_h$, $\bm{C},\bm{D}\in \bm{W}_h$ and $\hat{\bm{C}},\hat{\bm{D}}\in \hat{\bm{W}}_h$, 
We set  
\[
	\begin{aligned}
		&a_h((\bm{u},\hat{\bm{u}}),(\bm{v},\hat{\bm{v}})):=\sum_{K\in\mathcal{T}}\int_{K}\nu\nabla\bm{u}:\nabla\bm{v}\,dx
		-\sum_{K\in\mathcal{T}}\int_{\partial K}\nu\nabla\bm{u}:(\bm{v}-\hat{\bm{v}})\otimes\bm{n}\,ds\notag\\
		& \qquad\qquad\qquad\qquad \quad
		-\sum_{K\in\mathcal{T}}\int_{\partial K}\nu(\bm{u}-\hat{\bm{u}})\otimes\bm{n}:\nabla\bm{v}\,ds
		+\sum_{K\in\mathcal{T}}\int_{\partial K}\frac{\nu\alpha}{h_{K}}(\bm{u}-\hat{\bm{u}})\otimes\bm{n}:(\bm{v}-\hat{\bm{v}})\otimes\bm{n}\,ds,\\
		&b_h((p,\hat{p}),(\bm{v},\hat{\bm{v}})):=-\sum_{K\in\mathcal{T}}\int_{K}p\nabla\cdot\bm{v}\,dx+\sum_{K\in\mathcal{T}}\int_{\partial K}(\bm{v}-\hat{\bm{v}}_h)\cdot\bm{n}\hat{p}\,ds,
		\\
		&o_h(\bm{w};(\bm{u},\hat{\bm{u}}),(\bm{v},\hat{\bm{v}}))
		:=-\sum_{K\in\mathcal{T}}\int_{K}\bm{u}\otimes\bm{w}:\nabla\bm{v}\,dx
		+\sum_{K\in\mathcal{T}}\int_{\partial K}\frac{\bm{w}\cdot\bm{n}}{2}(\bm{u}+\hat{\bm{u}})(\bm{v}-\hat{\bm{v}})\,ds\notag\\
		&\qquad\qquad\qquad\qquad\qquad
		+\sum_{K\in\mathcal{T}}\int_{\partial K}\frac{|\bm{w}\cdot\bm{n}|}{2}(\bm{u}-\hat{\bm{u}})(\bm{v}-\hat{\bm{v}})\,ds,
			\end{aligned}
\]
\[
	\begin{aligned}
		&A_h((\bm{C},\hat{\bm{C}}),(\bm{D},\hat{\bm{D}}))
		:=\sum_{K\in\mathcal{T}}\int_{K}\epsilon\nabla\bm{C}:\nabla\bm{D}\,dx
		-\sum_{K\in\mathcal{T}}\int_{\partial K}\epsilon\nabla\bm{C}:(\bm{D}-\hat{\bm{D}})\otimes\bm{n}\,ds\notag\\
		&\qquad\qquad\qquad\qquad\qquad
		-\sum_{K\in\mathcal{T}}\int_{\partial K}\epsilon(\bm{C}-\hat{\bm{C}})\otimes\bm{n}:\nabla\bm{D}\,ds
		+\sum_{K\in\mathcal{T}}\int_{\partial K}\frac{\epsilon\beta}{h_{K}}(\bm{C}-\hat{\bm{C}})\otimes\bm{n}:(\bm{D}-\hat{\bm{D}})\otimes\bm{n}\,ds.
	\end{aligned}
\]
We divide the time interval $I=[0,T]$ into $N$ equidistant time series $0=t^0<t^1<\cdots<t^N=T$. Set $\tau={T}/{N}$ and $t^n=n\tau$. 
Given the initial value $(\bm{u}_h^0,\bm{C}_h^0)\in \bm{V}_h\times\bm{W}_h$, for $n=0,1,2,\cdots,N$, find $(\bm{u}_h^{n+1},\hat{\bm{u}}_h^{n+1},p_h^{n+1},\hat{p}_h^{n+1},\bm{C}_h^{n+1},\hat{\bm{C}}_h^{n+1})\in  \bm{V}_h\times \hat{\bm{V}}_h\times Q_h\times\hat{Q}_h\times \bm{W}_h\times\hat{\bm{W}}_h,$ s.t.
\begin{subequations}\label{eq:Full}
	\begin{align}
		&\quad b_h((q_h,\hat{q}_h),(\bm{u}_h^{n+1},\hat{\bm{u}}_h^{n+1}))=0,\quad\forall q_h\in Q_h,\hat{q}_h\in \hat{Q}_h,\label{eq:Full-a}
		\\
		& \begin{aligned}
		   &\int_{\Omega}\frac{\bm{u}_h^{n+1}-\bm{u}_h^{n}}{\tau}\cdot\bm{v}_h\,d x+
		   a_h((\bm{u}_h^{n+1},\hat{\bm{u}}_h^{n+1}),(\bm{v}_h,\hat{\bm{v}}_h))
		   \\
		   &
		   +o_h(\bm{u}_h^n;(\bm{u}_h^{n+1},\hat{\bm{u}}_h^{n+1}),(\bm{v}_h,\hat{\bm{v}}_h))
		   +b_h((p_h^{n+1},\hat{p}_h^{n+1}),(\bm{v}_h,\hat{\bm{v}}_h))\\
		   = & 
		   \sum_{K\in\mathcal{T}}\int_{\partial K}{\rm tr}\bm{C}_h^{n+1}\bm{C}_h^n:(\bm{v}_h
		   -\hat{\bm{v}}_h)\otimes\bm{n}\,ds-\sum_{K\in\mathcal{T}}\int_{K}{\rm tr}\bm{C}_h^{n+1}\bm{C}_h^n:\nabla\bm{v}_h\,dx\\
		   & +\int_{\Omega}\bm{f}^{n+1}\cdot\bm{v}_h\,dx,\quad\forall\bm{v}_h\in\bm{V}_h,\hat{\bm{v}}_h\in\hat{\bm{V}}_h,
		\end{aligned}\label{eq:Full-b}\\
		& \begin{aligned}
		&\int_{\Omega}\frac{\bm{C}_h^{n+1}-\bm{C}_h^{n}}{\tau}:\bm{D}_h\,d x
		+A_h((\bm{C}_h^{n+1},\hat{\bm{C}}_h^{n+1}),(\bm{D}_h,\hat{\bm{D}}_h))+o_h(\bm{u}_h^n;(\bm{C}_h^{n+1},\hat{\bm{C}}_h^{n+1}),(\bm{D}_h,\hat{\bm{D}}_h))\\
		& 
		-\int_{\Omega}((\nabla \bm{u}_h^{n+1})\bm{C}_h^{n}+\bm{C}_h^{n}(\nabla\bm{u}_h^{n+1})^\top):\bm{D}_h\,dx
		\\
		=& \int_{\Omega}{\rm tr}\bm{C}_h^{n}{\bm I}:\bm{D}_h\,dx
		-\int_{\Omega}({\rm tr}\bm{C}_h^{n})^2\bm{C}_h^{n+1}:\bm{D}_h\,dx,
		\quad\forall\bm{D}_h\in \bm{W}_h,\hat{\bm{D}}_h\in\hat{\bm{W}}_h,
		\end{aligned}\label{eq:Full-c}
	\end{align}
\end{subequations}
where $\bm{f}^{n+1}=\bm{f}(t^{n+1})$ and $\alpha,\beta>0$ are penalty coefficients.
\begin{remark}\label{re:linear}
From \eqref{eq:Full-a}, it is easy to see that \cite[Proposition 1]{Sander2018}
\begin{subequations}\label{eq:mass-full}
	\begin{align}
		&\nabla\cdot\bm{u}_h^{n+1}=0,&\quad\forall x\in K, K\in\mathcal{T},\label{eq:mass-full-a}\\
		&[\![\bm{u}_h^{n+1}]\!]=\bm{u}_h^{n+1+}\cdot\bm{n}^{+}+\bm{u}_h^{n+1-}\cdot\bm{n}^{-}=0,&\quad \forall x\in F, F\in\mathcal{F}_I,\label{eq:mass-full-b}\\
		&\bm{u}_h^{n+1}\cdot\bm{n}=\hat{\bm{u}}_h^{n+1}\cdot\bm{n},&\quad \forall x\in F, F\in\mathcal{F}_B.\label{eq:mass-full-c}
	\end{align}
	\end{subequations}
\end{remark}

%%%%%%%%%%%%%%%%%%%%%%%%%%%%%%%%%%%%
%%%%%%%%%%%%%%%%%%%%%%%%%%%%%%%%%%%%
%%%%%%%%%%%%%%%%%%%%%%%%%%%%%%%%%%%%

\section{Stability}\label{sec:3}
\subsection{Preliminaries}\label{sec:3-1}
We denote by $\|\cdot\|_{L^p}$, $\|\cdot\|_{H^k}$ the norm and $|\cdot|_k$ the semi-norm of $L^p(\Omega)$ (or $L^p(\Omega)^2$, $L^p(\Omega)^{2\times 2}$) and $H^k(\Omega)$ (or $H^k(\Omega)^2$, $H^k(\Omega)^{2\times 2}$) recording to \cite{FEM2008}, respectively. The inner-product of $L^2(\Omega)$ (or $L^2(\Omega)^2$, $L^2(\Omega)^{2\times 2}$) is presented by $(\cdot,\cdot)$. The notation with subscript means the symmetric matrix or matrix-valued function. The symbol $C$ (or $c_i,C_i$) expresses a generic constant independent of $\tau,h,\epsilon$ and $n$. Here, we define he Bochner space $l^2(0,T;X)$ and $l^\infty(0,T;X)$ with following norms in mesh skeleton and time interval $[0,T]$: 
\[
% \begin{subequations}\label{eq:norm}
	\begin{aligned}
		&\|\bm{v}-\hat{\bm{v}}\|_{L^2(\mathcal{F})}^2:=\sum_{K\in\mathcal{T}}\|\bm{v}-\hat{\bm{v}}\|_{L^2(\partial K)}^2\quad \forall \bm{v}\in L^2(\mathcal{T}),\hat{\bm{v}}\in L^2(\mathcal{F}),\\
		&\|\bm{v}\|_{l^2(0,T;X)}^2:=\sum_{n=0}^N \tau\|\bm{v}^n\|_X^2\quad\big(\bm{v}:=(\bm{v}^0,\bm{v}^1,\cdots,\bm{v}^N),\quad\bm{v}^n\in X,n=0,1,\cdots,N\big),\\
		&\|\bm{v}\|_{l^\infty(0,T;X)}:=\max_{n=0,\cdots,N}\|\bm{v}^n\|_X\quad\big(\bm{v}:=(\bm{v}^0,\bm{v}^1,\cdots,\bm{v}^N),\quad\bm{v}^n\in X,n=0,1,\cdots,N\big).
	\end{aligned}
% \end{subequations}
\]

%Further the function space $Z^m(0,T)$ can be defined as
%\begin{equation}\label{def:Z}
%	\begin{aligned}
%		&Z^m(0,T):\{\bm{u}\in H^j(0,T;H^{m-1}(\Omega)),j=0,1,\cdots,m;\|\bm{u}\|_{Z^m(0,T)}\leq\infty\},\\
%		&\|\bm{u}\|_{Z^m(0,T)}:=\max_{j=0,\cdots,m}\|\bm{u}\|_{H^j(0,T;H^{m-1}(\Omega))}.
%	\end{aligned}
%\end{equation}
% For convenience, we introduce some mesh-dependent norms \cite{Aycil2017}:
Following the definitions in \cite{Aycil2017}, we introduce some mesh-dependent norms:
\[
	\begin{aligned}
	   &
	   \interleave(\bm{v}_h,\hat{\bm{v}}_h)\interleave_{0,v}^2:=\|\bm{v}_h\|_{L^2}^2
 	   +\sum_{K\in\mathcal{T}} h_K(\|\hat{\bm{v}}_h\|_{L^2(\partial K)}^2
 	   +\|\bm{v}_h-\hat{\bm{v}}_h\|_{L^2(\partial K)}^2),&&\quad \forall \bm{v}_h\in \bm{V}_h,\hat{\bm{v}}_h\in \hat{\bm{V}}_h,
 	   \\
 	   &\interleave(\bm{v}_h,\hat{\bm{v}}_h)\interleave_{v}^2:=\|\nabla\bm{v}_h\|_{L^2}^2
 	   +\sum_{K\in\mathcal{T}}\alpha h_K^{-1}\|\bm{v}_h-\hat{\bm{v}}_h\|_{L^2(\partial K)}^2,&&\quad \forall \bm{v}_h\in \bm{V}_h,\hat{\bm{v}}_h\in \hat{\bm{V}}_h,\\
		&\interleave(\bm{D}_h,\hat{\bm{D}}_h)\interleave_{w}^2:=\|\nabla\bm{D}_h\|_{L^2}^2
		+\sum_{K\in\mathcal{T}}\beta h_K^{-1}\|\bm{D}_h-\hat{\bm{D}}_h\|_{L^2(\partial K)}^2, &&\quad \forall \bm{D}\in\bm{W}_h,\hat{\bm{D}}\in \hat{\bm{W}}_h,\\
		&\interleave(q_h,\hat{q}_h)\interleave_{q}^2:=\|q_h\|_{L^2}^2+\sum_{K\in\mathcal{T}}h_K\|\hat{q}\|_{L^2(\partial K)}^2,&&\quad \forall q_h\in Q_h,\hat{q}_h\in\hat{Q}_h,
		\\
		&\|\bm v_h\|_{1}^2:=\interleave(\bm v_h,\avg{\bm v_h})\interleave_v^2
		&&\quad \forall \bm{v}_h\in \bm{V}_h,
	\end{aligned}		
\]
where $\avg{\bm v_h}:=\frac{1}{2}(\bm v_h^+ + \bm v_h^-)$.
By the triangle inequality, for $e=\partial K^+\cap \partial K^-$, one can easily check that:
$$
\|\bm v_h^+ - \avg{\bm v_h}\|_{L^2(e)}
= \frac{1}{2}\|\bm v_h^+ - \bm v_h^-\|_{L^2(e)}
\le 
\frac{1}{2}(\|\bm v_h^+-\hat{\bm v}_h\|_{L^2(e)} + \|\bm v_h^- -\hat{\bm v}_h\|_{L^2(e)}).
$$
As a result, 
$$
\|\bm v_h\|_1 \le C\interleave(\bm v_h, \hat{\bm v}_h)\interleave_v.
$$
Next, we collect some inequalities and lemmas in the following which will be applied to the convergence analysis \cite{FEM2008, Hana2015}.

\begin{itemize}
\item \textit{Poinc\'{a}re’s inequality.}\label{th:KP} There exists a positive constant $C_P$ such that
	\begin{equation}\label{eq:KP}
		\|\nabla v\|_{L^2}^2\geq C_P\|v\|_{H^1}^2\quad(\forall v\in H_0^1(\Omega)).
	\end{equation}

\item \textit{Trace inequality.}\label{th:T} There exists a positive constant $C_T$ such that
	\begin{equation}\label{eq:T}
		\|v\|_{L^2(\partial K)}\leq C_T h_K^{-\frac{1}{2}}\left(\|v\|_{L^2(K)}^2+h_K^2\|\nabla v\|_{L^2(K)}^2\right)^{\frac{1}{2}}\quad(\forall v\in H^1(K)).
	\end{equation}

\item \textit{Inverse inequality.}\label{th:I}
Under the assumption that $q$ is quasi-uniform, there exists a positive constant $C_I$ such that
\begin{equation}\label{eq:I}
	\|\bm{v}_h\|_{W^{l,p}}\leq C_I h^{\min(0,\frac{d}{p}-\frac{d}{q})-l}\|\bm{v}_h\|_{L^q}\quad (p,q\in[1,\infty],\forall \bm{v}_h\in \bm{V}_h).
\end{equation}

\item \textit{Inverse trace inequality.}\label{th:IT} There exists a positive constant $c_I$ such that for all $K\in \mathcal T_h$
	\begin{equation}\label{eq:IT}
		\|v_h\|^2_{L^2(\partial K)}\leq c_I h_K^{-1}\|v_h\|^2_{L^2(K)}\quad(\forall v_h\in V_h).
	\end{equation}
By the triangle inequality,
\begin{equation*}
\begin{aligned}
\interleave(\bm v_h,\hat{\bm v}_h)\interleave_{0,v}^2
&\leq
\|\bm v_h\|^2_{L^2} + \sum_{K\in\mathcal{T}}h_K(\|\bm v_h\|_{L^2(\partial K)} + \|\bm v_h - \hat{\bm v}_h\|_{L^2(\partial K)})^2
+ \sum_{K\in\mathcal{T}}h_K\|\bm v_h - \hat{\bm v}_h\|_{L^2(\partial K)}^2 
\\
&\leq
C(\|\bm v_h\|^2_{L^2} +  \sum_{K\in\mathcal{T}}h_K\|\bm v_h - \hat{\bm v}_h\|_{L^2(\partial K)}^2) \quad (\text{by } \eqref{eq:IT}).
\end{aligned}
\end{equation*}
Therefore, the norm $\interleave\cdot\interleave_{0,v}^2$ has an equivalent definition
 \begin{equation}\label{eq:ITeq}
  \interleave(\bm v_h,\hat{\bm v}_h)\interleave_{0,v}^2 = \|\bm v_h\|^2_{L^2} +  \sum_{K\in\mathcal{T}}h_K\|\bm v_h - \hat{\bm v}_h\|_{L^2(\partial K)}^2 \quad (\forall \bm v_h \in V_h, \hat{\bm{v}}_h\in \hat{\bm{V}}_h).
 \end{equation}

\item \textit{Inf-sup condition.}\label{prop:inf-sup} There exists a positive constant $C_{IS}$ such that 
	\begin{equation}\label{eq:inf}
  	   \inf_{q_h\in Q_h,\hat{q}_h\in \hat{Q}_h}\sup_{\bm{v}_h\in \bm{V}_h,\hat{\bm{v}}_h\in \hat{\bm{V}}_h}
	   \frac{b_h((q_h,\hat{q}_h),(\bm{v}_h,\hat{\bm{v}}_h))}{\interleave(q_h,\hat{q}_h)\interleave_{q}\interleave(\bm{v}_h,\hat{\bm{v}}_h)\interleave_{v}}
	  \geq C_{IS}.
	\end{equation}
The proof is in \cite{Kirk2019} and \cite[Lemma 1]{Sander20182}.
Note that $C_P,C_T,C_I,c_I$ and $C_{IS}$ are independent of $h$.
\end{itemize}
%\textit{(Discrete Gronwall Inequality)}\label{le:Gron} Let $a_0,\Delta t$ be positive number s.t.
%	\[
%	\tau\leq\frac{1}{2a_0}
%	\]
%	Futher, $\{x^n\}_{n\geq0},\{y^n\}_{n\geq1},\{b^n\}_{n\geq 1},\{a_1^n\}_{n\geq 1}$ are some non-negative sequences, Suppuse
%	\begin{subequations}\label{eq:Gron}
%		\begin{align}
%			\bar{D}_{\tau}x^n+y^n&\leq a_0x^n+a_1x^{n-1}+b^n,\quad\forall n\geq1,\\
%			\intertext{Then for all $n\geq1$, it holds that:}
%			x^n+\Delta t\sum_{i=1}^{n}y^i&\leq\exp\{2a_0n\Delta t+\Delta t\sum_{i=1}^{n}a_i^i\}(x^0+\Delta t\sum_{i=1}^{n}b^i).
%		\end{align}
%	\end{subequations}
%	Where $\bar{D}_{\tau}$ stands for the backward difference operator
%	\[
%	\bar{D}_{\tau}x^n:=\frac{x^n-x^{n-1}}{\tau}.
%	\]

\begin{lemma}[{\cite[Proposition 3.6]{Hana2015}}]\label{le:terms}
	For all $\bm{u}_h\in\bm{V}_h,\bm{C}_h,\bm{D}_h\in\bm{W}_h$, the following equation holds true:
	\begin{equation}\label{eq:terms}
		({\rm tr}\bm{D}_h \bm{C}_h,\nabla\bm{u}_h)-\frac{1}{2}({\rm tr}\bm{D}_{h},{\rm tr}[(\nabla\bm{u}_h)\bm{C}_{h}+\bm{C}_{h}(\nabla\bm{u}_h)^{\top}])=0.
	\end{equation}
\end{lemma}

\begin{lemma}[Norm equivalence \cite{Hana2015}]\label{prop:N}
	 We denote all symmetric positive definite matrix functions in $\mathbb{R}^{d\times d}$ by $S_+^d$ $(d=2,3)$. For all $\bm{D}\in S_+^d$,
	\begin{equation}\label{eq:N}
		\|\bm{D}\|_{L^p(L^p)}^p\leq\|{\rm tr}\bm{D}\|_{L^p(L^p)}^p\leq d^{p-1}\|\bm{D}\|_{L^p(L^p)}^p \quad(p\geq2),
	\end{equation}
	where $\|\bm{D}\|_{L^p(L^p)}^p = \int_{\Omega}\|\bm D(x)\|_{L^p}^p~dx$ and $\|\bm D(x)\|_{L^p}$ represents the $p$-norm of the matrix $\bm D(x)$.
	For briefness, we will write $\|\bm{D}\|_{L^p}$ instead of $\|\bm{D}\|_{L^p(L^p)}$.
\end{lemma}

\begin{lemma}[{\cite[Lemma 4.3]{Sander2017}}]\label{le:cona}
	For all $\bm{u}_h,\bm{v}_h\in \bm{V}_h$ and  $\hat{\bm{u}}_h,\hat{\bm{v}}_h\in \hat{\bm{V}}_h$, 
	\begin{equation}\label{eq:cona}
		a_h((\bm{u}_h,\hat{\bm{u}}_h),(\bm{v}_h,\hat{\bm{v}}_h))
		\leq
		c\nu\interleave(\bm{u}_h,\hat{\bm{u}}_h)\interleave_{v}\interleave(\bm{v}_h,\hat{\bm{v}}_h)\interleave_{v}.
	\end{equation}
\end{lemma}

\begin{lemma}[{\cite[Propositions 3.4 and 3.5]{Aycil2017}}]\label{le:conv}
	For all $\bm{w}_{1h},\bm{w}_{2h},\bm{u}_h,\bm{v}_h\in \bm{V}_h$ and  $\hat{\bm{u}}_h,\hat{\bm{v}}_h\in \hat{\bm{V}}_h$,
	\begin{subequations}
	\begin{align}\label{eq:conv1}
	&o_h(\bm{w}_{1h};(\bm{u}_h,\hat{\bm{u}}_h),(\bm{v}_h,\hat{\bm{v}}_h))
		\leq
		c\|\bm{w}_{1h}\|_{L^\infty}\interleave(\bm{u}_h,\hat{\bm{u}}_h)\interleave_{0,v}\interleave(\bm{v}_h,\hat{\bm{v}}_h)\interleave_{v},
	\\
	&
	\begin{aligned}\label{eq:conv2}
		&|o_h(\bm{w}_{1h};(\bm{u}_h,\hat{\bm{u}}_h),(\bm{v}_h,\hat{\bm{v}}_h))
		-o_h(\bm{w}_{2h};(\bm{u}_h,\hat{\bm{u}}_h),(\bm{v}_h,\hat{\bm{v}}_h))|
		\\
		&\quad\leq
		c\|\bm{w}_{1h}-\bm{w}_{2h}\|_1
		\interleave(\bm{u}_h,\hat{\bm{u}}_h)\interleave_{v}\interleave(\bm{v}_h,\hat{\bm{v}}_h)\interleave_{v}.
	\end{aligned}
	\end{align}
	\end{subequations}
\end{lemma}
\begin{lemma}[{\cite[Proposition 3.6]{Aycil2017}}]
	For all $\bm{w}_h,\bm{u}_h,\in \bm{V}_h$ and  $\hat{\bm{u}}_h,\in \hat{\bm{V}}_h$,
	\begin{equation}\label{eq:cov}
	o_h(\bm{w}_h,(\bm{u}_h,\hat{\bm{u}}_h),(\bm{u}_h,\hat{\bm{u}}_h))
	=\sum_{K\in\mathcal{T}}\frac{1}{2}\int_{\partial K}|\bm{u}_h^{n}\cdot\bm{n}||\bm{u}_h-\hat{\bm{u}}_h|^2\,d s.
\end{equation}
\end{lemma}
\begin{proposition}\label{prop:EU}
	When $\bm C_h^{n}\in S^2_+$, there exists a unique solution $(\bm{u}_h^{n+1},\hat{\bm{u}}_h^{n+1},p_h^{n+1}$, $\hat{p}_h^{n+1},\bm{C}_h^{n+1},\hat{\bm{C}}_h^{n+1})$ to \eqref{eq:Full}.
\end{proposition}
\begin{proof}
We rewrite the scheme \eqref{eq:Full} as 
\begin{subequations}\label{eq:Full-eu}
	\begin{align}
		&b_h((q_h,\hat{q}_h),(\bm{u}_h^{n+1},\hat{\bm{u}}_h^{n+1}))=0,\quad\forall q_h\in Q_h,\hat{q}_h\in \hat{Q}_h,
		\label{eq:Full-eu-a}
		\\
		& \begin{aligned}
		   &\frac{1}{\tau}\int_{\Omega}\bm{u}_h^{n+1}\cdot\bm{v}_h\,d x+
		   a_h((\bm{u}_h^{n+1},\hat{\bm{u}}_h^{n+1}),(\bm{v}_h,\hat{\bm{v}}_h))
		   +o_h(\bm{u}_h^n;(\bm{u}_h^{n+1},\hat{\bm{u}}_h^{n+1}),(\bm{v}_h,\hat{\bm{v}}_h))
		   \\
		   &\quad +b_h((p_h^{n+1},\hat{p}_h^{n+1}),(\bm{v}_h,\hat{\bm{v}}_h))
		   -\sum_{K\in\mathcal{T}}\int_{\partial K}{\rm tr}\bm{C}_h^{n+1}\bm{C}_h^n:(\bm{v}_h
		   -\hat{\bm{v}}_h)\otimes\bm{n}\,ds
		   \\
		   &\quad +\sum_{K\in\mathcal{T}}\int_{K}{\rm tr}\bm{C}_h^{n+1}\bm{C}_h^n:\nabla\bm{v}_h\,dx
		   = \int_{\Omega}\bm{g}^{n+1}_h\cdot\bm{v}_h\,dx,\quad\forall\bm{v}_h\in\bm{V}_h,\hat{\bm{v}}_h\in\hat{\bm{V}}_h,
		\end{aligned}\label{eq:Full-eu-b}\\
		& \begin{aligned}
		&\frac{1}{\tau}\int_{\Omega}\bm{C}_h^{n+1}:\bm{D}_h\,d x
		+A_h((\bm{C}_h^{n+1},\hat{\bm{C}}_h^{n+1}),(\bm{D}_h,\hat{\bm{D}}_h))
		+o_h(\bm{u}_h^n;(\bm{C}_h^{n+1},\hat{\bm{C}}_h^{n+1}),(\bm{D}_h,\hat{\bm{D}}_h))\\
		& \quad -\int_{\Omega}((\nabla \bm{u}_h^{n+1})\bm{C}_h^{n}+\bm{C}_h^{n}(\nabla\bm{u}_h^{n+1})^\top):\bm{D}_h\,dx
		+\int_{\Omega}({\rm tr}\bm{C}_h^{n})^2\bm{C}_h^{n+1}:\bm{D}_h\,dx
		\\
		&=\int_{\Omega}\bm{G}_h^{n}:\bm{D}_h\,dx,
		\quad\forall\bm{D}_h\in \bm{W}_h,\hat{\bm{D}}_h\in\hat{\bm{W}}_h,
		\end{aligned}\label{eq:Full-eu-c}
	\end{align}
\end{subequations}
where $\bm{g}_h^{n} = \frac{1}{\tau}\bm u^{n}_h$ and $\bm{G}_h^{n} = \frac{1}{\tau}\bm C^{n}_h + \text{tr}\bm C_h^{n}\bm I$. 
The unique existence is equivalent to proving that $(\bm{u}_h^{n+1},\hat{\bm{u}}_h^{n+1},p_h^{n+1},\hat{p}_h^{n+1},\bm{C}_h^{n+1},\hat{\bm{C}}_h^{n+1}) = (\bm 0,\bm 0,0,0,\bm 0,\bm 0)$ when $(\bm{g}_h^{n},\bm{G}_h^{n}) = (\bm 0,\bm 0)$.
%
% We adopt a technique similar to that in Proposition 4.1 of \cite{M2} to establish our proof. 
%

First, we show $(\bm{u}_h^{n+1},\hat{\bm{u}}_h^{n+1}, \text{tr}\bm{C}_h^{n+1}, \text{tr}\hat{\bm{C}}_h^{n+1})= (\bm 0, \bm 0,\bm 0, \bm 0)$.
Substituting $(\bm v_h,\hat{\bm{v}}_h,q_h,\hat{q}_h,\bm{D}_h,\hat{\bm{D}}_h) =(\bm{u}_h^{n+1},\hat{\bm{u}}_h^{n+1},$ $-p^{n+1}_h,-\hat{p}^{n+1}_h,(\frac{1}{2}\text{tr}\bm{C}_h^{n+1})\bm I,(\frac{1}{2}\text{tr}\hat{\bm{C}}_h^{n+1})\bm I)$ into \eqref{eq:Full-eu}, together with \eqref{eq:cov}, we have
\begin{subequations}\label{eq:Full-eu-1}
	\begin{align}
	& \begin{aligned}\label{eq:Full-eu-1-a}
	& \frac{1}{\tau}\|\bm u_h^{n+1}\|_{L^2}^2
	+ \nu\|\nabla\bm u_h^{n+1}\|^2_{L^2} 
	- 2\sum_{K\in\mathcal{T}}\int_{\partial K}\nu\nabla\bm u_h^{n+1}:(\bm u_h^{n+1}-\hat{\bm u}_h^{n+1})\otimes\bm n ~ds
	\\ & 
	 + \sum_{K\in\mathcal{T}}\frac{\nu\alpha}{h_K}\|((\bm u_h^{n+1}-\hat{\bm u}_h^{n+1})\otimes\bm n)\|^2_{L^2(\partial K)}
	+ \sum_{K\in\mathcal{T}}\frac{|\bm{u}_h^{n}\cdot\bm{n}|}{2}\|\bm{u}_h^{n+1}-\hat{\bm{u}}_h^{n+1}\|^2_{L^2(\partial K)}
	\\ & 
	 -\sum_{K\in\mathcal{T}}\int_{\partial K}{\rm tr}\bm{C}_h^{n+1}\bm{C}_h^n:(\bm{u}^{n+1}_h -\hat{\bm u}^{n+1}_h)\otimes\bm{n}\,ds
	+ \sum_{K\in\mathcal{T}}\int_{K}{\rm tr}\bm{C}_h^{n+1}\bm{C}_h^n:\nabla\bm{u}^{n+1}_h\,dx = 0,
	\end{aligned}\\
	& \begin{aligned}\label{eq:Full-eu-1-b}
	& \frac{1}{2\tau}\|\text{tr}\bm C_h^{n+1}\|_{L^2}^2
	+ \frac{\ep}{2}\|\nabla\text{tr}\bm C_h^{n+1}\|^2_{L^2}
	- \sum_{K\in\mathcal{T}}\int_{\partial K}\ep\nabla\text{tr}\bm C_h^{n+1}\bm I:\text{tr}(\bm C_h^{n+1}-\hat{\bm C}_h^{n+1})\otimes\bm n~ds
	\\ &
	 + \sum_{K\in\mathcal{T}}\frac{\ep\beta}{2h_K}\|\text{tr}(\bm C_h^{n+1}-\hat{\bm C}_h^{n+1})\otimes\bm n\|^2_{L^2(\partial K)}
	+ \sum_{K\in\mathcal{T}}\frac{|\bm{u}_h^{n}\cdot\bm{n}|}{4}\|\text{tr}(\bm{C}_h^{n+1}-\hat{\bm{C}}_h^{n+1})\|^2_{L^2(\partial K)}
	\\ & 
	 - \frac{1}{2}({\rm tr}\bm{C}^{n+1}_{h},{\rm tr}[(\nabla\bm{u}^{n+1}_h)\bm{C}^{n}_{h}+\bm{C}^n_{h}(\nabla\bm{u}^{n+1}_h)^{\top}])
	+ \frac{1}{2}\|{\rm tr}\bm{C}^{n}_{h}{\rm tr}\bm{C}^{n+1}_{h}\|^2_{L^2}
	= 0.
	\end{aligned}
	\end{align}
\end{subequations}
Applying the Schwartz inequality, the Young inequality and the inverse trace inequality \eqref{eq:IT}, we get
\[
    \begin{aligned}
        &- 2\sum_{K\in\mathcal{T}}\int_{\partial K}\nu\nabla\bm u_h^{n+1}:(\bm u_h^{n+1}-\hat{\bm u}_h^{n+1})\otimes\bm n~ds
        -\sum_{K\in\mathcal{T}}\int_{\partial K}{\rm tr}\bm{C}_h^{n+1}\bm{C}_h^n:(\bm{u}^{n+1}_h -\hat{\bm u}^{n+1}_h)\otimes\bm{n}\,ds
		\\\geq&
		- 2\nu\sum_{K\in\mathcal{T}} \|\nabla\bm u_h^{n+1}\|_{L^2(\partial K)} \|(\bm u_h^{n+1}-\hat{\bm u}_h^{n+1})\otimes\bm n\|_{L^2(\partial K)}\\
		& - \sum_{K\in\mathcal{T}} \|{\rm tr}\bm{C}_h^{n+1}\bm{C}_h^n\|_{L^2(\partial K)} \|(\bm u_h^{n+1}-\hat{\bm u}_h^{n+1})\otimes\bm n\|_{L^2(\partial K)}
		\\\geq&
		-\gamma\nu\|\nabla\bm u_h^{n+1}\|^2_{L^2} 
		- c_I\sum_{K\in\mathcal{T}} \frac{\nu+\frac{1}{2}}{h_K \gamma}\|(\bm u_h^{n+1}-\hat{\bm u}_h^{n+1})\otimes\bm n\|^2_{L^2(\partial K)}
		- \frac{\gamma}{2}\|{\rm tr}\bm{C}_h^{n+1}{\rm tr}\bm{C}_h^n\|^2_{L^2}
		,
    \end{aligned}
\]
where $c_I$ is the constant of the inverse trace inequality \eqref{eq:IT},
and we have used $ \|{\rm tr}\bm{C}_h^{n+1}\bm{C}_h^n\|^2_{L^2}\leq  \|{\rm tr}\bm{C}_h^{n+1}{\rm tr}\bm{C}_h^n\|^2_{L^2}$ because of $\bm C_h^{n}\in S_+^2$.
Here and hereafter, $\gamma$ is an arbitrary positive constant.
\[
    \begin{aligned}
        &- \sum_{K\in\mathcal{T}}\int_{\partial K}\ep\nabla\text{tr}\bm C_h^{n+1}{\bm I}:\text{tr}(\bm C_h^{n+1}-\hat{\bm C}_h^{n+1})\otimes\bm n~ds
		\\ \geq &
		-\frac{\ep\gamma}{2}\|\nabla\text{tr}\bm C_h^{n+1}\|^2_{L^2} 
		- c_I\sum_{K\in\mathcal{T}}\frac{\ep}{2\gamma h_K}\|\text{tr}(\bm C_h^{n+1}-\hat{\bm C}_h^{n+1})\otimes\bm n\|^2_{L^2(\partial K)}
		.
    \end{aligned}
\]
Next, summing up \eqref{eq:Full-eu-1-a} and \eqref{eq:Full-eu-1-b} and using \eqref{eq:terms}, we get
\[
    \begin{aligned}
    &\frac{1}{\tau}\|\bm u_h^{n+1}\|_{L^2}^2
	+ \nu(1-\gamma)\|\nabla\bm u_h^{n+1}\|^2_{L^2} 
	+ \sum_{K\in\mathcal{T}}\frac{1}{h_K}(\nu\alpha-c_I(\frac{\nu}{\gamma}+\frac{1}{2\gamma}))\|((\bm u_h^{n+1}-\hat{\bm u}_h^{n+1})\otimes\bm n)\|^2_{L^2(\partial K)}
	\\ & 
	 + \sum_{K\in\mathcal{T}}\frac{|\bm{u}_h^{n}\cdot\bm{n}|}{2}\|\bm{u}_h^{n+1}-\hat{\bm{u}}_h^{n+1}\|^2_{L^2(\partial K)}
	%%%%%%%
	+\frac{1}{2\tau}\|\text{tr}\bm C_h^{n+1}\|_{L^2}^2
	+ \frac{\ep}{2}(1-\gamma)\|\nabla\text{tr}\bm C_h^{n+1}\|^2_{L^2}
	\\ &
	 + \sum_{K\in\mathcal{T}}\frac{\ep}{h_K}(\frac{\beta}{2}-\frac{c_I}{2\gamma})\|\text{tr}(\bm C_h^{n+1}-\hat{\bm C}_h^{n+1})\otimes\bm n\|^2_{L^2(\partial K)}
	+ \sum_{K\in\mathcal{T}}\frac{|\bm{u}_h^{n}\cdot\bm{n}|}{4}\|\text{tr}(\bm{C}_h^{n+1}-\hat{\bm{C}}_h^{n+1})\|^2_{L^2(\partial K)}
	\\ & 
	 + \frac{1}{2}(1-\gamma)\|{\rm tr}\bm{C}^{n}_{h}{\rm tr}\bm{C}^{n+1}_{h}\|^2_{L^2}
	= 0.
    \end{aligned}
\]
Take $\gamma<1$ and sufficiently large $\alpha$, $\beta$ such that $\nu\alpha-c_I(\frac{\nu}{\gamma}+\frac{1}{2\gamma})>0$ and $\frac{\beta}{2}-\frac{c_I}{2\gamma}>0$.
As a result, we have $\bm u_h^{n+1}=\bm 0$ and $\text{tr}\bm C_h^{n+1}=\bm 0$ which implies $\hat{\bm u}_h^{n+1}=\bm 0$ and $\text{tr}\hat{\bm C}_h^{n+1}=\bm 0$. 

The next step is to show that $(\bm C_h^{n+1},\hat{\bm C}_h^{n+1})=(\bm 0,\bm 0)$.
Substituting $(\bm D_h,\hat{\bm D}_h)=$ $(\bm C_h^{n+1},\hat{\bm C}_h^{n+1})$ in \eqref{eq:Full-eu-c}, similar to the calculation of \eqref{eq:Full-eu-1},
\[
	\begin{aligned}
	0 &= \frac{1}{\tau}\|\bm C_h^{n+1}\|_{L^2}^2
	+ \ep\|\nabla\bm C_h^{n+1}\|^2_{L^2}
	- \sum_{K\in\mathcal{T}}\int_{\partial K}2\ep\nabla\bm C_h^{n+1}:(\bm C_h^{n+1}-\hat{\bm C}_h^{n+1})\otimes\bm n~ds
	\\ &\quad
	+ \sum_{K\in\mathcal{T}}\frac{\ep\beta}{h_K}\|(\bm C_h^{n+1}-\hat{\bm C}_h^{n+1})\otimes\bm n\|^2_{L^2(\partial K)}
	+ \|{\rm tr}\bm{C}^{n}_{h}\bm{C}^{n+1}_{h}\|^2_{L^2} 
    +\sum_{K\in\mathcal{T}}\beta_0h_K\|\bm C^n_h\bm C^{n+1}_h\|^2_{L^2(\partial K)}
	\ (\text{by $\bm u_h^{n+1}=0$})
	\\ &
	\geq
	 \frac{1}{\tau}\|\bm C_h^{n+1}\|_{L^2}^2
	+ \ep(1-\gamma)\|\nabla\bm C_h^{n+1}\|^2_{L^2}
	+ \|{\rm tr}\bm{C}^{n}_{h}\bm{C}^{n+1}_{h}\|^2_{L^2}
	% \\ &\quad
	+ \sum_{K\in\mathcal{T}}\frac{\ep}{h_K}(\beta-\frac{c_I}{\gamma})\|(\bm C_h^{n+1}-\hat{\bm C}_h^{n+1})\otimes\bm n\|^2_{L^2(\partial K)}.
	\end{aligned}
\]
We can take $\gamma<1$ and select sufficiently large $\beta$ such that $\beta-\frac{c_I}{\gamma}>0$, which implies $(\bm C_h^{n+1},\hat{\bm C}_h^{n+1})=(\bm 0,\bm 0)$. 

Now, \eqref{eq:Full-eu-b} becomes $b_h( (p_h^{n+1},\hat{p}_h^{n+1}), (\bm v_h,\hat{\bm v}_h))= 0$, in view of \eqref{eq:inf}, we conclude $p_h^{n+1}=0$ and $\hat{p}_h^{n+1}=0$.
The proof is complete.

\end{proof}

\begin{remark}
According to norm equivalence (see Lemma~\ref{prop:N}), $\bm C_h^n\in S_+^2$ implies that $\|{\rm tr}\bm{C}_h^{n+1}\bm{C}_h^n\|_{L^2}\leq \|{\rm tr}\bm{C}_h^{n+1}{\rm tr}\bm{C}_h^n\|_{L^2}$, which is utilized in the above proof of the unique existence of $\bm C_h^{n+1}$.
However, we cannot guarantee that $\bm C_h^{n+1}\in S_+^2$.
    Similar to \cite{M1,M2} for the Oseen-type Peterlin model and \cite{Jiang2018, Xia2023, Han23, Zhangyz2024} for \eqref{eq:P}, the positive-definiteness preserving is not addressed.
    Note that, for Oseen-type model, \cite[Proposition 6.9]{Hana2015} shows the positive-definiteness of $\bm C_h^{n+1}$ when $\bm C_h^{n}$ is positive definite and we choose sufficiently small $\tau$ such that $\tau\leq c_{n}(\bm u_h^n,\bm C_h^{n})$.
    However, the constant $c_{n}(\bm u_h^n,\bm C_h^{n})$ dependents on the numerical solution and cannot be prescribed. 
    
    We emphasize that our focus is on addressing the convection-dominated problems that arise when $\epsilon$ is extremely small. Our numerical examples show that the HDG scheme performs significantly better in preserving positive definiteness compared to the ordinary FEM. 
    For the positive-definiteness preserving technique, one can consult \cite{ Becker2023,Fattal2004,Wittschieber2022}.     
\end{remark}

%%%%%%%%%%%%%%%%%%%%%%%%%%%%%%
\subsection{Stability of the HDG scheme}\label{sec:3-2}
We denote by $(\bm{u}_h,\hat{\bm{u}}_h,p_h,\hat{p}_h,\bm{C}_h,\hat{\bm{C}}_h)$ the solution of \eqref{eq:Full}, where $\bm{u}_h:=(\bm{u}_h^0,\bm{u}_h^1,\cdots,\bm{u}_h^N)$ and other terms are defined as the same way. We introduce the norms:
\[
\begin{aligned}
\|\bm{u}_h\|_{l^\infty(0,T;L^2)}:=\max_{n=0,\cdots,N}\|\bm{u}_h^n\|_{L^2}, \quad \|\bm{u}_h\|_{l^2(0,T;L^2)}^2:=\sum_{n=1}^{N}\tau\|\bm{u}_h^n\|_{L^2}^2.
\end{aligned}
\]

\begin{theorem}\label{th:sta-f} (Stability) Assume that $\alpha,\beta$ are sufficiently large (independent of $h,\tau$ and $\nu,\epsilon$). We have
\begin{equation}\label{eq:sta-f}
	\begin{aligned}
		&\|\bm{u}_h\|_{l^\infty(0,T;L^2)}^{2}
		+\sum_{n=1}^{N}\|\bm u_h^n - \bm u_h^{n-1}\|^2_{L2}
		+\|{\rm tr}\bm{C}_h\|_{l^\infty(0,T;L^2)}^{2}
		+\sum_{n=1}^{N}\|{\rm tr}\bm C_h^n - {\rm tr}\bm C_h^{n-1} \|^2_{L^2}
		\\
		& 
		\quad +\nu\|\nabla\bm{u}_h\|_{l^2(0,T;L^2)}^{2}
        +\frac{\nu}{h}\|\bm{u}_h-\hat{\bm{u}}_h\|_{l^2(0,T;L^2(\mathcal{F}))}^{2}
		+\frac{\epsilon}{h}\|{\rm tr}\bm{C}_h-{\rm tr}\hat{\bm{C}}_h\|_{l^2(0,T;L^2(\mathcal{F}))}^{2}
	    \\
		&
		\quad +\epsilon\|\nabla{\rm tr}\bm{C}_h\|_{l^2(0,T;L^2)}^{2}
		+\sum_{n=1}^{N}\tau\|\sqrt{|\bm{u}_h^{n-1}\cdot\bm{n}|}(\bm{u}_h^{n}-\hat{\bm{u}}_h^{n})\|_{L^2(\mathcal{F})}^{2}
		+\sum_{n=1}^{N}\tau\|{\rm tr}\bm{C}_h^{n}{\rm tr}\bm{C}_h^{n-1}\|_{L^2}^{2}
		\\
		&
		\quad +\sum_{n=1}^{N}\tau\|\sqrt{|\bm{u}_h^{n-1}\cdot\bm{n}|}({\rm tr}\bm{C}_h^{n}-{\rm tr}\hat{\bm{C}}_h^{n})\|_{L^2(\mathcal{F})}^{2}
		\leq C_1e^{C_2T}(\|\bm{u}_{h}^{0}\|_{L^2}^{2}+\|{\rm tr}\bm{C}_{h}^{0}\|_{L^2}^{2}
		+\|\bm{f}^{n+1}\|_{l^2(0,T;L^2)}^2),
	\end{aligned}
\end{equation}
where $C_1$ and $C_2$ are independent of $h$, $\tau$, $\nu$, $\ep$ and $n$.
Moreover,
\begin{equation}\label{eq:sta-fp}
\begin{aligned}
	& \tau\|p_h\|_{l^2(0,T;L^2)}^2+h\tau\|\hat{p}_h\|_{l^2(0,T;L^2(\mathcal{F}))}^2\\
	\leq &C_3e^{C_4T}(1+\|\bm{u}_h\|_{l^\infty(0,T;L^\infty)}^2)(\|\bm{u}_{h}^{0}\|_{L^2}^{2}+\|{\rm tr}\bm{C}_{h}^{0}\|_{L^2}^{2} +\|\bm{f}\|_{l^2(0,T;L^2)}^2),
\end{aligned}
\end{equation}
where $C_3$ and $C_4$ depends on $C_1$, $C_2$ and $\nu$.
\end{theorem}

\begin{proof}
Taking $(q_h,\hat{q}_h,\bm{v}_h,\hat{\bm{v}}_h)=(p_h^{n+1}, \hat{p}_h^{n+1},\bm{u}_h^{n+1},\hat{\bm{u}}_h^{n+1})$ into \eqref{eq:Full-a} and \eqref{eq:Full-b}, 
adding \eqref{eq:Full-a} and \eqref{eq:Full-b} yields
\begin{equation*}\label{eq:sta-set-a}
\begin{aligned}
	&\int_{\Omega}\frac{\bm{u}_h^{n+1}-\bm{u}_h^{n}}{\tau}\cdot\bm{u}_h^{n+1}\,d x
	+\sum_{K\in\mathcal{T}}\frac{|\bm{u}_h^{n}\cdot\bm{n}|}{2}\|\bm{u}_h^{n+1}-\hat{\bm{u}}_h^{n+1}\|^2_{L^2(\partial K)}
	+\nu\|\nabla\bm{u}_h^{n+1}\|^2_{L^2} 
	\\
	&
	\quad +\sum_{K\in\mathcal{T}}\frac{\nu\alpha}{h_{K}}\|\bm{u}_h^{n+1}-\hat{\bm{u}}_h^{n+1}\|^2_{L^2(\partial K)}
	-2\sum_{K\in\mathcal{T}}\int_{\partial K}\nu\nabla\bm{u}_h^{n+1}:(\bm{u}_h^{n+1}-\hat{\bm{u}}_h^{n+1})\otimes\bm{n}\,ds\\
	&=\sum_{K\in\mathcal{T}}\int_{\partial K}({\rm tr}\bm{C}_h^{n+1})\bm{C}_h^{n}:(\bm{u}_h^{n+1}-\hat{\bm{u}}_h^{n+1})\otimes\bm{n}\,ds
	-\sum_{K\in\mathcal{T}}\int_K  ({\rm tr}\bm{C}_h^{n+1})\bm{C}_h^{n}:\nabla\bm{u}_h^{n+1} \,d x
	+\int_{\Omega}\bm{f}:\bm{u}_h^{n+1}\,dx.
\end{aligned}
\end{equation*}
Analogously, taking $(\bm{D}_{h},{\hat{\bm{D}}}_{h})=(\frac{1}{2}{\rm tr}\bm{C}_{h}^{n+1}\bm I,\frac{1}{2}{\rm tr}{\hat{\bm{C}}}_h^{n+1}\bm I)$ into \eqref{eq:Full-c}, and taking trace of the conformation tensor, we derive
\begin{equation*}\label{eq:sta-set-b}
	\begin{aligned}
		&\frac{1}{2}\int_{\Omega}\frac{{\rm tr}\bm{C}_h^{n+1}-{\rm tr}\bm{C}_h^{n}}{\tau}{\rm tr}\bm{C}_h^{n+1}\,d x
		+\sum_{K\in\mathcal{T}}\frac{|\bm{u}_h^{n}\cdot\bm{n}|}{4}\|{\rm tr}\bm{C}_h^{n+1}-{\rm tr}\hat{\bm{C}}_h^{n+1}\|^2_{L^2(\partial K)} 
		+\frac{\ep}{2}\|{\rm tr}\nabla\bm{C}_h^{n+1}\|^2_{L^2}
		\\
		&
		\quad +\sum_{K\in\mathcal{T}}\frac{\epsilon\beta}{2h_K}\|{\rm tr}\bm{C}_h^{n+1}-{\rm tr}\hat{\bm{C}}_h^{n+1}\|^2_{L^2(\partial K)}
		-\sum_{K\in\mathcal{T}}\int_{\partial K}\epsilon({\rm tr}\bm{C}_h^{n+1}-{\rm tr}\hat{\bm{C}}_h^{n+1})(\bm{n}\cdot\nabla {\rm tr}\bm{C}_h^{n+1}) \,d s
		\\
		&\quad +\frac{1}{2}\|{\rm tr}\bm{C}_h^{n}{\rm tr}\bm{C}_h^{n+1}\|^2_{L^2}
		=\frac{1}{2}\int_{\Omega}{\rm tr}((\nabla \bm{u}_h^{n+1})\bm{C}_h^{n}+\bm{C}_h^{n}(\nabla\bm{u}_h^{n+1})^\top){\rm tr}\bm{C}_h^{n+1}\,dx
		+\frac{1}{2}\int_{\Omega}{\rm tr}\bm{C}_h^{n}{\rm tr}\bm{C}_h^{n+1}\,dx.
	\end{aligned}
\end{equation*}
%We can also infer that the convection term of $\bm{C}$ satisfies:
%
%\begin{equation}
%	o_h(\bm{u}_n^n,(\bm{u}_h^{n+1},\hat{\bm{u}}_h^{n+1}),(\bm{u}_h^{n+1},\hat{\bm{u}}_h^{n+1}))=\sum_{K\in\mathcal{T}}\frac{1}{2}\int_{\partial K}|\bm{u}_h^{n}\cdot\bm{n}||\bm{u}_h^{n+1}-\hat{\bm{u}}_h^{n+1}|^2\,d s.
%\end{equation}
Next, we sum up above equations, and obtain (by \eqref{eq:mass-full} and \eqref{eq:terms})
\begin{equation}\label{eq:sta-add1}
	    \begin{aligned}
	    	&\int_{\Omega}\frac{\bm{u}_h^{n+1}-\bm{u}_h^{n}}{\tau}\cdot\bm{u}_h^{n+1}\,d x
	    	+\frac{1}{2}\int_{\Omega}\frac{{\rm tr}\bm{C}_h^{n+1}-{\rm tr}\bm{C}_h^{n}}{\tau}{\rm tr}\bm{C}_h^{n+1}\,d x
	    	+\sum_{K\in\mathcal{T}}\frac{|\bm{u}_h^{n}\cdot\bm{n}|}{2}\|\bm{u}_h^{n+1}-\hat{\bm{u}}_h^{n+1}\|^2_{L^2(\partial K)}
	    	\\
	    	&\quad +\sum_{K\in\mathcal{T}}\frac{|\bm{u}_h^{n}\cdot\bm{n}|}{4}\|{\rm tr}\bm{C}_h^{n+1}-{\rm tr}\hat{\bm{C}}_h^{n+1}\|^2_{L^2(\partial K)}
	    	+\nu\|\nabla\bm{u}_h^{n+1}\|^2_{L^2}
	    	+\sum_{K\in\mathcal{T}}\frac{\nu\alpha}{h_{K}}\|\bm{u}_h^{n+1}-\hat{\bm{u}}_h^{n+1}\|^2_{L^2(\partial K)}
	    	\\
	    	&
	    	\quad +\frac{\ep}{2}\|{\rm tr}\nabla\bm{C}_h^{n+1}\|^2_{L^2}
	    	+\sum_{K\in\mathcal{T}}\frac{\epsilon\beta}{2h_K}\|{\rm tr}\bm{C}_h^{n+1}-{\rm tr}\hat{\bm{C}}_h^{n+1}\|^2_{L^2(\partial K)}
	    	+\frac{1}{2}\|{\rm tr}\bm{C}_h^{n}{\rm tr}\bm{C}_h^{n+1}\|^2_{L^2}
	    	\\
	    	=& 2\sum_{K\in\mathcal{T}}\int_{\partial K}\nu\nabla\bm{u}_h^{n+1}:(\bm{u}_h^{n+1}-\hat{\bm{u}}_h^{n+1})\otimes\bm{n}\,ds
	    	+\sum_{K\in\mathcal{T}}\int_{\partial K}\epsilon({\rm tr}\bm{C}_h^{n+1}-{\rm tr}\hat{\bm{C}}_h^{n+1})(\bm{n}\cdot\nabla {\rm tr}\bm{C}_h^{n+1}) \,d s
	    	\\
	    	&\quad +\sum_{K\in\mathcal{T}}\int_{\partial K}({\rm tr}\bm{C}_h^{n+1})\bm{C}_h^{n}:(\bm{u}_h^{n+1}-\hat{\bm{u}}_h^{n+1})\otimes\bm{n}\,ds
	    	+\frac{1}{2}\int_{\Omega}{\rm tr}\bm{C}_h^{n}{\rm tr}\bm{C}_h^{n+1}\,dx
	    	+\int_{\Omega}\bm{f}^{n+1}:\bm{u}_h^{n+1}\,dx=:\sum_{i=1}^{5}R_i.
	\end{aligned}
\end{equation}
%Firstly, we deal with the time difference terms on the left hand that
%\begin{subequations}\label{eq:sta-time}
%	\begin{align}
%		\int_{\Omega}\frac{\bm{u}_h^{n+1}-\bm{u}_h^{n}}{\tau}\cdot\bm{u}_h^{n+1}\,d x
%		&=\frac{1}{2}\int_{\Omega}\frac{|\bm{u}_h^{n+1}-\bm{u}_h^{n}|^2}{\tau}\,d x
%		+\frac{1}{2}\int_{\Omega}\frac{|\bm{u}_h^{n+1}|^2}{\tau}\,d x
%		-\frac{1}{2}\int_{\Omega}\frac{|\bm{u}_h^{n}|^2}{\tau}\,d x\\
%		\frac{1}{2}\int_{\Omega}\frac{{\rm tr}\bm{C}_h^{n+1}-{\rm tr}\bm{C}_h^{n}}{\tau}:{\rm tr}\bm{C}_h^{n+1}\,d x
%		&=\frac{1}{4}\int_{\Omega}\frac{|{\rm tr}\bm{C}_h^{n+1}-{\rm tr}\bm{C}_h^{n}|^2}{\tau}\,d x
%		+\frac{1}{4}\int_{\Omega}\frac{|{\rm tr}\bm{C}_h^{n+1}|^2}{\tau}\,d x
%		-\frac{1}{4}\int_{\Omega}\frac{|{\rm tr}\bm{C}_h^{n}|^2}{\tau}\,d x.
%	\end{align}
%\end{subequations}
Now we estimate the terms on the right hand. Using the inverse trace \eqref{eq:IT}, we have:
\[
% \label{eq:sta-norm}
	\begin{aligned}
		&R_1
		\leq
		\sum_{K\in\mathcal{T}}\frac{c_I\nu}{\gamma h_K}\|\bm{u}_h^{n+1}-\hat{\bm{u}}_h^{n+1}\|^2_{L^2(\partial K)}
		+\nu\gamma\|\nabla\bm{u}_h^{n+1}\|^2_{L^2},
		\\
		&R_2
		\leq
		\sum_{K\in\mathcal{T}}\frac{c_I\epsilon}{2\gamma h_K}\|{\rm tr}\bm{C}_h^{n+1}-{\rm tr}\hat{\bm{C}}_h^{n+1}\|^2_{L^2(\partial K)}
		+\frac{\ep\gamma}{2}\|\nabla {\rm tr}\bm{C}_h^{n+1}\|^2_{L^2},
		% \label{eq:sta-norm-b}
		\nonumber\\
		&R_3
		\leq 
		\sum_{K\in\mathcal{T}}\frac{c_I}{2\gamma h_K}\|\bm{u}_h^{n+1}-\hat{\bm{u}}_h^{n+1}\|^2_{L^2(\partial K)}
		+\frac{\gamma}{2}\|{\rm tr}\bm{C}_h^{n+1}{\rm tr}\bm{C}_h^n\|^2_{L^2} \ (\text{by } \eqref{eq:N}).
	\end{aligned}
\]
Similarly, $R_4,R_5$ are estimated as follows: 
\[
\begin{aligned}
		&R_4
		\leq
		\frac{1}{2}\|{\rm tr}\bm{C}_h^{n}\|_{L^2}^2
		+\frac{1}{8}\|{\rm tr}\bm{C}_h^{n+1}\|_{L^2}^2,
		\quad
		% \label{eq:sta-norm-d}
		R_5\leq
		\frac{1}{2}\|\bm{f}\|_{L^2}^2
		+\frac{1}{2}\|\bm{u}_h^{n+1}\|_{L^2}^2.
		% \label{eq:sta-norm-e}
		\nonumber\\
	\end{aligned}	
\]
Taking the above estimates of $\{R_i\}_{i=1}^5$ into the right hand side of \eqref{eq:sta-add1} and summing up with respect to $i$, we arrive:
\begin{equation}\label{eq:sta-add2}
	\begin{aligned}
		&\frac{1}{2\tau}(\|\bm{u}_h^{n+1}\|_{L^2}^2-\|\bm{u}_h^{n}\|_{L^2}^2)
		+\frac{1}{4\tau}(\|{\rm tr}\bm{C}_h^{n+1}\|_{L^2}^2-\|{\rm tr}\bm{C}_h^{n}\|_{L^2}^2)
		+\frac{1}{2\tau}\|\bm{u}_h^{n+1}-\bm{u}_h^{n}\|_{L^2}^2
		\\
		&
		\quad +\frac{1}{4\tau}\|{\rm tr}\bm{C}_h^{n+1}-{\rm tr}\bm{C}_h^{n}\|_{L^2}^2+\sum_{K\in\mathcal{T}}\frac{|\bm{u}_h^{n}\cdot\bm{n}|}{2}\|\bm{u}_h^{n+1}-\hat{\bm{u}}_h^{n+1}\|^2_{L^2(\partial K)}
		+\nu (1-\gamma)\|\nabla\bm{u}_h^{n+1}\|_{L^2}^2
		\\
		&
		\quad+\sum_{K\in\mathcal{T}}\frac{|\bm{u}_h^{n}\cdot\bm{n}|}{4}\|{\rm tr}\bm{C}_h^{n+1}-{\rm tr}\hat{\bm{C}}_h^{n+1}\|^2_{L^2(\partial K)}
		+\sum_{K\in\mathcal{T}}\frac{(\nu\alpha-c_I(\frac{\nu}{\gamma}+\frac{1}{2\gamma}))}{h_K}\|\bm{u}_h^{n+1}-\hat{\bm{u}}_h^{n+1}\|_{L^2(\partial K)}^2
		\\
		&
		\quad+\frac{\ep}{2}(1-\gamma)\|{\rm tr}\nabla\bm{C}_h^{n+1}\|_{L^2}^2 
		+\sum_{K\in\mathcal{T}}\frac{\epsilon }{2h_K}(\beta-\frac{c_I}{\gamma})\|{\rm tr}\bm{C}_h^{n+1}-{\rm tr}\hat{\bm{C}}_h^{n+1}\|_{L^2(\partial K)}^2
		\\
		&
		\quad+\frac{1}{2}(1-\gamma)\|{\rm tr}\bm{C}_h^{n+1}{\rm tr}\bm{C}_h^{n}\|_{L^2}^2
		\leq
		\frac{1}{2}\|{\rm tr}\bm{C}_h^{n}\|_{L^2}^2
		+\frac{1}{4}\|{\rm tr}\bm{C}_h^{n+1}\|_{L^2}^2
		+\frac{1}{2}\|\bm{u}_h^{n+1}\|_{L^2}^2
		+\frac{1}{2}\|\bm{f}^{n+1}\|_{L^2}^2.
	\end{aligned}
\end{equation}
Choose $\gamma<1$ and suﬀiciently large $\alpha$ and $\beta$ such that $\nu\alpha-c_I(\frac{\nu}{\gamma}+\frac{1}{2\gamma})>0$ and $\beta-\frac{c_I}{\gamma}>0$. 
Applying the discrete Gronwall inequality to \eqref{eq:sta-add2}, we conclude \eqref{eq:sta-f}.

Now we estimate $p_h$ and $\hat{p}_h$.
%$$
%\sigma\interleave(p_h^{n+1},\hat{p}_h^{n+1})\interleave_{q}
%\leq
%\sup_{\bm{v}_h\in \bm{V}_h,\hat{\bm{v}}_h\in \hat{\bm{V}}_h}
%\frac{b_h((p_h^{n+1},\hat{p}_h^{n+1}),(\bm{v}_h,\hat{\bm{v}}_h))}{\interleave(\bm{v}_h,\hat{\bm{v}}_h)\interleave_{v}},
%$$
In view of
$$
	\begin{aligned}
		b_h((p_h^{n+1},\hat{p}_h^{n+1}),(\bm{v}_h,\hat{\bm{v}}_h))
		&=
		\int_{\Omega}\frac{\bm{u}_h^{n+1}-\bm{u}_h^{n}}{\tau}:\bm{v}_h\,d x
		+a_h((\bm{u}_h^{n+1},\hat{\bm{u}}_h^{n+1}),(\bm{v}_h,\hat{\bm{v}}_h))\\
		&\quad+o_h(\bm{u}_h^n;(\bm{u}_h^{n+1},\hat{\bm{u}}_h^{n+1}),(\bm{v}_h,\hat{\bm{v}}_h))
		+\sum_{K\in\mathcal{T}}\int_{K}{\rm tr}\bm{C}_h^{n+1}\bm{C}_h^n:\nabla\bm{v}_h\,dx\\
		&\quad+\sum_{K\in\mathcal{T}}\int_{\partial K}{\rm tr}\bm{C}_h^{n+1}\bm{C}_h^n:(\bm{v}_h-\hat{\bm{v}}_h)\otimes\bm{n}\,ds
		+\sum_{K\in\mathcal{T}}\int_{K}\bm{f}^{n+1}:\nabla\bm{v}_h\,dx.\\
%		&\leq
%		\bigg(\frac{1}{\tau}\|\bm{u}_h^{n+1}-\bm{u}_h^{n}\|_{L^2}
%		+(\nu+c\|\nabla\bm{u}_h^{n}\|_{L^2})\interleave(\bm{u}_h^{n+1},\hat{\bm{u}}_h^{n+1})\interleave
%		_{v}\\
%	   &\quad+\|{\rm tr}\bm{C}_h^{n+1}{\rm tr}\bm{C}_h^{n}\|_{L^2}
%	   +\|\bm{f}^{n+1}\|_{L^2}\bigg)
%	\interleave(\bm{v}_h,\hat{\bm{v}}_h)\interleave_{v}.
	\end{aligned}
$$
By using Lemmas~\ref{le:cona} and \ref{le:conv}, one can validate that
$$
\begin{aligned}
	b_h((p_h^{n+1},\hat{p}_h^{n+1}),(\bm{v}_h,\hat{\bm{v}}_h))
	&\leq
	\Big(\frac{1}{\tau}\|\bm{u}_h^{n+1}-\bm{u}_h^{n}\|_{L^2}
	+\nu\interleave(\bm{u}_h^{n+1},\hat{\bm{u}}_h^{n+1})\interleave_{v}
	+c\|\bm{u}_h^{n}\|_{L^\infty}
	\interleave(\bm{u}_h^{n+1},\hat{\bm{u}}_h^{n+1})\interleave_{0,v}
	\\
	&\quad+\|{\rm tr}\bm{C}_h^{n+1}{\rm tr}\bm{C}_h^{n}\|_{L^2}
	+\|\bm{f}^{n+1}\|_{L^2}\Big)
	\interleave(\bm{v}_h,\hat{\bm{v}}_h)\interleave_{v}.
\end{aligned}
$$
Using the inf-sup condition \eqref{eq:inf}, we have
$$
\begin{aligned}
	\interleave(p_h^{n+1},\hat{p}_h^{n+1})\interleave_{q}^2
	&\leq C\Big(\frac{1}{\tau^2}\|\bm{u}_h^{n+1}-\bm{u}_h^{n}\|_{L^2}^2
	+\nu^2\interleave(\bm{u}_h^{n+1},\hat{\bm{u}}_h^{n+1})\interleave
	_{v}^2
	+\|\bm{u}_h^{n}\|_{L^\infty}^2
	\interleave(\bm{u}_h^{n+1},\hat{\bm{u}}_h^{n+1})\interleave_{0,v}^2
	\\
	&\quad
	+\|{\rm tr}\bm{C}_h^{n+1}{\rm tr}\bm{C}_h^{n}\|_{L^2}^2
	+\|\bm{f}^{n+1}\|_{L^2}^2\Big).\\
%	\Rightarrow\nu^{-1}\sum_{n=1}^{N}\tau\|p_h\|_{L^2}^2+\nu^{-1}h\sum_{n=1}^{N}\tau\|\hat{p}_h\|_{L^2(\mathcal{F})}^2
%	&\leq\sigma^{-2}
%	\sum_{i=0}^{n}\tau\left(\frac{\|\bm{u}_h^{i+1}-\bm{u}_h^{i}\|_{L^2}^2}{\tau}
%	+\nu\interleave(\bm{u}_h^{i+1},\hat{\bm{u}}_h^{i+1})\interleave
%	_{v}^2\right.\\
%	&\left.+\|\sqrt{\bm{u}_h^i\cdot\bm{n}}(\bm{u}_h^{i+1}-\hat{\bm{u}}_h^{i+1})\|_{L^2(\mathcal{F})}^2
%	+\|{\rm tr}\bm{C}_h^{i+1}{\rm tr}\bm{C}_h^{i}\|_{L^2}^2\right).
\end{aligned}
$$
In use of \eqref{eq:sta-f} and \eqref{eq:ITeq}, we see that
$$
\begin{aligned}
	\sum_{i=0}^{N}\tau^2\interleave(p_h^{i+1},\hat{p}_h^{i+1})\interleave_{q}^2
	&\leq
	C(1+\|\bm{u}_h\|_{l^\infty(0,T;L^\infty)}^2)(\|\bm{u}_{h}^{0}\|_{L^2}^{2}+\|{\rm tr}\bm{C}_{h}^{0}\|_{L^2}^{2}
	+\|\bm{f}\|_{l^2(0,T;L^2)}^2).
	% &\leq
	% C(1+\nu^{-2})(\|\bm{u}_{h}^{0}\|_{L^2}^{2}+\|{\rm tr}\bm{C}_{h}^{0}\|_{L^2}^{2}
	% +\|\bm{f}\|_{l^2(0,T;L^2)}^2)^2.
\end{aligned}
$$
Hence we conclude \eqref{eq:sta-fp}.
\end{proof}
% We can infer the boundary of $\|\bm{C}_h\|_{l^\infty(0,T;L^2)}$ with the stability of $\|{\rm tr}\bm{C}_h\|_{l^\infty(0,T;L^2)}$ and norm inequality \eqref{eq:N}, because of the symmetric positive definite matrices $\{\bm{C}_h^n\}_{n=1}^N$. Moreover, it is easy to certificate the linear HDG scheme exists an unique solution.
%\begin{remark}
%	We have also prove the stability of $theta$-method for time discretization. We partition the time interval I into an ordered series of time levels$0=t^0<t^1<\cdots <t^N$. The difference between each time level is denoted by$\Delta t^n=t^{n+1}-t^n$. To discretize in time, we consider the $\theta$-method and denote midpoint values of a function $y$ by $y^{n+\theta}=(1-\theta)y^n+\theta y^{n+1}$. A sufficient condition for the stability of the method is: $\alpha,\beta$ are large enough and $\theta\geq \frac{1}{2}$.
%\end{remark}
%%%%%%%%%%%%%%%%%%%%%%%%%%%%%%%%%%%%
%%%%%%%%%%%%%%%%%%%%%%%%%%%%%%%%%%%%
%%%%%%%%%%%%%%%%%%%%%%%%%%%%%%%%%%%%
\section{Error Analysis}\label{sec:4}
This section is devoted to the error estimates.
We assume that the exact solution of \eqref{eq:P} satisfies:
\begin{equation}\label{def:error-space}
    \begin{aligned}
        & \bm{u} \in C(0,T;H^2(\Omega)^2\cap W^{1,\infty}(\Omega)^2) \cap C^1(0, T ; H^{k+1}(\Omega)^2)\cap  H^2(0, T ; L^2(\Omega)^2), \\
        & p\in C(0,T;H^{k}(\Omega)\cap L^2(\Omega)), \\
        & \bm{C} \in  C(0,T;H^2(\Omega)_{\rm sym}^{2\times2}\cap W^{1,\infty}(\Omega)_{\rm sym}^{2\times2}) \cap C^1(0,T;H^{k+1}(\Omega)_{\rm sym}^{2\times2}) \cap H^2(0, T ; L^2(\Omega)_{\rm sym}^{2 \times 2}).
	\end{aligned}
\end{equation}

% \subsection{Projection operators}\label{sec:4-1}
%$$
%(\bm{u}_h,p_h,\bm{C}_h):=\{(\bm{u}_h^{n+1},p_h^{n+1},\bm{C}_h^{n+1})\}_{n=0}^{N-1}\subset\bm{V}_h\times Q_h\times \bm{W}_h
%$$
%and
%$$
%{(\bm{u},p,\bm{C})(t);t\in(0,T)}\subset\bm{V}\times Q\times \bm{W}
%$$
%where the function spaces are

%\textbf{Hypotheses 2} Assume $\bm{u}\in L^\infty(0,T;L^\infty(\Omega)^2)$ ,and
%\begin{equation}\label{eq:hypo2}
%	\|\bm{u}_h^0\|_{L^\infty}\leq C+\|\bm{u}\|_{L^\infty(0,T;L^\infty)}.
%\end{equation}

%We want to prove a proposition with hypothesis 2.
%
%\begin{proposition}\label{prop:l-infty}
%	($L^\infty$ norm boundary) Suppose hypothesis 2, and $\{\bm{u}_h^{n+1}\}$ are solutions of \eqref{eq:Full}. Let $k$ of the function space \eqref{eq:space} and $\tau$ satisfy that
%	\begin{equation}\label{eq:con}
%		\begin{aligned}
%			&k\geq1,\\
%			&\tau<h,
%		\end{aligned}
%	\end{equation}
%   we have the boundary estimation that
%	\begin{equation}
%		\|\bm{u}_h^{n+1}\|_{L^\infty}
%		\leq C+\|\bm{u}\|_{L^\infty(0,T;L^\infty)},
%	\end{equation}
%    which means $\|\bm{u}_h^{n+1}\|_{L^\infty}$ is boundary.
%\end{proposition}
%The proof of proposition \ref{prop:l-infty} is in section \ref{sec:4-2}.
We introduce the following projection operators for error analysis.
For all $(\bm{w},r)\in H_{0}^{1}(\Omega)^2\times L_{0}^{2}(\Omega)$ with $\nabla\cdot\bm{w}=0$, there exists a unique  $(\bm{w}_h,\hat{\bm{w}}_h,r_h,\hat{r}_h)\in\bm{V}_h\times\hat{\bm{V}}_h\times Q_h\times\hat{Q}_h$ such that, for all $(\bm v_h,\hat{\bm v}_h, q_h,\hat{q}_h)\in \bm{V}_h\times\hat{\bm{V}}_h\times Q_h\times\hat{Q}_h$,
\[
% \begin{subequations}\label{def:ritz-u}
	\begin{aligned}
		\nu(\nabla\bm{w},\nabla\bm{v}_h)
		-(r,\nabla\cdot\bm{v}_h)
		&= a_h((\bm{w}_h,\hat{\bm{w}}_h),(\bm{v}_h,\hat{\bm{v}}_h))
		-b_h((r_h,\hat{r}_h),(\bm{v}_h,\hat{\bm{v}}_h)),\nonumber
		\\	
		(q_h,\nabla\cdot\bm{w})
		&=b_h((q_h,\hat{q}_h),(\bm{w}_h,\hat{\bm{w}}_h)).
	\end{aligned}
% \end{subequations}
\]
We define the Ritz projections $(\mathcal{R} _h^{u},\mathcal{R}_h^{p}):H_{0}^{1}(\Omega)^2\times L_{0}^{2}(\Omega)\rightarrow\bm{V}_h\times Q_h$ with $(\mathcal{R} _h^{u}\bm{w},\mathcal{R}_h^{p}r)=(\bm{w}_h,r_h),$ and $(\hat{\mathcal{R}} _h^{u},\hat{\mathcal{R}}_h^{p}):H_{0}^{1}(\Omega)^2\times L_{0}^{2}(\Omega)\rightarrow\hat{\bm{V}}_h\times \hat{Q}_h$ with $(\hat{\mathcal{R}} _h^{u}\bm{w},\hat{\mathcal{R}}_h^{p}r)=(\hat{\bm{w}}_h,\hat{r}_h)$. According to  \cite[Lemma 5.2]{Sander2017}, we have the error estimates: for $(\bm{w},r)\in(H^{k+1}(\Omega)^2\times H^{k}(\Omega))\cap(H_{0}^{1}(\Omega)^2\times L_{0}^{2}(\Omega))$,
\[ 
%\begin{subequations}\label{eq:ritz-error-u}
	\begin{aligned}
		&\interleave(\bm{w}-\mathcal{R} _h^{u}\bm{w},\bm{w}-\hat{\mathcal{R}} _h^{u}\bm{w})\interleave_{v}^2
		\leq C h^{2k}|\bm{w}|_{k+1}^2,\nonumber\\
		&\interleave(r-\mathcal{R}_h^{p}r,r-\hat{\mathcal{R}}_h^{p}r)\interleave_{q}^2\leq
		C(h^{2k}|r|_{k}^2+ h^{2k}|\bm{w}|_{k+1}^2).
	\end{aligned}
% \end{subequations}
\]
Furthermore, we have the $L^2$-norm estimate:
\begin{equation*}
	\|\bm{w}-\mathcal{R} _h^{u}\bm{w}\|_{L^2}^2\leq C(h^{2k+2}|r|_{k}^2+ h^{2k+2}|\bm{w}|_{k+1}^2).
\end{equation*}
We also define the Ritz projections of HDG method for the Poisson equation \cite[Theorem 3,Theorem 4]{Oikawa2015}. For any $\bm{C}\in H_{\rm sym}^{1}(\Omega)^{2\times 2}$, there exists a unique symmetric matrix-valued function $(\bm{C}_h,\hat{\bm{C}}_h)\in \bm{W}_h\times\hat{\bm{W}}_h$ such that
\begin{equation*}\label{def:ritz-C}
		\epsilon(\nabla\bm{C},\nabla \bm{D}_h)
		=A_h
		((\bm{C},\hat{\bm{C}}),(\bm{D}_h,\hat{\bm{D}}_h))
		\quad\forall \bm{D}_h\in \bm{W}_h,\hat{\bm{D}}_h\in\hat{\bm{W}}_h.
\end{equation*}
We define $\mathcal{R}_h^C:H_{\rm sym}^{1}(\Omega)^{2\times 2}\rightarrow \bm{W}_h$ with $\mathcal{R}_h^C \bm{C}=\bm{C}_h$
and $\hat{\mathcal{R}}_h^C:H_{\rm sym}^{1}(\Omega)^{2\times 2}\rightarrow \hat{\bm{W}_h}$ with $\hat{\mathcal{R}}_h^C \bm{C}=\hat{\bm{C}}_h$. The Ritz projection satisfies: for $\bm{C}\in H_{\rm sym}^{k+1}(\Omega)^{2\times 2}\cap H_{\rm sym}^{1}(\Omega)^{2\times 2}$
% \begin{subequations}\label{eq:ritz-error-C}
   \begin{equation*}
   	\interleave(\bm{C}-\mathcal{R}_h^C\bm{C},\hat{\mathcal{R}}_h^C\bm{C}-\bm{C})\interleave_{w}^2\leq Ch^{2k}|\bm{C}|_{k+1}^2,\quad
   	\|\bm{C}-\mathcal{R}_h^C\bm{C}\|_{L^2}^2\leq Ch^{2k+2}|\bm{C}|_{k+1}^2.
   \end{equation*}
% \end{subequations}
% With the above preparation, we now proceed error analysis.

% \subsection{Error estimate of the HDG scheme}\label{sec:4-2}
With the above preparation, we now proceed error analysis.
\begin{theorem}\label{th:error} 
	(Error estimate) 
	% Let $P_k$ denote the space of polynomials of degree $k$.
%\begin{equation}\label{eq:con}
%		\begin{cases}
%			\|\bm{u}_h^n\|_{L^\infty}\leq C,\quad n=1,2\cdots,N;&k=1;\\
%			\tau=h^{\lambda+1},\quad\lambda>0;&k>1;
%		\end{cases}
%\end{equation}
     We assume the regularity \eqref{def:error-space} and the $L^\infty$-boundness of the numerical solution:
    \begin{equation}\label{eq:con1}
    	\|\bm{u}_h^n\|_{L^\infty} + \|\bm{C}_h^n\|_{L^\infty}\leq \bar{C}\quad n=0,1,2\cdots,N.
    \end{equation}
    % When $k>1$, we set up
    % \begin{equation}\label{eq:con2}
    % 	\tau=h^{\lambda+1},\quad\lambda>0.
    % \end{equation}
    % We also assume that the given initial value $\|\bm{u}_h^0\|_{L^\infty}\leq C$.
    Supposing that $\alpha,\beta$ are large enough and $\tau\leq\frac{1}{2}$, we have the error estimates
\begin{equation}\label{eq:error-uC}
	\begin{aligned}
		&\|\bm{u}-\bm{u}_h\|_{l^\infty(0,T;L^2)}^2
		+\|{\rm tr}\bm{C}-{\rm tr}\bm{C}_h\|_{l^\infty(0,T;L^2)}^2
		+\nu\|\nabla\bm{u}-\nabla\bm{u}_h\|_{l^2(0,T;L^2)}^2
		\\
		&
		+\epsilon\|\nabla{\rm tr}\bm{C}-\nabla{\rm tr}\bm{C}_h\|_{l^2(0,T;L^2)}^2
		+\frac{\nu}{h}\|(\bm{u}-\bm{u}_h)-(\bm{u}-\hat{\bm{u}}_h)\|_{l^2(0,T;L^2(\mathcal{F}))}^2
		\\
		&
		+\frac{\epsilon }{h}\|({\rm tr}\bm{C}-{\rm tr}\bm{C}_h)-({\rm tr}\bm{C}-{\rm tr}\hat{\bm{C}}_h)\|_{l^2(0,T;L^2(\mathcal{F}))}^2
		\leq
		C(h^{2k}+\tau^2),
	\end{aligned}
\end{equation}
	and
\begin{equation}\label{eq:error-p}
	\begin{aligned}
		\tau\|p-p_h\|_{l^2(0,T;L^2)}^2
		+h\tau\|p-\hat{p}_h\|_{l^2(0,T;L^2(\mathcal{F}))}^2
		\leq
		C(h^{2k}+\tau^2),
	\end{aligned}
\end{equation}
where $C$ depends on the norm of $(\bm u,p,\bm C)$ and $\bar C$, but is independent of $h$, $\tau$, $n$ and $\epsilon^{-1}$.
\end{theorem}
%%%%%%%%%%%%%%%%%%%%%%%%%%%%ss
%%%%%%%%%%%%%%%%%%%%%%%%%%%%
%%%%%%%%%%%%%%%%%%%%%%%%%%%%
\begin{proof}
	For briefness, we set the notations:
	\[% \begin{equation}\label{eq:bri}
		\begin{aligned}
		&\bm{u}^n:=\bm{u}(t^n),\quad
		\bm{e_u}^n:=\bm{u}^n-\bm{u}_h^n,\quad
		\boldsymbol{\hat{e}_u}^n:=\bm{u}^n-\hat{\bm{u}}_h^n,\quad
		p^n:=p(t^n),\quad
	    e_p^n:=p^n-p_h^n,\quad
		\hat{e}_p^n:=p^n-\hat{p}_h^n,
		\\
		&\bm{C}^n:=\bm{C}(t^n),\quad
	    \bm{e_C}^n:=\bm{C}^n-\bm{C}_h^n,\quad
		\boldsymbol{\hat{e}_C}^n:=\bm{C}^n-\hat{\bm{C}}_h^n,
		\end{aligned}
		\]
	% \end{equation}\]
	and we make the decomposition:
	\[
	% \begin{subequations}\label{eq:decomp}
		\begin{aligned}
			&\bm{e_u}^n:=(\bm{u}^n-\mathcal{R} _h^{u}\bm{u}^n)+(\mathcal{R} _h^{u}\bm{u}^n-\bm{u}_h^n)
			=:\boldsymbol{\eta_u}^n+\boldsymbol{\delta_u}^n,
			\quad
			\boldsymbol{\hat{e}_u}^n:=(\bm{u}^n-\hat{\mathcal{R}} _h^{u}\bm{u}^n)+(\hat{\mathcal{R}} _h^{u}\bm{u}^n-\hat{\bm{u}}_h^n)
			=:\boldsymbol{\hat{\eta}_u}^n+\boldsymbol{\hat{\delta}_u}^n,
			\\
			&e_p^n:=(p^n-\mathcal{R} _h^{p}p^n)+(\mathcal{R} _h^{p}p^n-p_h^n):=\eta_p^n+\delta_p^n,
			\quad
			\hat{e}_p^n:=(p^n-\hat{\mathcal{R}} _h^{p}p^n)+(\hat{\mathcal{R}} _h^{p}p^n-\hat{p}_h^n):=\hat{\eta}_p^n+\hat{\delta}_p^n,
			\\
			&\bm{e_C}^n:=(\bm{C}^n-\mathcal{R}_h^C \bm{C}^n)+(\mathcal{R}_h^C \bm{C}^n-\bm{C}_h^n):=\boldsymbol{\eta_C}^n+\boldsymbol{\delta_C}^n,
			\\
			&\boldsymbol{\hat{e}_C}^n:=(\bm{C}^n-\hat{\mathcal{R}}_h^C \bm{C}^n)+(\hat{\mathcal{R}}_h^C \bm{C}^n-\hat{\bm{C}}_h^n):=\boldsymbol{\hat{\eta}_C}^n+\boldsymbol{\hat{\delta}_C}^n.
		\end{aligned}
	% \end{subequations}
	\]
    We set  $\bm{e}:=(\bm{e}^0,\bm{e}^1,\cdots,\bm{e}^N)=\boldsymbol{\delta}+\boldsymbol{\eta}:=(\bm{\delta}^0,\bm{\delta}^1,\cdots,\bm{\delta}^N)+(\bm{\eta}^0,\bm{\eta}^1,\cdots,\bm{\eta}^N)$ for errors of different moment. When $t=t^{n+1}$, the variational problem of \eqref{eq:P} can be written as:
\begin{subequations}
    \begin{align}
        & b_h((q_h,\hat{q}_h),(\bm{u}^{n+1},\bm{u}^{n+1}))=0,\quad\forall q_h\in Q_h,\hat{q}_h\in \hat{Q}_h, \label{eq:O1-a}
        \\
        & \begin{aligned} 
            &(\partial_t \bm{u}^{n+1},\bm{v}_h)+
    		a_h((\bm{u}^{n+1},\bm{u}^{n+1}),(\bm{v}_h,\hat{\bm{v}}_h))
    		+o_h(\bm{u}^{n+1};(\bm{u}^{n+1},\bm{u}^{n+1}),(\bm{v}_h,\hat{\bm{v}}_h))\\
    		&\quad+b_h((p^{n+1},p^{n+1}),(\bm{v}_h,\hat{\bm{v}}_h))
    		+\sum_{K\in\mathcal{T}}\int_{K}{\rm tr}\bm{C}^{n+1}\bm{C}^{n+1}:\nabla\bm{v}_h\,dx\\
    		&\quad+\sum_{K\in\mathcal{T}}\int_{\partial K}{\rm tr}\bm{C}^{n+1}\bm{C}^{n+1}:(\bm{v}_h-\hat{\bm{v}}_h)\otimes\bm{n}\,ds
    		=\int_{\Omega}\bm{f}^{n+1}\cdot\bm{v}_h\,dx,\quad\forall\bm{v}_h\in\bm{V}_h,\hat{\bm{v}}_h\in\hat{\bm{V}}_h,\\
        \end{aligned} \label{eq:O1-b} \\
        &\begin{aligned} 
            &(\partial_t \bm{C}^{n+1},\bm{D}_h)
    		+A_h((\bm{C}^{n+1},\bm{C}^{n+1}),(\bm{D}_h,\hat{\bm{D}}_h))
    		+o_h(\bm{u}^{n+1};(\bm{C}^{n+1},\bm{C}^{n+1}),(\bm{D}_h,\hat{\bm{D}}_h))\\
    		&\quad=\int_{\Omega}((\nabla \bm{u}^{n+1})\bm{C}^{n+1}+\bm{C}^{n+1}(\nabla\bm{u}^{n+1})^\top):\bm{D}_h\,dx
    		+\int_{\Omega}{\rm tr}\bm{C}^{n+1}\bm I:\bm{D}_h\,dx\\
    		&\quad -\int_{\Omega}({\rm tr}\bm{C}^{n+1})^2\bm{C}^{n+1}:\bm{D}_h\,dx
    		,\quad \forall\bm{D}_h\in \bm{W}_h,\hat{\bm{D}}_h\in\hat{\bm{W}}_h.
        \end{aligned} \label{eq:O1-c}
    \end{align}
\end{subequations}

The proof is divided into five steps. In (\textbf{Step 1}), we subtract the original variational equation from the HDG scheme \eqref{eq:Full}. In (\textbf{Step 2}) and (\textbf{Step 3}), we estimate the terms related to $\boldsymbol{\delta}_u^{n+1},\boldsymbol{\hat{\delta}}_u^{n+1}$ and $\boldsymbol{\delta}_C^{n+1},\boldsymbol{\hat{\delta}}_C^{n+1}$ respectively. We apply the discrete Gronwall inequality in (\textbf{Step 4}) to derive the estimates of $\boldsymbol{\delta}_u^{n+1},\boldsymbol{\hat{\delta}}_u^{n+1}$ and $\boldsymbol{\delta}_C^{n+1},\boldsymbol{\hat{\delta}}_C^{n+1}$, furthermore derive the error estimate of $\bm u,\hat{\bm u}$ and $\bm C,\hat{\bm C}$. Finally, we address the error estimate of $p$ by using the inf-sup condition in (\textbf{Step 5}).

%%%%%%%%%%%%%%%%%
%%%%%%%%%%%%%%%%%
%%%%%%%%%%%%%%%%%
(\textbf{Step 1})\label{step1} Subtracting \eqref{eq:Full-a}--\eqref{eq:Full-c} from \eqref{eq:O1-a}--\eqref{eq:O1-c}, we have 
\begin{subequations}
    \begin{align}
        & 0=b_h((q_h,\hat{q}_h),(\bm{u}^{n+1},\bm{u}^{n+1}))
		-b_h((q_h,\hat{q}_h),(\bm{u}_h^{n+1},\hat{\bm{u}}_h^{n+1})),\label{eq:diff-a}
        \\
        & \begin{aligned} 
            &0=\Big(\partial_t\bm{u}^{n+1}-\frac{\bm{u}_h^{n+1}-\bm{u}_h^{n}}{\tau},\bm{v}_h\Big)
			+a_h((\bm{u}^{n+1},\bm{u}^{n+1}),(\bm{v}_h,\hat{\bm{v}}_h))
			-a_h((\bm{u}_h^{n+1},\hat{\bm{u}}_h^{n+1}),(\bm{v}_h,\hat{\bm{v}}_h))
			\\
			&\quad+o_h(\bm{u}^{n+1};(\bm{u}^{n+1},\bm{u}^{n+1}),(\bm{v}_h,\hat{\bm{v}}_h))
			-o_h(\bm{u}_h^{n};(\bm{u}_h^{n+1},\hat{\bm{u}}_h^{n+1}),(\bm{v}_h,\hat{\bm{v}}_h))
			\\
			&\quad+b_h((p^{n+1},p^{n+1}),(\bm{v},\hat{\bm{v}}))
			-b_h((p_h^{n+1},\hat{p}_h^{n+1}),(\bm{v},\hat{\bm{v}}))\\
			&\quad+\sum_{K\in\mathcal{T}}\int_{K}(({\rm tr}\bm{C}^{n+1})\bm{C}^{n+1}-({\rm tr}\bm{C}_h^{n+1})\bm{C}_h^{n}):\nabla\bm{v}_h\, dx\\
			&\quad-\sum_{K\in\mathcal{T}}\int_{\partial K}(({\rm tr}\bm{C}^{n+1})\bm{C}^{n+1}-({\rm tr}\bm{C}_{h}^{n+1})\bm{C}_{h}^{n}):(\bm{v}_h-\hat{\bm{v}}_h)\otimes\bm{n}\,ds,
        \end{aligned} \label{eq:diff-b} \\
        &\begin{aligned} 
            &\Big(\partial_t\bm{C}^{n+1}-\frac{\bm{C}_h^{n+1}-\bm{C}_h^{n}}{\tau},\bm{D}_h\Big)
			+A_h((\bm{C}^{n+1},\bm{C}^{n+1}),(\bm{D}_h,\hat{\bm{D}}_h))
			\\
			&\quad
			-A_h((\bm{C}_h^{n+1},\hat{\bm{C}}_h^{n+1}),(\bm{D}_h,\hat{\bm{D}}_h))
			+o_h(\bm{u}^{n+1};(\bm{C}^{n+1},\bm{C}^{n+1}),(\bm{D}_h,\hat{\bm{D}}_h))
		    \\
			&\quad
            -o_h(\bm{u}_h^{n};(\bm{C}_h^{n+1},\hat{\bm{C}}_h^{n+1}),(\bm{D}_h,\hat{\bm{D}}_h))
			+\int_{\Omega}(({\rm tr}\bm{C}^{n+1})^2\bm{C}^{n+1}-({\rm tr}\bm{C}_h^{n})^2\bm{C}_h^{n+1}):\bm{D}_h\, dx
			\\
			&=\int_{\Omega}(\nabla \bm{u}^{n+1})\bm{C}^{n+1}+\bm{C}^{n+1}(\nabla\bm{u}^{n+1})^\top:\bm{D}_h\, dx
			-\int_{\Omega}(\nabla \bm{u}_h^{n+1})\bm{C}_h^{n}+\bm{C}_h^{n}(\nabla\bm{u}_h^{n+1})^\top:\bm{D}_h\,dx\\
			&\quad+\int_{\Omega}({\rm tr}\bm{C}^{n+1}\bm I-{\rm tr}\bm{C}_h^{n}\bm I):\bm{D}_h\, dx.
        \end{aligned} \label{eq:diff-c}
    \end{align}
\end{subequations}

Then we apply different test functions into \eqref{eq:diff-a}-\eqref{eq:diff-c} in next steps.

	(\textbf{Step 2})\label{step2} Taking $(\bm{v}_h,\hat{\bm{v}}_h,q_h,\hat{q}_h)=(\boldsymbol{\delta_u}^{n+1},\boldsymbol{\hat{\delta}_u}^{n+1},\delta_p^{n+1},\hat{\delta}_p^{n+1})$ into \eqref{eq:diff-a}, \eqref{eq:diff-b} and adding together. Because of the Ritz projection and \eqref{eq:diff-a}, we have
	$$
	\begin{aligned}
		&a_h((\boldsymbol{e_u}^{n+1},\boldsymbol{\hat{e}_u}^{n+1}),(\boldsymbol{\delta_u}^{n+1},\boldsymbol{\hat{\delta}_u}^{n+1}))
		+b_h((e_p^{n+1},\hat{e}_p^{n+1}),(\boldsymbol{\delta_u}^{n+1},\boldsymbol{\hat{\delta}_u}^{n+1}))
		=a_h((\boldsymbol{\delta_u}^{n+1},\boldsymbol{\hat{\delta}_u}^{n+1}),(\boldsymbol{\delta_u}^{n+1},\boldsymbol{\hat{\delta}_u}^{n+1})).
	\end{aligned}
	$$
	Furthermore, the summing equation is 
	\begin{equation*}
		\begin{aligned}
			&\sum_{i=1}^{3}L_i^u:=\Big(\frac{\boldsymbol{\delta_u}^{n+1}-\boldsymbol{\delta_u}^{n}}{\tau},\boldsymbol{\delta_u}^{n+1}\Big)
			+a_h((\boldsymbol{\delta_u}^{n+1},\boldsymbol{\hat{\delta}_u}^{n+1}),(\boldsymbol{\delta_u}^{n+1},\boldsymbol{\hat{\delta}_u}^{n+1}))\\
			&\quad+\bigg(o_h(\bm{u}^{n+1};(\bm{u}^{n+1},\bm{u}^{n+1}),(\boldsymbol{\delta_u}^{n+1},\boldsymbol{\hat{\delta}_u}^{n+1}))
			-o_h(\bm{u}_h^{n};(\bm{u}_h^{n+1},\hat{\bm{u}}_h^{n+1}),(\boldsymbol{\delta_u}^{n+1},\boldsymbol{\hat{\delta}_u}^{n+1}))\bigg)\\
            &=\Big(\frac{\mathcal{R}_{h}^{u}\bm{u}^{n+1}-\mathcal{R}_{h}^{u}\bm{u}^{n}}{\tau}-\partial_t\bm{u}^{n+1},\boldsymbol{\delta_u}^{n+1}\Big)
			+\sum_{K\in\mathcal{T}}\int_{K}(({\rm tr}\bm{C}^{n+1})\bm{C}^{n+1}-({\rm tr}\bm{C}_h^{n+1})\bm{C}_h^{n}):\nabla\boldsymbol{\delta_u}^{n+1}\, dx\\
			&\quad+\sum_{K\in\mathcal{T}}\int_{\partial K}(({\rm tr}\bm{C}^{n+1})\bm{C}^{n+1}-({\rm tr}\bm{C}_{h}^{n+1})\bm{C}_{h}^{n}):(\boldsymbol{\delta_u}^{n+1}-\boldsymbol{\hat{\delta}_u}^{n+1})\otimes\bm{n}\,d s
			=:\sum_{i=1}^{3}R_i^u.
		\end{aligned}
	\end{equation*}
%	\begin{equation}
%		\begin{aligned}
%			&\underbrace{(\frac{\boldsymbol{\delta_u}^{n+1}-\boldsymbol{\delta_u}^{n}}{\tau},\boldsymbol{\delta_u}^{n+1})}_{L_1^u}
%			+\underbrace{a_h((\bm{u}^{n+1},\bm{u}^{n+1}),(\boldsymbol{\delta_u}^{n+1},\boldsymbol{\hat{\delta}_u}^{n+1}))
%			-a_h((\bm{u}_h^{n+1},\hat{\bm{u}}_h^{n+1}),(\boldsymbol{\delta_u}^{n+1},\boldsymbol{\hat{\delta}_u}^{n+1}))}_{L_2^u}
%			\\
%			&+\underbrace{o_h(\bm{u}^{n+1};(\bm{u}^{n+1},\bm{u}^{n+1}),(\boldsymbol{\delta_u}^{n+1},\boldsymbol{\hat{\delta}_u}^{n+1}))
%			-o_h(\bm{u}_h^{n};(\bm{u}_h^{n+1},\hat{\bm{u}}_h^{n+1}),(\boldsymbol{\delta_u}^{n+1},\boldsymbol{\hat{\delta}_u}^{n+1}))}_{L_3^u}
%			\\
%			&=\underbrace{(\frac{\mathcal{R}_{h}^{u}\bm{u}^{n+1}-\mathcal{R}_{h}^{u}\bm{u}^{n}}{\tau}-\partial_t\bm{u}^{n+1},\boldsymbol{\delta_u}^{n+1})}_{R_1^u}
%			-\underbrace{\sum_{K\in\mathcal{T}}\int_{K}(({\rm tr}\bm{C}^{n+1})\bm{C}^{n+1}-({\rm tr}\bm{C}_h^{n+1})\bm{C}_h^{n}):\nabla\boldsymbol{\delta_u}^{n+1}\, dx}_{R_2^u}\\
%			&+\underbrace{\sum_{K\in\mathcal{T}}\int_{\partial K}(({\rm tr}\bm{C}^{n+1})\bm{C}^{n+1}-({\rm tr}\bm{C}_{h}^{n+1})\bm{C}_{h}^{n}):(\boldsymbol{\delta_u}^{n+1}-\boldsymbol{\hat{\delta}_u}^{n+1})\otimes\bm{n}\,d s,}_{R_3^u}
%		\end{aligned}
%	\end{equation}
Firstly we see that
	\begin{equation}\label{eq:L1u}
		L_1^u
		=\frac{1}{2\tau}(\|\boldsymbol{\delta_u}^{n+1}\|_{L^2}^2-\|\boldsymbol{\delta_u}^{n}\|_{L^2}^2+\|\boldsymbol{\delta_u}^{n+1}-\boldsymbol{\delta_u}^{n}\|_{L^2}^2).
	\end{equation}
For $L_2^u$, we have
	\begin{equation*}\label{eq:L2u}
		\begin{aligned}
			L_2^u&=\nu\interleave(\boldsymbol{\delta_u}^{n+1},\boldsymbol{\hat{\delta}_u}^{n+1})\interleave_{v}^2
			-2\nu\sum_{K\in\mathcal{T}}\int_{\partial K}\nabla\boldsymbol{\delta_u}^{n+1}:(\boldsymbol{\delta_u}^{n+1}-\boldsymbol{\hat{\delta}_u}^{n+1})\otimes\bm{n}\,ds\\
			&\geq
			\nu(1-\gamma)\|\nabla\boldsymbol{\delta_u}^{n+1}\|_{L^2}^2
			+\sum_{K\in\mathcal{T}}\frac{\nu}{h_K}(\alpha-\frac{c_I}{\gamma})\|\boldsymbol{\delta_u}^{n+1}-\boldsymbol{\hat{\delta}_u}^{n+1}\|_{L^2(\partial K)}^2.
		\end{aligned}
	\end{equation*}
%%%%%%%%
%%%%%%%%
    We decompose $L_3^u$ as follows:
    $$
    	\begin{aligned}
    		L_3^u&=o_h(\bm{u}^{n+1}-\bm{u}^{n};(\bm{u}^{n+1},\bm{u}^{n+1}),(\boldsymbol{\delta_u}^{n+1},\boldsymbol{\hat{\delta}_u}^{n+1}))\\
    		&\quad+o_h(\boldsymbol{\eta_u}^{n};(\bm{u}^{n+1},\bm{u}^{n+1}),(\boldsymbol{\delta_u}^{n+1},\boldsymbol{\hat{\delta}_u}^{n+1}))
    		+o_h(\bm{u}_h^{n};(\boldsymbol{\eta_u}^{n+1},\boldsymbol{\hat{\eta}_u}^{n+1}),(\boldsymbol{\delta_u}^{n+1},\boldsymbol{\hat{\delta}_u}^{n+1}))\\
    		&\quad+o_h(\boldsymbol{\delta_u}^{n};(\bm{u}^{n+1},\bm{u}^{n+1}),(\boldsymbol{\delta_u}^{n+1},\boldsymbol{\hat{\delta}_u}^{n+1}))
    		+o_h(\bm{u}_h^{n};(\boldsymbol{\delta_u}^{n+1},\boldsymbol{\hat{\delta}_u}^{n+1}),(\boldsymbol{\delta_u}^{n+1},\boldsymbol{\hat{\delta}_u}^{n+1}))=:\sum_{i=1}^{5}L_{3i}^u.
    	\end{aligned}
    $$
    % The boundary of $\{L_{3i}^u\}_{i=1}^5$ can be estimated by the Green's formula and the interpolation inequality. 
    % Because the third integral term is equal to $0$, and by using mass conservation \eqref{eq:mass-full}, we have 
    Next we estimate $L_{3i}^u$ ($i = 1,\ldots,5$).
    In view of 
    $$
    \begin{aligned}
    &\sum_{K\in\mathcal{T}}\int_{\partial K}\frac{|(\bm u^{n+1}-\bm u^{n})\cdot\bm{n}|}{2}(\bm{u}^{n+1}-\bm{u}^{n+1})(\boldsymbol{\delta}_u^{n+1}-\boldsymbol{\hat{\delta}_u}^{n+1})\,ds=0,\\
    &\sum_{K\in\mathcal{T}}\int_{\partial K}(\bm{u}^{n+1}-\bm{u}^{n})\cdot\bm{n}(\bm{u}^{n+1}\cdot\boldsymbol{\hat{\delta}_u}^{n+1})\,ds = 0,
    \end{aligned}
    $$
we have

    $$
    \begin{aligned}
    	L_{31}^u&=
    	-\sum_{K\in\mathcal{T}}\int_{K}\bm{u}^{n+1}\otimes(\bm{u}^{n+1}-\bm{u}^{n}):\nabla\boldsymbol{\delta_u}^{n+1}\,dx
    	+\sum_{K\in\mathcal{T}}\int_{\partial K}(\bm{u}^{n+1}-\bm{u}^{n})\cdot\bm{n}(\bm{u}^{n+1}\cdot\boldsymbol{\delta_u}^{n+1})\,ds\\
    	&
    	=
    	\sum_{K\in\mathcal{T}}\int_{K}\nabla\cdot(\bm{u}^{n+1}-\bm{u}^{n})\bm{u}^{n+1}\cdot\boldsymbol{\delta_u}^{n+1}\,dx
    	+\sum_{K\in\mathcal{T}}\int_{K}(\bm{u}^{n+1}-\bm{u}^{n})^{\top}\nabla\bm{u}^{n+1}\boldsymbol{\delta_u}^{n+1}\,dx
    	\\
    	& =	\sum_{K\in\mathcal{T}}\int_{K}(\bm{u}^{n+1}-\bm{u}^{n})^{\top}\nabla\bm{u}^{n+1}\boldsymbol{\delta_u}^{n+1}\,dx \quad(\text{by } \nabla\cdot\bm u= 0),
    \end{aligned}
    $$
    where we have used the integration by parts in second equality.
    Therefore,
    \begin{equation*}
    	\begin{aligned}
    		|L_{31}^u|&
    		\leq
    		\|\nabla\bm{u}^{n+1}\|_{L^\infty}\|\bm{u}^{n+1}-\bm{u}^{n}\|_{L^2}\|\boldsymbol{\delta_u}^{n+1}\|_{L^2}
    		% &\leq c_0\|\bm{u}^{n+1}\|_{H^2}^{\frac{1}{2}}\|\bm{u}^{n+1}\|_{H^1}^{\frac{1}{2}}
    		% \tau^{\frac{1}{2}}\|\bm{u}\|_{H^1(t^n,t^{n+1};L^2)}^{\frac{1}{2}}\|\bm{u}\|_{H^1(t^n,t^{n+1};H^1)}^{\frac{1}{2}}\|\boldsymbol{\delta_u}^{n+1}\|_{L^2}\\
    		\leq C
    		\|\bm{u}^{n+1}\|_{W^{1,\infty}}
    		\tau^{\frac{1}{2}}\|\bm{u}\|_{H^1(t^n,t^{n+1};L^2)}\|\boldsymbol{\delta_u}^{n+1}\|_{L^2}
    		\\&
    		\leq C\tau\gamma^{-1}\|\bm{u}\|_{H^1(t^n,t^{n+1};L^2)}^2
    		+\gamma\|\boldsymbol{\delta_u}^{n+1}\|_{L^2}^2,
    	\end{aligned}
    \end{equation*}
    where $C$ depends on $\|\bm{u}^{n+1}\|_{H^2}$. We can estimate $\{L_{3i}^u\}_{i=2}^5$ in the same manner:
    \begin{subequations}\label{eq:L3iu}
    	  \begin{align}
    	  &\begin{aligned}
    		|L_{32}^u| &= 
    		\sum_{K\in\mathcal{T}}\int_{K}\nabla\cdot\boldsymbol{\eta_u}^{n}\bm{u}^{n+1}\cdot\boldsymbol{\delta_u}^{n+1}\,dx
    		+\sum_{K\in\mathcal{T}}\int_{K}(\boldsymbol{\eta_u}^{n})^{\top}\nabla\bm{u}^{n+1}\boldsymbol{\delta_u}^{n+1}\,dx
    		\nonumber
    		\\ &
    		\leq
    		C\gamma^{-1}\|\boldsymbol{\eta_u}^{n}\|_{L^2}^2
    		+\gamma\|\nabla\bm u^{n+1}\|_{L^\infty}^2\|\boldsymbol{\delta_u}^{n+1}\|_{L^2}^2 \quad(\text{by } \nabla\cdot\boldsymbol{\eta_u}^{n}=0),
    		\end{aligned}
    		\nonumber
    		\\
    		& \begin{aligned}
    		|L_{33}^u|
    		&\leq\left|
    		\sum_{K\in\mathcal{T}}\int_{K}\bm{u}_h^{n\top}\nabla\boldsymbol{\eta_u}^{n+1}\boldsymbol{\delta_u}^{n+1}\,dx
    		-\sum_{K\in\mathcal{T}}\int_{\partial K}\bm{u}_h^{n}\cdot\bm{n}(\boldsymbol{\eta_u}^{n+1}-\boldsymbol{\hat{\eta}_u}^{n+1}):\boldsymbol{\delta_u}^{n+1}\,ds\right.
    		\nonumber
    		\\
    		&\quad\left.+\sum_{K\in\mathcal{T}}\int_{\partial K}\frac{\bm{u}_h^{n}\cdot\bm{n}+|\bm{u}_h^{n}\cdot\bm{n}|}{2}(\boldsymbol{\eta_u}^{n+1}-\boldsymbol{\hat{\eta}_u}^{n+1})(\boldsymbol{\delta_u}^{n+1}-\boldsymbol{\hat{\delta}_u}^{n+1})\,ds
    		\right|
    		\nonumber
    		\\
    		&\leq C\gamma^{-1}\|\bm{u}_h^{n}\|_{L^\infty}^2\interleave(\boldsymbol{\eta_u}^{n+1},\boldsymbol{\hat{\eta}_u}^{n+1})\interleave_{v}^2
    		+\gamma\|\boldsymbol{\delta_u}^{n+1}\|_{L^2}^2
    		\nonumber
    		\\
    		&\quad+C\gamma^{-1}\|\bm{u}_h^{n}\|_{L^\infty}^2\sum_{K\in\mathcal{T}}\|\boldsymbol{\eta_u}^{n+1}-\boldsymbol{\hat{\eta}_u}^{n+1}\|_{L^2(\partial K)}^2
    		+\gamma\sum_{K\in\mathcal{T}}\int_{\partial K}\frac{|\bm{u}_h^{n}\cdot\bm{n}|}{2}|\boldsymbol{\delta_u}^{n+1}-\boldsymbol{\hat{\delta}_u}^{n+1}|^2\, ds,
    		\end{aligned}
    		\nonumber
    		\\
    		&
    		\begin{aligned}
    		|L_{34}^u| &= 
    		\sum_{K\in\mathcal{T}}\int_{K}\nabla\cdot\boldsymbol{\delta_u}^{n}\bm{u}^{n+1}\cdot\boldsymbol{\delta_u}^{n+1}\,dx
    		+\sum_{K\in\mathcal{T}}\int_{K}(\boldsymbol{\delta_u}^{n})^{\top}\nabla\bm{u}^{n+1}\boldsymbol{\delta_u}^{n+1}\,dx
    		\nonumber
    		\\ &
    		\leq
    		C\|\nabla\bm u^{n+1}\|_{L^\infty}(\|\boldsymbol{\delta_u}^{n}\|_{L^2}^2 + \|\boldsymbol{\delta_u}^{n+1}\|_{L^2}^2)
    		\quad(\text{by } \nabla\cdot\boldsymbol{\delta_u}^{n}=0),
    		\end{aligned}
    		\nonumber
    		\\
    		&L_{35}^u
    		=\sum_{K\in\mathcal{T}}\int_{\partial K}\frac{|\bm{u}_h^{n}\cdot\bm{n}|}{2}|\boldsymbol{\delta_u}^{n+1}-\boldsymbol{\hat{\delta}_u}^{n+1}|^2\, ds.
    		\nonumber
    	\end{align}
\end{subequations}
    So we have
    \begin{equation}\label{eq:L3u}
    	\begin{aligned}
    		L_3^u
    		&\geq
    		(1-\gamma)\sum_{K\in\mathcal{T}}\int_{\partial K}\frac{|\bm{u}_h^{n}\cdot\bm{n}|}{2}|\boldsymbol{\delta_u}^{n+1}-\boldsymbol{\hat{\delta}_u}^{n+1}|^2\, ds
    		-C(\|\boldsymbol{\delta_u}^{n}\|_{L^2}^2 + \|\boldsymbol{\delta_u}^{n+1}\|_{L^2}^2)
    		\\
    		& \quad
    		-C\gamma^{-1}\big(
    		\tau\|\bm{u}\|_{H^1(t^n,t^{n+1};L^2)}^2
    		+\|\boldsymbol{\eta_u}^{n}\|_{L^2}^2
    		+\interleave(\boldsymbol{\eta_u}^{n+1},\boldsymbol{\hat{\eta}_u}^{n+1})\interleave_{v}^2
    		+\sum_{K\in\mathcal{T}}\|\boldsymbol{\eta_u}^{n+1}-\boldsymbol{\hat{\eta}_u}^{n+1}\|_{L^2(\partial K)}^2
    		\big)	
    		.
    	\end{aligned}
    \end{equation}
%%%%%%%%%%%
%%%%%%%%%%%

Now we deal with $R_i^u$ $(i=1,2,3)$. Firstly, we estimate $R_1^u$:
	\begin{equation*}\label{eq:R1u}
		\begin{aligned}
			|R_1^u|
			&\leq C\left(h^{k}\|\bm{u}\|_{C^1(0,T;H^{k+1})}+\tau^{\frac{1}{2}}\|\bm{u}_{tt}\|_{L^2(t^n,t^{n+1};L^2)}\right)\|\boldsymbol{\delta_u}^{n+1}\|_{L^2}.
			% \\
			% &\leq C(h^{2k}\|\bm{u}\|_{C^1(0,T,H^{k+1})}^2+\tau\|\bm{u}\|_{H^2(t^n,t^{n+1};L^2)}^2)
			% +c\|\boldsymbol{\delta_u}^{n+1}\|_{L^2}^2.
		\end{aligned}
	\end{equation*}
	For $R_2^u$, we have the decomposition
	$$
	\begin{aligned}
		R_2^u
		&=\sum_{K\in\mathcal{T}}\int_{K}({\rm tr}\bm{C}^{n+1})(\bm{C}^{n+1}-\bm{C}^{n}):\nabla\boldsymbol{\delta_u}^{n+1}\, dx\\
		&\quad+\sum_{K\in\mathcal{T}}\int_{K}({\rm tr}\bm{C}^{n+1})\boldsymbol{\eta_C}^{n}:\nabla\boldsymbol{\delta_u}^{n+1}\, dx
		+\sum_{K\in\mathcal{T}}\int_{K}({\rm tr}\bm{C}^{n+1})\boldsymbol{\delta_C}^{n}:\nabla\boldsymbol{\delta_u}^{n+1}\, dx\\
		&\quad+\sum_{K\in\mathcal{T}}\int_{K}({\rm tr}\boldsymbol{\eta_C}^{n+1})\bm{C}_h^{n}:\nabla\boldsymbol{\delta_u}^{n+1}\, dx
   	    +\sum_{K\in\mathcal{T}}\int_{K}({\rm tr}\boldsymbol{\delta_C}^{n+1})\bm{C}_h^{n}:\nabla\boldsymbol{\delta_u}^{n+1}\, dx=:\sum_{i=1}^{5}R_{2i}^u.
	\end{aligned}
	$$
	Now we estimate $\{R_{2i}^u\}_{i=1}^{5}$.
	$$
	\begin{aligned}
			|R_{21}^u|&\leq
				\nu\gamma\|\nabla\boldsymbol{\delta_u}^{n+1}\|_{L^2}^2
				+C\tau\gamma^{-1}\nu^{-1}\|\bm{C}^{n+1}\|_{L^\infty}^2\|\bm{C}_t\|_{L^2(t^n,t^{n+1};L^2)}^2.
		% 		\leq
		% \nu\gamma_3\|\nabla\boldsymbol{\delta_u}^{n+1}\|_{L^2}^2+\frac{c_2 \tau}{4\nu\gamma_3}\|\bm{C}_t\|_{L^2(t^n,t^{n+1};L^2)}^2,
	\end{aligned}
	$$
	Analogously, we have
	$$
	\begin{aligned}
%		|R_{21}^u|&\leq
%		\nu\gamma_3\|\nabla\boldsymbol{\delta_u}^{n+1}\|_{L^2}^2
%	    +\frac{ \tau}{4\nu\gamma_3}\|\bm{C}^{n+1}\|_{L^\infty}\|\bm{C}_t\|_{L^2(t^n,t^{n+1};L^2)}^2\\
%	&\leq
%		\nu\gamma_3\|\nabla\boldsymbol{\delta_u}^{n+1}\|_{L^2}^2+\frac{c_2 \tau}{4\nu\gamma_3}\|\bm{C}_t\|_{L^2(t^n,t^{n+1};L^2)}^2,\\
		|R_{22}^u|&\leq
		\nu\gamma\|\nabla\boldsymbol{\delta_u}^{n+1}\|_{L^2}^2
		+C\gamma^{-1}\nu^{-1}\|\boldsymbol{\eta_C}^{n}\|_{L^2}^2,\quad
		|R_{23}^u| \leq
		\nu\gamma\|\nabla\boldsymbol{\delta_u}^{n+1}\|_{L^2}^2
		+C\gamma^{-1}\nu^{-1}\|\boldsymbol{\delta_C}^{n}\|_{L^2}^2\\
		|R_{24}^u|&\leq
		\nu\gamma\|\nabla\boldsymbol{\delta_u}^{n+1}\|_{L^2}^2
		+C\gamma^{-1}\nu^{-1}\|\bm{C}_h^n\|_{L^\infty}^2\|\boldsymbol{\eta_C}^{n+1}\|_{L^2}^2,\,
		|R_{25}^u|\leq
		\nu\gamma\|\nabla\boldsymbol{\delta_u}^{n+1}\|_{L^2}^2
		+C\gamma^{-1}\nu^{-1}\|\bm{C}_h^n\|_{L^\infty}^2\|\boldsymbol{\delta_C}^{n+1}\|_{L^2}^2.
	\end{aligned}
	$$
	Summing up the above estimates of $R_{2i}^u$, we arrive
	\begin{equation*}\label{eq:R2u}
		\begin{aligned}
			R_2^u
			&\leq
			5\nu\gamma\|\nabla\boldsymbol{\delta_u}^{n+1}\|_{L^2}^2
			+C\gamma^{-1}\nu^{-1}
			\bigg(
			\tau\|\bm{C}_h^{n+1}\|_{L^\infty}^2\|\bm{C}_t\|_{L^2(t^n,t^{n+1};L^2)}^2
			+\|\boldsymbol{\eta_C}^{n}\|_{L^2}^2
			+\|\boldsymbol{\delta_C}^{n}\|_{L^2}^2\\
			&\quad+\|\bm{C}_h^n\|_{L^\infty}^2\|\boldsymbol{\eta_C}^{n+1}\|_{L^2}^2
			+\|\bm{C}_h^n\|_{L^\infty}^2\|\boldsymbol{\delta_C}^{n+1}\|_{L^2}^2
			\bigg).
		\end{aligned}
	\end{equation*}
%%%%%%%%%%%%
%%%%%%%%%%%%
    Similarly, we divide $R_3^u$ into five terms:
    $$
    	\begin{aligned}
    		R_3^u
    		&=\sum_{K\in\mathcal{T}}\int_{\partial K}({\rm tr}\bm{C}^{n+1})(\bm{C}^{n+1}-\bm{C}^{n}):(\boldsymbol{\delta_u}^{n+1}-\boldsymbol{\hat{\delta}_u}^{n+1})\otimes\bm{n}\,d s
    		\\
    		&\quad+\sum_{K\in\mathcal{T}}\int_{\partial K}({\rm tr}\bm{C}^{n+1})\boldsymbol{\eta_C}^{n}:(\boldsymbol{\delta_u}^{n+1}-\boldsymbol{\hat{\delta}_u}^{n+1})\otimes\bm{n}\,d s
    		+\sum_{K\in\mathcal{T}}\int_{\partial K}({\rm tr}\bm{C}^{n+1})\boldsymbol{\delta_C}^{n}:(\boldsymbol{\delta_u}^{n+1}-\boldsymbol{\hat{\delta}_u}^{n+1})\otimes\bm{n}\,d s\\
    		&\quad+\sum_{K\in\mathcal{T}}\int_{\partial K}({\rm tr}\boldsymbol{\eta_C}^{n+1})\bm{C}_h^{n}:(\boldsymbol{\delta_u}^{n+1}-\boldsymbol{\hat{\delta}_u}^{n+1})\otimes\bm{n}\,d s
    		+\sum_{K\in\mathcal{T}}\int_{\partial K}({\rm tr}\boldsymbol{\delta_C}^{n+1})\bm{C}_h^{n}:(\boldsymbol{\delta_u}^{n+1}-\boldsymbol{\hat{\delta}_u}^{n+1})\otimes\bm{n}\,d s=:\sum_{i=1}^{5}R_{3i}^u.
    	\end{aligned}
    $$
The estimates of $\{R_{3i}^u\}_{i=1}^{5}$ are given as follows:
    $$
    \begin{aligned}
    	&
    	\begin{aligned}
    	|R_{31}^u|&\leq
    	\nu\gamma\sum_{K\in\mathcal{T}}h_K^{-1}\|\boldsymbol{\delta_u}^{n+1}-\boldsymbol{\hat{\delta}_u}^{n+1}\|_{L^2(\partial K)}^2
    	+C\tau\nu^{-1}\gamma^{-1}\|\text{tr}\bm C^{n+1}\|_{L^\infty}^2
    	\sum_{K\in\mathcal{T}}h_K
    	\|\bm{C}\|_{H^1(t^n,t^{n+1};L^2(\partial K))}^2
    	\\ &\leq
    	\nu\gamma\sum_{K\in\mathcal{T}}h_K^{-1}\|\boldsymbol{\delta_u}^{n+1}-\boldsymbol{\hat{\delta}_u}^{n+1}\|_{L^2(\partial K)}^2
    	\\ &
    	\quad+
    	C\tau\nu^{-1}\gamma^{-1}\|\text{tr}\bm C^{n+1}\|_{L^\infty}^2
    	(\|\bm{C}\|_{H^1(t^n,t^{n+1};L^2)}^2 
    	+ h^2\|\nabla\bm{C}\|_{H^1(t^n,t^{n+1};L^2)}^2) \quad
    	(\text{by } \eqref{eq:T}), 
    	\end{aligned}\\
    	&
    	|R_{32}^u|\leq
    	\nu\gamma\sum_{K\in\mathcal{T}}h_K^{-1}\|\boldsymbol{\delta_u}^{n+1}-\boldsymbol{\hat{\delta}_u}^{n+1}\|_{L^2(\partial K)}^2
    	+C\nu^{-1}\gamma^{-1}
    	(\|\boldsymbol{\eta_C}^{n}\|_{L^2}^2 
    	+h^2\|\nabla\boldsymbol{\eta_C}^{n}\|_{L^2}^2 ),\\
    	&|R_{33}^u|\leq
    	\nu\gamma\sum_{K\in\mathcal{T}}h_K^{-1}\|\boldsymbol{\delta_u}^{n+1}-\boldsymbol{\hat{\delta}_u}^{n+1}\|_{L^2(\partial K)}^2
    	+C\nu^{-1}\gamma^{-1}\|\boldsymbol{\delta_C}^{n}\|_{L^2}^2 
    	\quad (\text{by }\eqref{eq:IT}),\\
    	&|R_{34}^u|\leq
    	\nu\gamma\sum_{K\in\mathcal{T}}h_K^{-1}\|\boldsymbol{\delta_u}^{n+1}-\boldsymbol{\hat{\delta}_u}^{n+1}\|_{L^2(\partial K)}^2
    	+C\nu^{-1}\gamma^{-1}
    	(\|\boldsymbol{\eta_C}^{n+1}\|_{L^2}^2 
    	+h^2\|\nabla\boldsymbol{\eta_C}^{n+1}\|_{L^2}^2 )
    	\quad (\text{by }\eqref{eq:T}),\\
    	&|R_{35}^u|\leq
    	\nu\gamma\sum_{K\in\mathcal{T}}h_K^{-1}\|\boldsymbol{\delta_u}^{n+1}-\boldsymbol{\hat{\delta}_u}^{n+1}\|_{L^2(\partial K)}^2
    	+C\nu^{-1}\gamma^{-1}\|\bm{C}_h^n\|_{L^\infty}^2\|\boldsymbol{\delta_C}^{n+1}\|_{L^2}^2
    	\quad (\text{by }\eqref{eq:IT}).\\
    \end{aligned}
    $$
    Hence, 
    \begin{equation*}\label{eq:R3u}
    	\begin{aligned}
    		|R_3^u|
    		&\leq
    		5\nu\gamma\sum_{K\in\mathcal{T}}h_K^{-1}\|\boldsymbol{\delta_u}^{n+1}-\boldsymbol{\hat{\delta}_u}^{n+1}\|_{L^2(\partial K)}^2
    		+C\nu^{-1}\gamma^{-1}
    		\bigg(\tau\|\text{tr}\bm C^{n+1}\|_{L^\infty}^2
    		(\|\bm{C}\|_{H^1(t^n,t^{n+1};L^2)}^2
    		+ h^2\|\nabla\bm{C}\|_{H^1(t^n,t^{n+1};L^2)}^2)\\
    		&\quad
    		+\|\boldsymbol{\eta_C}^{n}\|_{L^2}^2
    		+h^2\|\nabla\boldsymbol{\eta_C}^{n}\|_{L^2}^2
    		+\|\boldsymbol{\delta_C}^{n}\|_{L^2}^2
    		+\|\boldsymbol{\eta_C}^{n+1}\|_{L^2}^2
    		+h^2\|\nabla\boldsymbol{\eta_C}^{n+1}\|_{L^2}^2
    		+\|\bm{C}_h^n\|_{L^\infty}^2\|\boldsymbol{\delta_C}^{n+1}\|_{L^2}^2\bigg).
    	\end{aligned}
    \end{equation*}
    Above all we have the estimates of terms relative to $\boldsymbol{\delta}_u^{n+1},\boldsymbol{\hat{\delta}}_u^{n+1}$. In \textbf{Step 4}, we will sum up all the terms $\sum_{i=1}^3L_i^u,\sum_{i=0}^3R_i^u$ together with the estimate of $\boldsymbol{\delta}_C^{n+1},\boldsymbol{\hat{\delta}}_C^{n+1}$, which comes from the result of \textbf{Step 3}.

	(\textbf{Step 3})\label{step3} We estimate $\boldsymbol{\delta_C}^{n+1},\boldsymbol{\hat{\delta}_C}^{n+1}$ in this step. Taking $(\bm{D}_h,\hat{\bm{D}}_h)=(\frac{1}{2}\boldsymbol{\delta_C}^{n+1},\frac{1}{2}\boldsymbol{\hat{\delta}_C}^{n+1})$ into \eqref{eq:diff-c}, we have
	\begin{equation*}
		\begin{aligned}
			\sum_{i=1}^{4}L_i^C&:=\frac{1}{2}(\frac{\boldsymbol{\delta_C}^{n+1}-\boldsymbol{\delta_C}^{n}}{\tau},\boldsymbol{\delta_C}^{n+1})\\
			&\quad+\frac{1}{2}\bigg(A_h((\bm{C}^{n+1},\bm{C}^{n+1}),(\boldsymbol{\delta_C}^{n+1},\boldsymbol{\hat{\delta}_C}^{n+1}))
			-A_h((\bm{C}_h^{n+1},\hat{\bm{C}}_h^{n+1}),(\boldsymbol{\delta_C}^{n+1},\boldsymbol{\hat{\delta}_C}^{n+1}))\bigg)\\
			&\quad+\frac{1}{2}\bigg(o_h(\bm{u}^{n+1};(\bm{C}^{n+1},\bm{C}^{n+1}),(\boldsymbol{\delta_C}^{n+1},\boldsymbol{\hat{\delta}_C}^{n+1}))
			-o_h(\bm{u}_h^{n};(\bm{C}_h^{n+1},\hat{\bm{C}}_h^{n+1}),(\boldsymbol{\delta_C}^{n+1},\boldsymbol{\hat{\delta}_C}^{n+1}))\bigg)\\
			&\quad+\frac{1}{2}\int_{\Omega}(({\rm tr}\bm{C}^{n+1})^2\bm{C}^{n+1}-({\rm tr}\bm{C}_h^{n})^2\bm{C}_h^{n+1}):\boldsymbol{\delta_C}^{n+1}\, dx\\
            &=\frac{1}{2}\int_{\Omega}((\nabla \bm{u}^{n+1})\bm{C}^{n+1}+\bm{C}^{n+1}(\nabla\bm{u}^{n+1})^\top
			-(\nabla \bm{u}_h^{n+1})\bm{C}_h^{n}+\bm{C}_h^{n}(\nabla\bm{u}_h^{n+1})^\top):\boldsymbol{\delta_C}^{n+1}\,dx\\
			&\quad+\frac{1}{2}(\frac{\mathcal{R}_h^C\bm{C}^{n+1}-\mathcal{R}_h^C \bm{C}^{n}}{\tau}-\partial_t\bm{C}^{n+1},\boldsymbol{\delta_C}^{n+1})
			+\frac{1}{2}\int_{\Omega}(\bm{C}^{n+1}-\bm{C}_h^{n})\bm I:\boldsymbol{\delta_C}^{n+1}\, dx
			:=\sum_{i=1}^{3}R_i^C,
		\end{aligned}
	\end{equation*}
 The calulation of $L^C_1$, $L^C_2$ and $L^C_3$ is similar to \textbf{Step 2}.
 \begin{equation}
 L_1^C=
 \frac{1}{4\tau}(\|\boldsymbol{\delta_C}^{n+1}\|_{L^2}^2-\|\boldsymbol{\delta_C}^{n}\|_{L^2}^2+\|\boldsymbol{\delta_C}^{n+1}-\boldsymbol{\delta_C}^{n}\|_{L^2}^2),
 \end{equation}
 \begin{equation}
L_2^C\geq
\frac{\ep}{2}(1-\gamma)\|\nabla\boldsymbol{\delta_C}^{n+1}\|_{L^2}^2
+\sum_{K\in\mathcal{T}}\frac{\epsilon}{2h_K}(\beta-\frac{c_I}{\gamma})\|\boldsymbol{\delta_C}^{n+1}-\boldsymbol{\hat{\delta}_C}^{n+1}\|_{L^2(\partial K)}^2.
 \end{equation}
	% \begin{subequations}
    % \begin{align}
	% 	&L_1^C=
	% 	\frac{1}{4\tau}(\|\boldsymbol{\delta_C}^{n+1}\|_{L^2}^2-\|\boldsymbol{\delta_C}^{n}\|_{L^2}^2+\|\boldsymbol{\delta_C}^{n+1}-\boldsymbol{\delta_C}^{n}\|_{L^2}^2),
	% 	\\
	% 	&L_2^C
	% 		\geq\frac{\ep}{2}(1-\gamma)\|\nabla\boldsymbol{\delta_C}^{n+1}\|_{L^2}^2
	% 		+\sum_{K\in\mathcal{T}}\frac{\epsilon}{2h_K}(\beta-\frac{c_I}{\gamma})\|\boldsymbol{\delta_C}^{n+1}-\boldsymbol{\hat{\delta}_C}^{n+1}\|_{L^2(\partial K)}^2.\label{eq:L2C}
	%   \end{align}
	% \end{subequations}
	$L_3^C$ can be written as
	$$
	\begin{aligned}
		L_3^C&=\frac{1}{2}o_h(\bm{u}^{n+1}-\bm{u}^{n};(\bm{C}^{n+1},\bm{C}^{n+1}),(\boldsymbol{\delta_C}^{n+1},\boldsymbol{\hat{\delta}_C}^{n+1}))\\
		&\quad+\frac{1}{2}o_h(\boldsymbol{\eta_u}^{n};(\bm{C}^{n+1},\bm{C}^{n+1}),(\boldsymbol{\delta_C}^{n+1},\boldsymbol{\hat{\delta}_C}^{n+1}))
		+\frac{1}{2}o_h(\bm{u}_h^{n};(\boldsymbol{\eta_C}^{n+1},\boldsymbol{\hat{\eta}_C}^{n+1}),(\boldsymbol{\delta_C}^{n+1},\boldsymbol{\hat{\delta}_C}^{n+1}))\\
		&\quad+\frac{1}{2}o_h(\boldsymbol{\delta_u}^{n};(\bm{C}^{n+1},\bm{C}^{n+1}),(\boldsymbol{\delta_C}^{n+1},\boldsymbol{\hat{\delta}_C}^{n+1}))
		+\frac{1}{2}o_h(\bm{u}_h^{n};(\boldsymbol{\delta_C}^{n+1},\boldsymbol{\hat{\delta}_C}^{n+1}),(\boldsymbol{\delta_C}^{n+1},\boldsymbol{\hat{\delta}_C}^{n+1}))
		:=\sum_{i=1}^{5}L_{3i}^C.
	\end{aligned}
	$$ 
	As the same way of estimating $L_3^u$ in \eqref{eq:L3u}, we have the estimates of $\{L_{3i}^C\}_{i=1}^{5}$ that
	\begin{subequations}\label{eq:L3iC}
		\begin{align}
			&|L_{31}^C|
			\leq C\tau\gamma^{-1}
			\|\bm{u}\|_{H^1(t^n,t^{n+1};L^2)}^2
			+\gamma\|\boldsymbol{\delta_C}^{n+1}\|_{L^2}^2,
			\quad
			% \nonumber\\
			% &
			|L_{32}^C|
			\leq
			C\gamma^{-1}\|\boldsymbol{\eta_u}^{n}\|_{L^2}^2
			+\gamma\|\nabla\bm{C}^{n+1}\|_{L^\infty}^2\|\boldsymbol{\delta_C}^{n+1}\|_{L^2}^2,
			\nonumber\\
			& \begin{aligned}
			|L_{33}^C|
			&\leq
			C\gamma^{-1}\|\bm{u}_h^{n}\|_{L^\infty}^2\interleave(\boldsymbol{\eta_C}^{n+1},\boldsymbol{\hat{\eta}_C}^{n+1})\interleave_{w}^2
			+\gamma\|\boldsymbol{\delta_C}^{n+1}\|_{L^2}^2
			\nonumber\\
			&\quad
			+C\gamma^{-1}\|\bm{u}_h^{n}\|_{L^\infty}^2\sum_{K\in\mathcal{T}}\|\boldsymbol{\eta_C}^{n+1}-\boldsymbol{\hat{\eta}_C}^{n+1}\|_{L^2(\partial K)}^2
			+\gamma\sum_{K\in\mathcal{T}}\int_{\partial K}\frac{|\bm{u}_h^n\cdot\bm{n}|}{2}|\boldsymbol{\delta_C}^{n+1}-\boldsymbol{\hat{\delta}_C}^{n+1}|^2\, ds,
			\end{aligned}
			\nonumber\\
			&|L_{34}^C|
			\leq
			C\|\nabla\bm C^{n+1}\|_{L^\infty}^2
			(\|\boldsymbol{\delta_C}^{n+1}\|_{L^2}^2+\|\boldsymbol{\delta_u}^{n}\|_{L^2}^2),
			\quad
			% \nonumber\\
			% &
			|L_{35}^C|=
			\sum_{K\in\mathcal{T}}\int_{\partial K}\frac{|\bm{u}_h^n\cdot\bm{n}|}{4}|\boldsymbol{\delta_C}^{n+1}-\boldsymbol{\hat{\delta}_C}^{n+1}|^2\, ds,
			\nonumber
		\end{align}
	\end{subequations}
	where $C$ depends on $\|\bm{C}^{n+1}\|_{H^2}$. So we have
	\begin{equation}\label{eq:L3C}
		\begin{aligned}
			L_3^C
			&\geq
			(1-\gamma)\sum_{K\in\mathcal{T}}\int_{\partial K}\frac{|\bm{u}_h^n\cdot\bm{n}|}{4}|\boldsymbol{\delta_C}^{n+1}-\boldsymbol{\hat{\delta}_C}^{n+1}|^2\, ds
			-C(\|\boldsymbol{\delta_C}^{n+1}\|_{L^2}^2+ \|\boldsymbol{\delta_u}^{n}\|_{L^2}^2)
			\\
			&\quad
			-C\gamma^{-1}
			\bigg(\tau\|\bm{u}\|_{H^1(t^n,t^{n+1};L^2)}^2
			+\|\boldsymbol{\eta_u}^{n}\|_{L^2}^2
			+\interleave(\boldsymbol{\eta_C}^{n+1},\boldsymbol{\hat{\eta}_C}^{n+1})\interleave_{w}^2
			+\sum_{K\in\mathcal{T}}\|\boldsymbol{\eta_C}^{n+1}-\boldsymbol{\hat{\eta}_C}^{n+1}\|_{L^2(\partial K)}^2\bigg)
			.
		\end{aligned}
	\end{equation}
Now we define the decomposition of $L_4^C$:
$$
\begin{aligned}
	L_4^C
	&=
	\frac{1}{2}\int_{\Omega}({\rm tr}\bm{C}^{n}+{\rm tr}\bm{C}_h^{n})({\rm tr}\bm{C}^{n}-{\rm tr}\bm{C}_h^{n})\bm C^{n+1}:\boldsymbol{\delta_C}^{n+1}\, dx
	+\frac{1}{2}\int_{\Omega}({\rm tr}\bm{C}_h^{n})^2(\bm{C}^{n+1}-\bm{C}_h^{n+1}):\boldsymbol{\delta_C}^{n+1}\, dx\\
    &\quad
    +\frac{1}{2}\int_{\Omega}({\rm tr}\bm{C}^{n+1}-{\rm tr}\bm{C}^{n})({\rm tr}\bm{C}^{n+1}\bm{C}^{n+1} + {\rm tr}\bm{C}^{n+1}\bm{C}^{n}):\boldsymbol{\delta_C}^{n+1}\, dx
	:=\sum_{i=1}^{3}L_{4i}^C.
\end{aligned}
$$
It is easy to show
\[
% \begin{subequations}\label{eq:L4iC}
	\begin{aligned}
		L_{41}^C
		&\leq
        C(\|\bm{C}_h^{n}\|_{L^\infty} + \|\bm{C}^{n}\|_{L^\infty})
        \|\bm C^{n+1}\|_{L^\infty}
        \|\boldsymbol{\delta_C}^{n}+\boldsymbol{\eta_C}^{n}\|_{L^2}
        \|\boldsymbol{\delta_C}^{n+1}\|_{L^2}
		\\
		&\leq
		C(\|\bm{C}_h^{n}\|_{L^\infty} + \|\bm{C}^{n}\|_{L^\infty})^2
		 \|\bm C^{n+1}\|_{L^\infty}^2
		\|\boldsymbol{\delta_C}^{n+1}\|_{L^2}^2
		+C\|\boldsymbol{\delta_C}^{n}\|_{L^2}^2
		+C\|\boldsymbol{\eta_C}^{n}\|_{L^2}^2,
		\\
		L_{42}^C
		&\leq
		C\|\bm{C}_h^{n}\|_{L^\infty}^2\|\boldsymbol{\delta_C}^{n+1} + \boldsymbol{\eta_C}^{n+1}\|_{L^2}
		\|\boldsymbol{\delta_C}^{n+1}\|_{L^2}
		\le
		C\|\bm{C}_h^{n}\|_{L^\infty}^2 
		(\|\boldsymbol{\delta_C}^{n+1}\|^2_{L^2}
		+\| \boldsymbol{\eta_C}^{n+1}\|_{L^2}^2)
		\\
		L_{43}^C
		&\leq
		C(\|\bm{C}^{n+1}\|_{L^\infty} + \|\bm{C}^{n}\|_{L^\infty})^2\tau
		\|\bm{C}\|_{H^1(t^n,t^{n+1};L^2)}^2
		+
		C\|\boldsymbol{\delta_C}^{n+1}\|_{L^2}^2,
	\end{aligned}
% \end{subequations}
\]
% where the constant $C$ depends on $\|\bm{C}^{n+1}\|_{L^\infty}$ and $\|\bm{C}^{n}\|_{L^\infty}$. Above all the boundary is
Thus,
\begin{equation*}\label{eq:L4C}
	|L_4^C|
	\leq
	C(1+\|\bm C_h^n\|^2_{L^\infty})(\tau\|\bm{C}\|_{H^1(t^n,t^{n+1};L^2)}^2
	+\|\boldsymbol{\eta_C}^{n+1}\|_{L^2}^2
	+\|\boldsymbol{\eta_C}^{n}\|_{L^2}^2
	+\|\boldsymbol{\delta_C}^{n+1}\|_{L^2}^2
	+\|\boldsymbol{\delta_C}^{n}\|_{L^2}^2
	).
\end{equation*}
%%%%%%%%%%%%
%%%%%%%%%%%%
    Then we need to offset the positive terms $\sum_{i=1}^{3}R_i^C$ on the right hand. Firstly, there is the decomposition that
    $$
    \begin{aligned}
    	R_1^C&=\frac{1}{2}\int_{\Omega}(\nabla \bm{u}^{n+1})(\bm{C}^{n+1}-\bm{C}^{n})+(\bm{C}^{n+1}-\bm{C}^{n})(\nabla\bm{u}^{n+1})^\top:\boldsymbol{\delta_C}^{n+1}\, dx\\
    	&\quad+\frac{1}{2}\int_{\Omega}(\nabla \bm{u}^{n+1})\boldsymbol{\eta_C}^{n}+\boldsymbol{\eta_C}^{n}(\nabla\bm{u}^{n+1})^\top:\boldsymbol{\delta_C}^{n+1}\, dx
    	+\frac{1}{2}\int_{\Omega}(\nabla \bm{u}^{n+1})\boldsymbol{\delta_C}^{n}+\boldsymbol{\delta_C}^{n}(\nabla\bm{u}^{n+1})^\top:\boldsymbol{\delta_C}^{n+1}\, dx\\
    	&\quad+\frac{1}{2}\int_{\Omega}(\nabla \boldsymbol{\eta_u}^{n+1})\bm{C}_h^{n}+\bm{C}_h^{n}(\nabla\boldsymbol{\eta_u}^{n+1})^\top:\boldsymbol{\delta_C}^{n+1}\, dx
    	+\frac{1}{2}\int_{\Omega}(\nabla \boldsymbol{\delta_u}^{n+1})\bm{C}_h^{n}+\bm{C}_h^{n}(\nabla\boldsymbol{\delta_u}^{n+1})^\top:\boldsymbol{\delta_C}^{n+1}\, dx
    	=:\sum_{i=1}^{5}R_{1i}^C.
    \end{aligned}
    $$
    We estimates each terms as follows, 
\begin{subequations}\label{eq:R1iC}
    \begin{align}
        &\begin{aligned}
    	|R_{11}^C|
    	&\leq\|\boldsymbol{\delta_C}^{n+1}\|_{L^2}
    	\|\nabla \bm{u}^{n+1}\|_{L^\infty}
    	\|\bm{C}^{n+1}-\bm{C}^{n}\|_{L^2}
    	\leq
    	C(\|\boldsymbol{\delta_C}^{n+1}\|_{L^2}^2 
    	+\tau\|\bm{C}\|_{H^1(t^n,t^{n+1};L^2)}^2),
    	\end{aligned}
    	\nonumber
    	\\
    	&|R_{12}^C|\leq
    	C\|\nabla \bm{u}^{n+1}\|_{L^\infty}
    	(\|\boldsymbol{\delta_C}^{n+1}\|_{L^2}^2
    	+\|\boldsymbol{\eta_C}^{n}\|_{L^2}^2),
        % \nonumber
    	% \\
    	% &
    	\quad
    	|R_{13}^C|\leq
    	C\|\nabla \bm{u}^{n+1}\|_{L^\infty}
    	(\|\boldsymbol{\delta_C}^{n+1}\|_{L^2}^2
    	+\|\boldsymbol{\delta_C}^{n}\|_{L^2}^2),\nonumber
    	\\
    	&|R_{14}^C|\leq
    	C(\|\bm{C}_h^n\|_{L^\infty}^2
    	\|\boldsymbol{\delta_C}^{n+1}\|_{L^2}^2
    	+\|\nabla\boldsymbol{\eta_u}^{n+1}\|_{L^2}^2),
    	\quad
    	% \nonumber
    	% \\
    	% &
    	|R_{15}^C|\leq
    	C\gamma^{-1}\nu^{-1}\|\bm{C}_h^n\|_{L^\infty}^2\|\boldsymbol{\delta_C}^{n+1}\|_{L^2}^2
    	+\gamma\nu\| \nabla\boldsymbol{\delta_u}^{n+1}\|_{L^2}^2.
    	\nonumber
    \end{align}
\end{subequations}
As a result, 
%     Also, with lemma \eqref{eq:terms}, there is 
%    $$
%    R_{26}^C=\frac{1}{2}\int_{\Omega}(\nabla \boldsymbol{\delta_u}^{n+1}){\rm tr}\boldsymbol{\delta_C}^{n+1}+{\rm tr}\boldsymbol{\delta_C}^{n+1}(\nabla\boldsymbol{\delta_u}^{n+1})^\top:{\rm tr}\boldsymbol{\delta_C}^{n+1}\, dx
%    =\sum_{K\in\mathcal{T}}\int_{K}({\rm tr}\boldsymbol{\delta_C}^{n+1})\boldsymbol{\delta_C}^{n}:\nabla\boldsymbol{\delta_u}^{n+1}\, dx
%    =L_{26}^u.
%    $$
%    Finally the boundary is

    \begin{equation}\label{eq:R1C}
    	\begin{aligned}
    		|R_1^C|&\leq
    		C(1+\|\bm{C}_h^n\|_{L^\infty}^2)
    		(
    		\tau\|\bm{C}\|_{H^1(t^n,t^{n+1};L^2)}^2
    		+\gamma^{-1}\nu^{-1}\|\boldsymbol{\delta_C}^{n+1}\|_{L^2}^2
    		+\|\boldsymbol{\delta_C}^{n}\|_{L^2}^2
    		+\|\boldsymbol{\eta_C}^{n}\|_{L^2}^2\\
    		&\quad
    		+\|\nabla\boldsymbol{\eta_u}^{n+1}\|_{L^2}^2)
    		+\gamma\nu\| \nabla\boldsymbol{\delta_u}^{n+1}\|_{L^2}^2.
    	\end{aligned}
    \end{equation}
%%%%%%%%%%%%
%%%%%%%%%%%%
For $R_2^C$ and $ R_3^C$, we get
 \begin{equation}\label{eq:R2C}
	|R_2^C|
	\leq C(h^{2k}\|\bm{C}\|_{C^1(0,T;H^{k+1})}^2+\tau\|\bm{C}\|_{H^2(t^n,t^{n+1};L^2)}^2+\|\boldsymbol{\delta_C}^{n+1}\|_{L^2}^2),
\end{equation}
%%%%%%%%%%%%%%
%%%%%%%%%%%%%%
% and 
    \begin{equation}\label{eq:R3C}
    	|R_3^C|\leq
    	C(\tau\|\bm{C}\|_{H^1(t^n,t^{n+1};L^2)}
    	+\|\boldsymbol{\eta_C}^{n}\|_{L^2}
    	+\|\boldsymbol{\delta_C}^{n}\|_{L^2}
    	+\|\boldsymbol{\delta_C}^{n+1}\|_{L^2}).
    \end{equation}
    
    Above all we have the estimates of terms relative to $\boldsymbol{\delta}_C^{n+1},\boldsymbol{\hat{\delta}}_C^{n+1}$. In \textbf{Step 4}, we will sum up all the terms $\sum_{i=1}^4L_i^C,\sum_{i=0}^3R_i^C$ together with terms of $\boldsymbol{\delta}_u^{n+1},\boldsymbol{\hat{\delta}}_u^{n+1}$.
    %Now we have all the estimation of $\boldsymbol{\delta_C},\boldsymbol{\hat{\delta}_C}$.
%%%%%%%%%%%%%%%%%%%%%%%%%
%%%%%%%%%%%%%%%%%%%%%%%%%
%%%%%%%%%%%%%%%%%%%%%%%%%
%%%%%%%%%%%%%%%%%%%%%%%%%

(\textbf{Step 4})\label{step 4} 
Combining the estimate of (\textbf{Step 2}) and (\textbf{Step 3}), we have: 
%Because of the boundary of $\bm{u}$ \eqref{def:error-space},
%\begin{equation}
%	\begin{aligned}
%		&\sum_{K\in\mathcal{T}}\int_{\partial K}|\boldsymbol{\eta_u}^{n+1}-\boldsymbol{\hat{\eta}_u}^{n+1}|^2\, ds
%		\leq Ch^{2k+1}|\bm{u}|_{k+1}^2,\\
%		&\sum_{K\in\mathcal{T}}\int_{\partial K}|{\rm tr}\boldsymbol{\eta_C}^{n+1}-{\rm tr}\boldsymbol{\hat{\eta}_C}^{n+1}|^2\, ds
%		\leq Ch^{2k+1}|{\rm tr}\bm{C}|_{k+1}^2.
%	\end{aligned}
%\end{equation}
$$
	LHS\leq RHS,
$$
where
\begin{equation}\label{eq:LHS}
	\begin{aligned}
		LHS&=
		\frac{1}{2\tau} (\|\boldsymbol{\delta_u}^{n+1}\|_{L^2}^2 - \|\boldsymbol{\delta_u}^{n}\|_{L^2}^2)
		+ \frac{1}{4\tau} (\|\boldsymbol{\delta_C}^{n+1}\|_{L^2}^2 - \|\boldsymbol{\delta_C}^{n}\|_{L^2}^2)
		+\frac{1}{2\tau}\|\boldsymbol{\delta_u}^{n+1}-\boldsymbol{\delta_u}^{n}\|_{L^2}^2
		\\&\quad
		+\frac{1}{4\tau}\|\boldsymbol{\delta_C}^{n+1}-\boldsymbol{\delta_C}^{n}\|_{L^2}^2
		+\nu(1-7\gamma)\|\nabla\boldsymbol{\delta_u}^{n+1}\|_{L^2}^2
		+\frac{\ep}{2}(1-\gamma)\|\nabla\boldsymbol{\delta_C}^{n+1}\|_{L^2}^2
		\\
		&\quad
		+\sum_{K\in\mathcal{T}}\frac{\nu}{h_K}(\alpha-\frac{c_I}{\gamma}-5\gamma)\|\boldsymbol{\delta_u}^{n+1}-\boldsymbol{\hat{\delta}_u}^{n+1}\|_{L^2(\partial K)}^2
		\\
		&\quad
		+\sum_{K\in\mathcal{T}}\frac{\epsilon}{2h_K}(\beta-\frac{c_I}{\gamma})\|\boldsymbol{\delta_C}^{n+1}-\boldsymbol{\hat{\delta}_C}^{n+1}\|_{L^2(\partial K)}^2
		+\frac{1-\gamma}{2}\||\bm{u}_h^{n}\cdot\bm{n}|^{\frac{1}{2}}(\boldsymbol{\delta_u}^{n+1}-\boldsymbol{\hat{\delta}_u}^{n+1})\|_{L^2(\mathcal{F})}^2
		\\
		&\quad
		+\frac{1-2\gamma}{4}\||\bm{u}_h^{n}\cdot\bm{n}|^{\frac{1}{2}}(\boldsymbol{\delta_C}^{n+1}-\boldsymbol{\hat{\delta}_C}^{n+1})\|_{L^2(\mathcal{F})}^2,
	\end{aligned}
\end{equation}
and
\begin{equation}\label{eq:RHS}
	\begin{aligned}
		RHS&=
		2C(\|\boldsymbol{\delta_u}^{n}\|_{L^2}^2 + \|\boldsymbol{\delta_u}^{n+1}\|_{L^2}^2)
		+ C(3\nu^{-1}\gamma^{-1}+4)\|\boldsymbol{\delta_C}^{n+1}\|_{L^2}^2
		+C(2\nu^{-1}\gamma^{-1}+3)\|\boldsymbol{\delta_C}^{n}\|_{L^2}^2
		\\
		&\quad
		+2\gamma^{-1}Ch^{2k+2}|\bm{u}^{n}|_{k+1}^2
		+3\gamma^{-1}Ch^{2k}|\bm{u}^{n+1}|_{k+1}^2
		+4\nu^{-1}\gamma^{-1}Ch^{2k+2}|\bm{C}^{n+1}|_{k+1}^2
		\\
		&\quad
		+2\gamma^{-1}Ch^{2k}|\bm{C}^{n+1}|_{k+1}^2
		+6\nu^{-1}\gamma^{-1}h^{2k+2}|\bm{C}^{n}|_{k+1}^2\\
		&\quad
		+Ch^{2k}\|\bm{u}\|_{C^1(0,T;H^{k+1})}^2
		+2C\gamma^{-1}\tau\|\bm{u}\|_{H^1(t^n,t^{n+1};L^2)}^2
		+C\tau\|\bm{u}_{tt}\|_{L^2(t^n,t^{n+1};L^2)}^2\\
		&\quad
		+Ch^{2k}\|\bm{C}\|_{C^1(0,T;H^{k+1})}^2
		+C\tau(3+\nu^{-1}\gamma^{-1})\|\bm{C}\|_{H^1(t^n,t^{n+1};L^2)}^2
		+C\tau\|\bm{C}\|_{H^2(t^n,t^{n+1};L^2)}^2\\
		&\quad
		+C\nu^{-1}\gamma^{-1}\tau\|\bm{C}_{t}\|_{L^2(t^n,t^{n+1};L^2)}^2
		+C\nu^{-1}\gamma^{-1}h^2\|\nabla\bm{C}\|_{H^1(t^n,t^{n+1};L^2)}^2
		.
	\end{aligned}
\end{equation}
First, we choose $\gamma\in(0,\frac{1}{7})$ and  then take suﬀiciently large $\alpha$ and $\beta$ such that $\alpha-\frac{c_I}{\gamma}-5\gamma>0$, and $\beta-\frac{c_I}{\gamma}>0$. 
Applying the discrete Gronwall inequality, we have:
\begin{equation}\label{eq:error-delta}
	\begin{aligned}
		&\|\boldsymbol{\delta_u}\|_{l^\infty(0,T;L^2)}^2
		+\|\boldsymbol{\delta_C}\|_{l^\infty(0,T;L^2)}^2
		+\sum_{i=1}^{N}\|\boldsymbol{\delta_u}^{i+1}-\boldsymbol{\delta_u}^{i}\|_{L^2}^2
		+\sum_{i=1}^{N}\|\boldsymbol{\delta_C}^{i+1}-\boldsymbol{\delta_C}^{i}\|_{L^2}^2\\
		&+\nu\|\nabla\boldsymbol{\delta_u}\|_{l^2(0,T;L^2)}^2
		+\frac{\nu}{h}\|\boldsymbol{\delta_u}-\boldsymbol{\hat{\delta}_u}\|_{l^2(0,T;L^2(\mathcal{F}))}^2
		+\epsilon\|\nabla\boldsymbol{\delta_C}\|_{l^2(0,T;L^2)}^2
		+\frac{\epsilon }{h}\|\boldsymbol{\delta_C}-\boldsymbol{\hat{\delta}_C}\|_{l^2(0,T;L^2(\mathcal{F}))}^2
		\\
		&+\sum_{i=1}^{N}\tau\||\bm{u}_h^{i}\cdot\bm{n}|^{\frac{1}{2}}(\boldsymbol{\delta_u}^{i+1}-\boldsymbol{\hat{\delta}_u}^{i+1})\|_{L^2(\mathcal{F})}^2
		+\sum_{i=1}^{N}\tau\||\bm{u}_h^{i}\cdot\bm{n}|^{\frac{1}{2}}(\boldsymbol{\delta_C}^{i+1}-\boldsymbol{\hat{\delta}_C}^{i+1})\|_{L^2(\mathcal{F})}^2
		\leq
		C(h^{2k}+\tau^2),
	\end{aligned}
\end{equation}
where $C$ depends on $\|\bm u_h^n\|_{L^\infty}$, $\|\bm C_h^n\|_{L^\infty}$ and the norm of $(\bm u,p,\bm C)$, but is independent of $h$, $\tau$, $n$ and $\epsilon^{-1}$.
By the triangle inequality and the estimate of $\bm{\eta_u}$ and $\bm{\eta_C}$, we obtain \eqref{eq:error-uC}.
%%%%%%%%%%%%
%%%%%%%%%%%%

(\textbf{Step 5})\label{step5} Now we estimate the error of $p$. By the inf-sup condition \eqref{eq:inf} and \eqref{eq:diff-b}, we have:
$$
	\begin{aligned}
		&\interleave(p^{n+1},p^{n+1})-(p_h^{n+1},\hat{p}_h^{n+1})\interleave_{q}
		\leq
		\interleave(\eta_p^{n+1},\hat{\eta}_p^{n+1})\interleave_{q}
		+\interleave(\delta_p^{n+1},\hat{\delta}_p^{n+1})\interleave_{q}\\
		\leq
		&\interleave(\eta_p^{n+1},\hat{\eta}_p^{n+1})\interleave_{q}
		+\sigma^{-1}\sup_{\bm{v}_h\in \bm{V}_h,\hat{\bm{v}}_h\in \hat{\bm{V}}_h}
		\frac{ |b_h((\delta_p^{n+1},\hat{\delta}_p^{n+1}),(\bm{v}_h,\hat{\bm{v}}_h))}{\interleave(\bm{v}_h,\hat{\bm{v}}_h)\interleave_v}|,\\
	\end{aligned}
$$
where 
$$
	\begin{aligned}
		&|b_h((\delta_p^{n+1},\hat{\delta}_p^{n+1}),(\bm{v}_h,\hat{\bm{v}}_h))|
		\leq(\frac{\boldsymbol{\delta_u}^{n+1}-\boldsymbol{\delta_u}^{n}}{\tau},\bm{v}_h)
		+a_h((\boldsymbol{\delta_u}^{n+1},\boldsymbol{\hat{\delta}_u}^{n+1}),(\bm{v}_h,\hat{\bm{v}}_h))
		\\
		& \quad +\bigg(o_h(\bm{u}^{n+1};(\bm{u}^{n+1},\bm{u}^{n+1}),(\bm{v}_h,\hat{\bm{v}}_h))
		-o_h(\bm{u}_h^{n};(\bm{u}_h^{n+1},\hat{\bm{u}}_h^{n+1}),(\bm{v}_h,\hat{\bm{v}}_h))\bigg)
		\\
		& \quad +(\frac{\mathcal{R}_{h}^{u}\bm{u}^{n+1}-\mathcal{R}_{h}^{u}\bm{u}^{n}}{\tau}-\partial_t\bm{u}^{n+1},\boldsymbol{\delta_u}^{n+1})
		+\bigg(\sum_{K\in\mathcal{T}}\int_{K}(({\rm tr}\bm{C}^{n+1})\bm{C}^{n+1}-({\rm tr}\bm{C}_h^{n+1})\bm{C}_h^{n}):\nabla\bm{v}_h\, dx
		\\
		&\quad +\sum_{K\in\mathcal{T}}\int_{\partial K}(({\rm tr}\bm{C}^{n+1})\bm{C}^{n+1}-({\rm tr}\bm{C}_{h}^{n+1})\bm{C}_{h}^{n}):(\bm{v}_h-\hat{\bm{v}}_h)\otimes\bm{n}\,d s
		\bigg)
		\\\leq &C \left(\frac{1}{\tau}\|\boldsymbol{\delta_u}^{n+1}-\boldsymbol{\delta_u}^{n}\|_{L^2}
		+\nu\interleave(\boldsymbol{\delta_u}^{n+1},\boldsymbol{\hat{\delta}_u}^{n+1})\interleave_{v}
		+\bigg( C\|\bm{u}^{n+1} - \bm{u}^{n} + \bm{e_u}^{n}\|_1
		\interleave(\bm u^{n+1},\bm u^{n+1})
		\interleave_v
		\right.
		\\
		&\left.\quad  + \|\bm{u}_h^n\|_{L^\infty}
		\interleave(\bm{e_u}^{n+1},\boldsymbol{\hat{e_u}}^{n+1})\interleave_{0,v}
		\bigg)
		+h^{k}\|\bm{u}\|_{C^1(0,T;H^{k+1})}
		+\tau^{\frac{1}{2}}\|\bm{u}\|_{H^2(t^n,t^{n+1};L^2)}\right.\\
		&\left. \quad + \tau^{\frac{1}{2}}\|\bm{C}\|_{H^1(t^n,t^{n+1};L^2)}
		+\|\bm{e_C}^{n+1}\|_{L^2}
		+\|\bm{e_C}^{n}\|_{L^2}
		\right)
		\interleave(\bm{v}_h,\hat{\bm{v}}_h)\interleave_{v},
		\quad (\text{by Lemma}~\ref{le:conv})
	\end{aligned}
$$
where $C$ depends on $\|\bm u_h^n\|_{L^\infty}$  and $\|\bm C_h^n\|_{L^\infty}$.
By \eqref{eq:error-uC}, the projection error and the triangle inequality, we obtain
\begin{equation}
	\begin{aligned}
		\sum_{i=1}^{n}\tau^2\interleave(p^{i+1},p^{i+1})-(p_h^{i+1},\hat{p}_h^{i+1})\interleave_{q}^2
		&\leq C(h^{2k}+\tau^2),
	\end{aligned}
\end{equation}
where $C$ depends on $\|\bm u_h^i\|_{L^\infty}$  and $\|\bm C_h^i\|_{L^\infty}$ $(i = 0,\ldots,n)$, but is independent of $h$, $\tau$, $n$ and $\epsilon^{-1}$.
Hence, the proof is complete.
\end{proof}
\begin{remark}
The above proof is under the assupmtion that $\|\bm{u}_h^n\|_{L^\infty} + \|\bm{C}_h^n\|_{L^\infty}\leq \bar{C}$ with the constant $\bar{C}$ independent of $n$. For $k\ge2$, suppose we have the projection error in $L^\infty$-norm, i.e.,
$$
\|\mathcal{R} _h^{u}\bm{u}-\bm{u}\|_{L^\infty}
\le Ch^a\|\bm u\|_{H^2}
\text{ and }
\|\mathcal{R} _h^{C}\bm{C}-\bm{C}\|_{L^\infty}
\le Ch^a\|\bm C\|_{H^2}
\text{ for some } a>0.
$$
Then, we can remove the assumption \eqref{eq:con1} and prove the same convergence result by using the induction argument.
In fact, with the initial values satisfying $\|\bm{u}_h^0\|_{L^\infty} + \|\bm{C}_h^0\|_{L^\infty}\leq \bar{C}:=2(\|\bm u\|_{L^\infty(0,T;L^\infty)} + \|\bm C\|_{L^\infty(0,T;L^\infty)}\leq \bar{C})$, we first make the induction assumption that 
$$
\|\bm{u}_h^i\|_{L^\infty} + \|\bm{C}_h^i\|_{L^\infty}\leq \bar{C}
\text{ for } i = 1,2,\ldots,n.
$$
Then we can prove the error estimates \eqref{eq:error-uC} and \eqref{eq:error-p}.
It then suffices to validate $\|\bm{u}_h^{n+1}\|_{L^\infty} + \|\bm{C}_h^{n+1}\|_{L^\infty}\leq \bar{C}$. 
This can be achieved using the triangle inequality and inverse inequality:
$$
\begin{aligned}
\|\bm u_h^{n+1}\|_{L^\infty}
&\le \|\bm{\delta_u}^{n+1}\|_{L^\infty} 
+ \|\mathcal{R} _h^{u}\bm{u}-\bm u^{n+1}\|_{L^\infty} + \|\bm u^{n+1}\|_{L^\infty}
\\
&\le 
Ch^{-1}\|\bm{\delta_u}^{n+1}\|_{L^2}
+ Ch^a\|\bm u^{n+1}\|_{H^2}
+ \|\bm u^{n+1}\|_{L^\infty}
\le 
C(h^{k-1} + \tau h^{-1} + h^a) + \|\bm u^{n+1}\|_{L^\infty}.
\end{aligned}
$$
Similarly
$$
\|\bm C_h^{n+1}\|_{L^\infty}
\le 
C(h^{k-1} + \tau h^{-1} + h^a) + \|\bm C^{n+1}\|_{L^\infty}.
$$
Hence, for sufficiently small $h$ and $\tau=h^{1+\ep_0}$ for some $\ep_0>0$, we can ensure that $\|\bm u_h^{n+1}\|_{L^\infty} + \|\bm C_h^{n+1}\|_{L^\infty}\le \bar{C}$.
\end{remark}
%%%%%%%%%%%%%%%%%%%%%%%%%%%%%%%%%%%%
% % %%%%%%%%%%%%%%%%%%%%%%%%%%%%%%%%%%%%
\section{Numerical experiments}\label{sec:ex}
In this section, two numerical examples are carried out to investigate the applicability of the proposed scheme. 
In Example~\ref{sec5.1}, we study the experimental errors and convergence rates in three cases ($\ep = 1$, $10^{-3}, 0$). Then, we perform a comparative simulation of the ordinary finite element discretization to show the stability of our proposed scheme in Example~\ref{sec5.2}.
\subsection{Example 1}\label{sec5.1}
We consider the problem \eqref{eq:P} with the external forces $\bm f$ and $\bm F$ added in \eqref{eq:P-a} and \eqref{eq:P-c}, respectively. 
Set a unit square domain $\Omega=(0,1)^2$, the final time $T=0.2$ and $\nu = 1$. The analytical solutions are given as follows:
\begin{equation}\label{eq:exa}
	\begin{aligned}
		&\bm{u}(x,t)=\left(-\frac{\partial \phi }{\partial x_2},\frac{\partial \phi }{\partial x_1}\right),
		\quad\text{with }\phi(x,t)=\frac{\sqrt{3}}{2\pi}\sin^2(\pi x_1)\sin^2(\pi x_2)\sin(\pi(x_1+x_2+t)),\\
		&p(x,t)=\sin(\pi(x_1+2x_2+t)),\quad C_{11}(x,t)=\frac{1}{2}\sin^2(\pi x_1)\sin^2(\pi x_2)\sin(\pi(x_1+t))+1,\\
		&C_{22}(x,t)=\frac{1}{2}\sin^2(\pi x_1)\sin^2(\pi x_2)\sin(\pi(x_2+t))+1,\quad C_{12}(x,t)=\frac{\pi}{\sqrt{3}}\phi(x,t)=C_{21}(x,t),
	\end{aligned}
\end{equation}
where $\bm f$ and $\bm F$, initial value $\bm u^0$ and $\bm C^0$ can be calculated when given the value of $\ep$.

We conduct the simulation for three cases: (i)$\epsilon=1$, (ii) $\epsilon=10^{-3}$ and (iii) $\epsilon=0$.  
We set $k =1$(the P1 element for $(\bm u_h,  \hat{\bm u}_h, p_h, \hat{p}_h, \bm C_h, \hat{\bm C}_h)$) and compute the experimental errors at $T=0.2$. 
%%%%%%%%%%%%%%%%%%%%%%%%%%%%%%%%%%%% ep=1
\subsubsection{$\epsilon=1$}\label{sec5.1.1}
\begin{figure}[htbp] 
\begin{tabular}{ccc}
\centering
\includegraphics[width=0.3\textwidth]{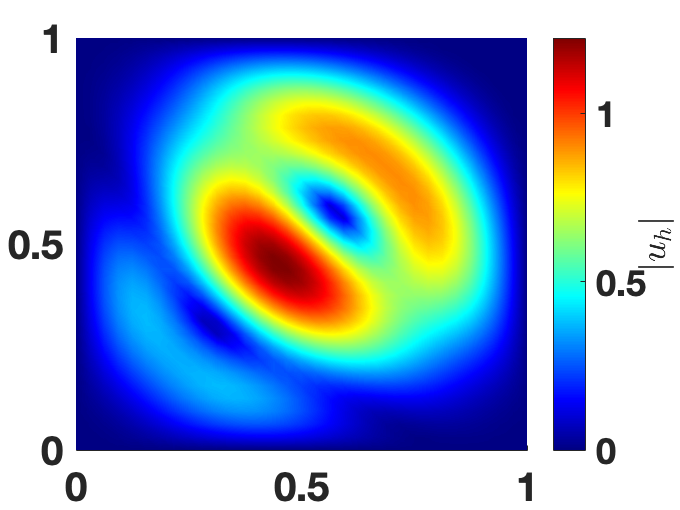} &
\includegraphics[width=0.3\textwidth]{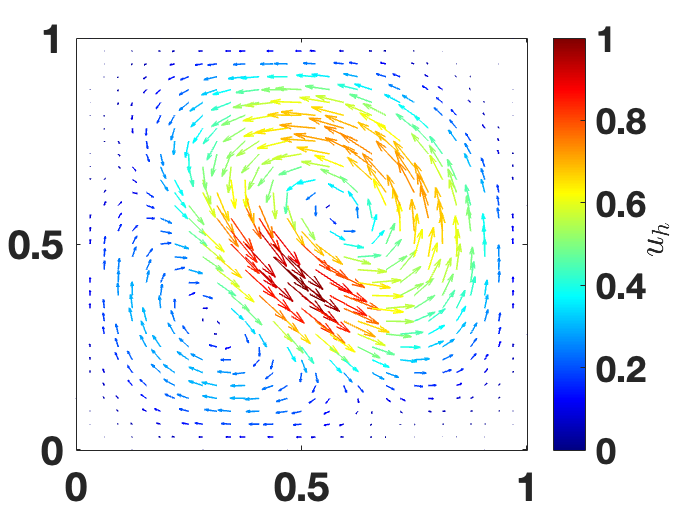} &
\includegraphics[width=0.3\textwidth]{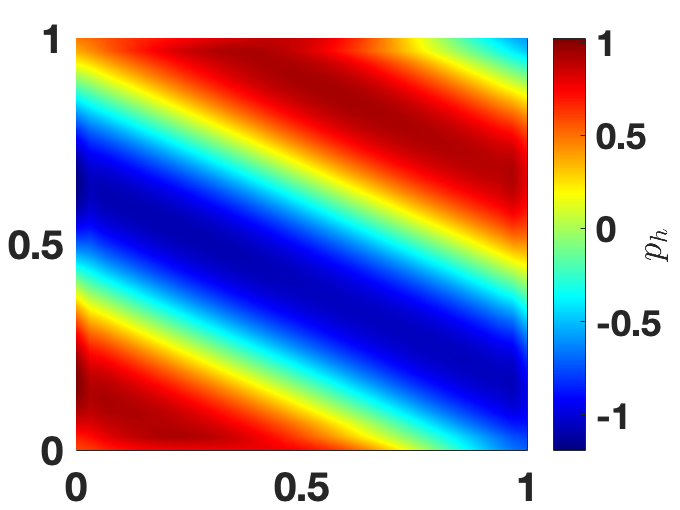} 
\\
(a) $|\bm u_h^N|$ &
(b) $\bm u_h^N$ &
(c) $p_h^N$
\\
\includegraphics[width=0.3\textwidth]{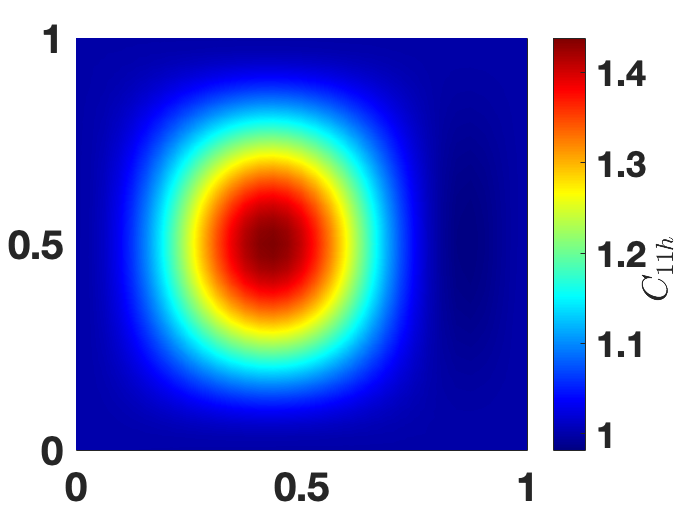} &
\includegraphics[width=0.3\textwidth]{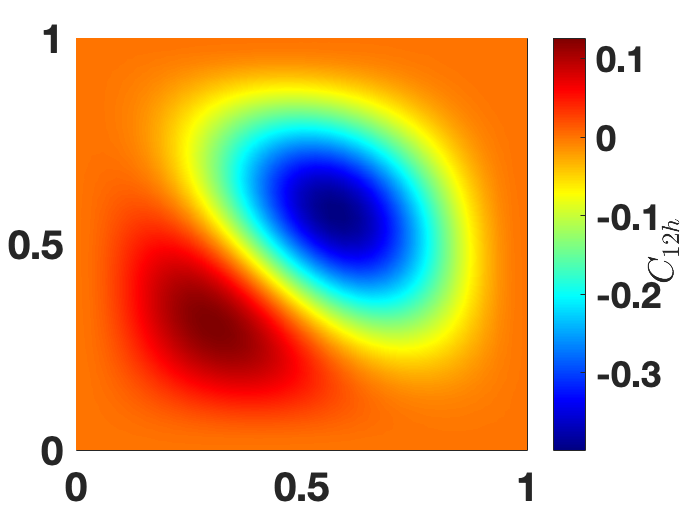} &
\includegraphics[width=0.3\textwidth]{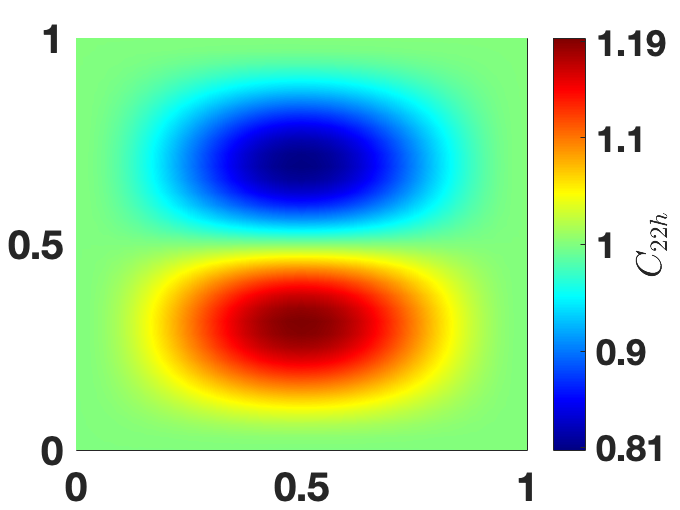} 
\\
(d) $C_{11h}^N$ &
(e) $C_{12h}^N$ &
(f) $C_{22h}^N$
\\
\end{tabular}
\caption{$\epsilon=1$: The numerical solution $(\bm{u}_h,p_h,\bm{C}_h)$ at $t=0.2$}\label{fig1}
\end{figure}

We plot the numerical solution $(\bm{u}_h,p_h,\bm{C}_h = (C_{ijh})_{1\le i,j\le 2})$ in Figure~\ref{fig1}. 
Then we verify the convergence behavior of mesh size $h$ and time-step increment $\tau$ with the penalty coefficients $\alpha = 8$ and $\beta = 10$. 

We first fix the time-step increment $\tau = 2^{-12}$ and compute the experimental errors and convergence rates with mesh sizes $h = 2^{-2}, 2^{-3}, 2^{-4},2^{-5}$. The results can be found in Table~\ref{table1} and Table~\ref{table2}, which implies the $O(h)$-convergence for $L^2$-norm error of $p_h$ and $H^1$-norm error of $(\bm{u}_h, \bm{C}_h)$.
Meanwhile, we observe $O(h^2)$-convergence for $\|\bm{u}_{h}^{N} - \bm{u}(T)\|_{L^2}$ and $\|\bm{C}_h^{N}-\bm{C}(T)\|_{L^2}$, which are better than the theoretical predictions, indicating the possibility of obtaining a higher order convergence rate for the case $\epsilon=0$.

In addition, under fixed $h = 2^{-6}$, we calculate the experimental errors in Table~\ref{table3} with different time-step size ($\tau=10^{-1},20^{-1},40^{-1},80^{-1}$). The results shows $O(\tau)$-convergence for $\|\bm{u}_{h}^{N} - \bm{u}(T)\|_{L^2}$ and $\|\bm{C}_h^{N}-\bm{C}(T)\|_{L^2}$, confirming the theoretical convergence rate of Theorem~\ref{th:error}.
%%%%%%%%%%%%%%%%%%%%%%%%%%%%%%%%%%%%
%%%%%%%%%%%%%%%%%%%%%%%%%%%%%%%%%%% 
\begin{table}[htbp] 
	\caption{$\epsilon=1$: The errors of $(\bm u_h ,p_h)$ with $\tau = 2^{-12}$}
	\begin{tabular}{ccccccc}
        \toprule
        $h$ & $\|\bm{u}_{h}^{N} - \bm{u}(T)\|_{L^2}$ & Rate & $\|\bm{u}_{h}^{N} - \bm{u}(T)\|_{H^1} $ & Rate & $\|p_{h}^{N} - p(T)\|_{L^2} $ & Rate \\
        \midrule
        $2^{-2}$ & $1.01\times 10^{-1}$ & - & $2.66\times 10^{0}$ & - & $1.72\times 10^{0}$ & - 
        \\
        $2^{-3}$& $2.73\times 10^{-2}$ & 1.88 & $1.40\times 10^{0}$ & 0.93 & $8.48\times 10^{-1}$ & 1.02
        \\
        $2^{-4}$ &$6.96\times 10^{-3}$ & 1.97 & $6.98\times 10^{-1}$ & 1.00 & $4.16\times 10^{-1}$ & 1.03 
        \\
        $2^{-5}$ &$1.75\times 10^{-3}$ & 1.99 & $3.48\times 10^{-1}$ & 1.00 & $2.03\times 10^{-1}$ & 1.04
        \\
        \bottomrule
    \end{tabular}
	\label{table1}
\end{table}
\begin{table}[htbp] 
	\caption{$\epsilon=1$: The errors of $\bm{C}_h$ with $\tau = 2^{-12}$ }
	\begin{tabular}{ccccc}
        \toprule
        $h$ &$\|\bm{C}_h^{N}-\bm{C}(T)\|_{L^2}$ & Rate & $\sqrt{\epsilon}\|\bm{C}_h^{N}-\bm{C}(T)\|_{H^1}$ & Rate\\
        \midrule
        $2^{-2}$ &$5.60\times 10^{-2}$ &- &$9.76\times 10^{-1}$ &-
        \\
        $2^{-3}$ &$1.39\times 10^{-2}$ & 2.01 &$4.81\times 10^{-1}$ & 1.02
        \\
        $2^{-4}$ &$3.49\times 10^{-3}$ & 2.00 &$2.39\times 10^{-1}$ & 1.01
        \\
        $2^{-5}$ &$8.87\times 10^{-4}$ & 1.98 & $1.19\times 10^{-1}$ & 1.00
        \\
        \bottomrule
    \end{tabular}
	\label{table2}
\end{table}
%%%%%%%%%%%%%%%%%%%%%%%%%%%%%%%%%%%% h
\begin{table}[htbp] 
	\caption{$\epsilon=1$: The errors of $(\bm u_h ,\bm C_h)$ with $h = 2^{-6}$}
	\begin{tabular}{ccccc}
        \toprule
        $\tau$ & $\|\bm{u}_{h}^{N} - \bm{u}(T)\|_{L^2}$ & Rate & $\|\bm{C}_h^{N}-\bm{C}(T)\|_{L^2}$ & Rate \\
        \midrule
        $10^{-1}$& $5.24\times 10^{-3}$ & - & $2.53\times 10^{-2}$  & -
        \\
        $20^{-1}$ &$2.81\times 10^{-3}$ & 0.90 & $1.24\times 10^{-2}$ & 1.03
        \\
        $40^{-1}$ &$1.51\times 10^{-3}$ & 0.90 & $6.16\times 10^{-3}$ & 1.01
        \\
        $80^{-1}$ &$8.74\times 10^{-4}$ & 0.79 & $3.11\times 10^{-3}$ & 0.99
        \\
        \bottomrule
    \end{tabular}
	\label{table3}
\end{table}
%%%%%%%%%%%%%%%%%%%%%%%%%%%%%%%%%%%% ep=10^-3
\subsubsection{$\epsilon=10^{-3}$}\label{sec5.1.2}
In this case, we take a smaller $\epsilon$ and check the corresponding numerical result.
Set the penalty coeﬀicients $\alpha = 8$ and $\beta = 300$.
We omit the figures of numerical solutions, which is similar to the case $\ep =1$.

To study the convergence rates, we fix $\tau= 2^{-12}$ and conduct the simulation with various mesh sizes in Table~\ref{table4} and Table~\ref{table5}. 
The $L^2$-norm error of $p_h$ and $H^1$-norm error of $\bm u_h$ are first-order convergence.
Furthermore, we observe that the $L^2$-norm error of $(\bm u_h,\bm C_h)$ exhibits nearly $O(h^2)$-convergence, and the $H^1$-norm error of $\bm C_h$ is better than $O(h)$-convergence, which is superior to the theoretical prediction.
Next, we test the convergence rates of $\tau$ and list the experimental results in Table~\ref{table6}. The $L^2$-norm errors of $(\bm u_h,\bm C_h)$ are about the first-order convergence.
%%%%%%%%%%%%%%%%%%%%%%%%%%%%%%%%%%%% h
\begin{table}[htbp] 
	\caption{$\epsilon=10^{-3}$: The errors of $(\bm u_h ,p_h)$ with $\tau = 2^{-12}$}
	\begin{tabular}{ccccccc}
		\toprule
		$h$ & $\|\bm{u}_{h}^{N} - \bm{u}(T)\|_{L^2}$ & Rate & $\|\bm{u}_{h}^{N} - \bm{u}(T)\|_{H^1} $ & Rate & $\|p_{h}^{N} - p(T)\|_{L^2} $ & Rate \\
		\midrule
		$2^{-2}$ & $9.63\times 10^{-2}$ & - & $2.65\times 10^{0}$ & - & $1.67\times 10^{0}$ & - 
		\\
		$2^{-3}$& $2.53\times 10^{-2}$ & 1.93 & $1.40\times 10^{0}$ & 0.92 & $8.46\times 10^{-1}$ & 0.98
		\\
		$2^{-4}$ &$6.37\times 10^{-3}$ & 1.99 & $6.98\times 10^{-1}$ & 1.00 & $4.16\times 10^{-1}$ & 1.02 
		\\
		$2^{-5}$ &$1.60\times 10^{-3}$ & 2.00 & $3.48\times 10^{-1}$ & 1.00 & $2.03\times 10^{-1}$ & 1.03
		\\
		\bottomrule
	\end{tabular}
	\label{table4}
\end{table}
\begin{table}[htbp] 
	\caption{$\epsilon=10^{-3}$: The errors of $\bm C_h$ with $\tau = 2^{-12}$}
	\begin{tabular}{ccccc}
		\toprule
		$h$ &$\|\bm{C}_h^{N}-\bm{C}(T)\|_{L^2}$ & Rate & $\sqrt{\epsilon}\|\bm{C}_h^{N}-\bm{C}(T)\|_{H^1}$ & Rate\\
		\midrule
		$2^{-2}$ &$3.52\times 10^{-1}$ &- &$1.79\times 10^{-1}$ &-
		\\
		$2^{-3}$ &$7.29\times 10^{-2}$ & 2.27 &$5.44\times 10^{-2}$ & 1.71
		\\
		$2^{-4}$ &$1.72\times 10^{-2}$ & 2.09 &$1.91\times 10^{-2}$ & 1.51
		\\
		$2^{-5}$ &$4.30\times 10^{-3}$ & 2.00 &$7.78\times 10^{-3}$ & 1.29
		\\
		\bottomrule
	\end{tabular}
	\label{table5}
\end{table}
%%%%%%%%%%%%%%%%%%%%%%%%%%%%%%%%%%%% h
\begin{table}[htbp] 
	\caption{$\epsilon=10^{-3}$: The errors of $(\bm u_h ,\bm C_h)$ with $h = 2^{-6}$}
	\begin{tabular}{ccccc}
        \toprule
        $\tau$ & $\|\bm{u}_{h}^{N} - \bm{u}(T)\|_{L^2}$ & Rate & $\|\bm{C}_h^{N}-\bm{C}(T)\|_{L^2}$ & Rate \\
        \midrule
        $10^{-1}$ & $5.69\times 10^{-3}$ & - & $5.65\times 10^{-2}$ & - 
        \\
        $20^{-1}$ & $3.15\times 10^{-3}$ & 0.85& $2.80\times 10^{-2}$ & 1.01 
        \\
        $40^{-1}$ &$1.61\times 10^{-3}$ & 0.97 & $1.40\times 10^{-2}$ & 1.00
        \\
        $80^{-1}$ &$8.57\times 10^{-4}$ & 0.91 & $7.05\times 10^{-3}$ & 0.98
        \\
        \bottomrule
    \end{tabular}
	\label{table6}
\end{table}
%%%%%%%%%%%%%%%%%%%%%%%%%%%%%%%%%%%% ep=0
\subsubsection{$\epsilon=0$}\label{sec5.1.3}
The figures of the discrete solutions are almost same to the case $\ep =1$.
The experimental errors and convergence rates are investigated with $\alpha=8$ and $\beta=10$.   
We fix $\tau=2^{-12}$ and compute the HDG scheme on different $h$, similar to the above cases. The reuslts are displayed in Table~\ref{table7}, which indicates the $O(h)$-convergence for $\|\bm{u}_{h}^{N}- \bm{u}(T)\|_{H^1} $ and $\|p_{h}^{N} - p(T)\|_{L^2} $.
Similarly, the $L^2$-norm error of $(\bm u_h, \bm C_h)$ is nearly $O(h^2)$-convergence, better than the theoretical prediction of Theorem~\ref{th:error}. The convergence rates of $\|\bm{C}_h^{N}-\bm{C}(T)\|_{L^2}$ decrease when $h$ is too small.  
And Table~\ref{table8} presents the $O(\tau)$-convergence for the $L^2$-norm error of $(\bm u_h,\bm C_h)$.
%%%%%%%%%%%%%%%%%%%%%%%%%%%%%%%%%%%%
\begin{table}[htbp] 
	\caption{$\epsilon=0$: The errors of $(\bm u_h ,p_h,\bm C_h)$ with $\tau = 2^{-12}$}
	\setlength{\tabcolsep}{0.45mm}
	\begin{tabular}{ccccccccc}
        \toprule
        $h$ & $\|\bm{u}_{h}^{N} - \bm{u}(T)\|_{L^2}$ & Rate & $\|\bm{u}_{h}^{N}- \bm{u}(T)\|_{H^1} $ & Rate & $\|p_{h}^{N} - p(T)\|_{L^2} $ & Rate &$\|\bm{C}_h^{N}-\bm{C}(T)\|_{L^2}$ & Rate\\
        \midrule
        $2^{-2}$ & $9.69\times 10^{-2}$ & - & $2.68\times 10^{0}$ & - & $1.71\times 10^{0}$ & - &$3.22\times10^{-1}$ &- 
        \\
        $2^{-3}$& $2.52\times 10^{-2}$ & 1.94 & $1.40\times 10^{0}$ & 0.94 & $8.52\times 10^{-1}$ & 1.01 &$6.70\times10^{-2}$ & 2.20 
        \\
        $2^{-4}$ &$6.35\times 10^{-3}$ & 1.99 & $6.98\times 10^{-1}$ & 1.00 & $4.16\times 10^{-1}$ & 1.03   &$1.82\times10^{-2}$ &1.94
        \\
        $2^{-5}$ &$1.59\times 10^{-3}$ & 1.99 & $3.48\times 10^{-1}$ & 1.00 & $ 2.03\times 10^{-1}$ & 1.03   &$5.33\times10^{-3}$ &1.77
        \\
        \bottomrule
    \end{tabular}
	\label{table7}
\end{table}
%%%%%%%%%%%%%%%%%%%%%%%%%%%%%%%%%%%%
\begin{table}[H] 
\centering
    \caption{$\epsilon=0$: The errors of $(\bm u_h , \bm C_h)$ with $h = 2^{-6}$}
    \setlength{\tabcolsep}{0.45mm}
    \begin{tabular}{ccccc}
        \toprule
        $\tau$ & $\|\bm{u}_{h}^{N} - \bm{u}(T)\|_{L^2}$ & Rate & $\|\bm{C}_h^{N}-\bm{C}(T)\|_{L^2}$ & Rate\\
        \midrule
        $10^{-1}$ & $5.70\times 10^{-3}$ & - & $5.71\times 10^{-2}$ & -
        \\
        $20^{-1}$ & $3.16\times 10^{-3}$ & 0.85 & $2.83\times 10^{-2}$ & 1.01
        \\
        $40^{-1}$ &$1.61\times 10^{-3}$ & 0.97 & $1.41\times 10^{-2}$ & 1.00
        \\
        $80^{-1}$ &$8.58\times 10^{-4}$ & 0.91 & $7.23\times 10^{-3}$ & 0.97
        \\
        \bottomrule
    \end{tabular}
    \label{table8}
\end{table}

\subsection{Example 2}\label{sec5.2}
In this example, we focused on the stability of the proposed HDG scheme when $\ep$ is small. Set the same domain $\Omega=(0,1)^2$ and $T=1$. The initial conditions and the external forces are given by
\begin{equation}
	\begin{aligned}
		&\bm{u}^0(x)=\left(-\frac{\partial \psi }{\partial x_2}, \frac{\partial \psi }{\partial x_1}\right),
		\quad\text{with }\psi(x,t)=-200(x_1(1-x_1)x_2(1-x_2))^2,\\
		&\bm{C}^0(x)=\left(
        \begin{array}{cc}
        \frac{\sqrt{2}}{2} & 0 \\
        0 & \frac{\sqrt{2}}{2}\\
        \end{array} \right),\quad 
        \bm{f}(x)=\left(-70(x_2-\frac{1}{2}),70(x_1-\frac{1}{2})\right),\quad\bm{F}(x)=0.
	\end{aligned}
\end{equation}
We choose $\epsilon=10^{-4}$ and $\nu=10^{-2}$. Note that both simulations below are computed with $h=2^{-6}$ and $\tau = \frac{1}{100}$.

\begin{figure}[htbp] 
\begin{tabular}{ccc}
\centering
\includegraphics[width=0.3\textwidth]{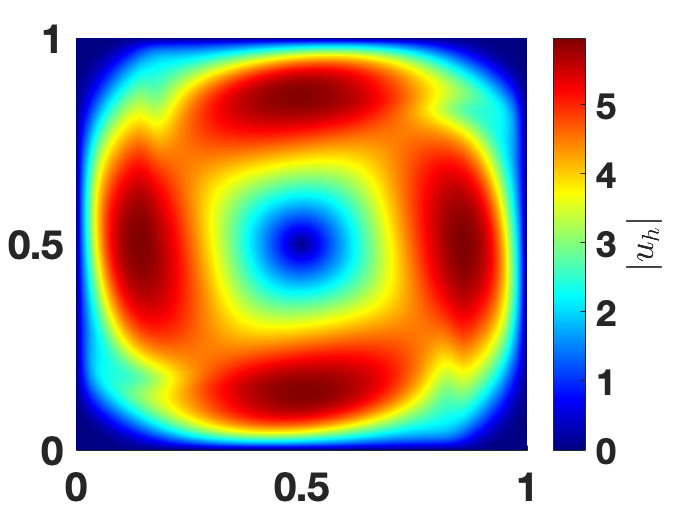} &
\includegraphics[width=0.3\textwidth]{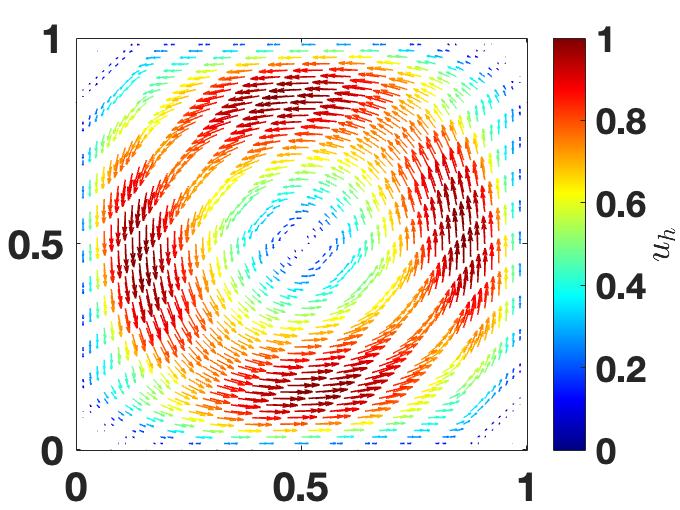} &
\includegraphics[width=0.3\textwidth]{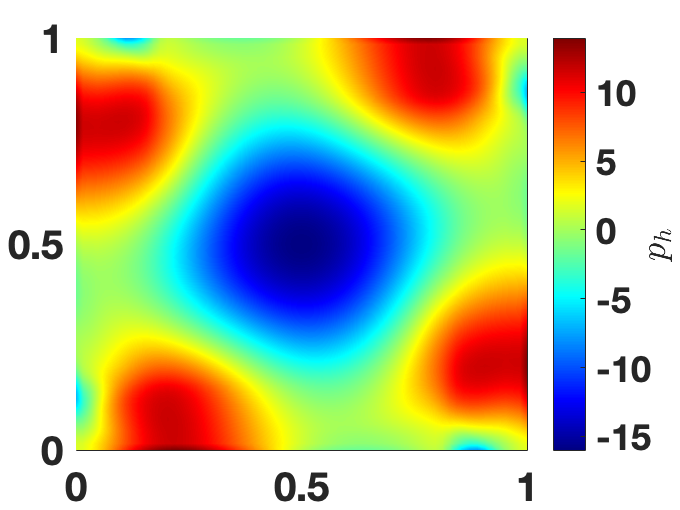} 
\\
(a) $|\bm u_h^N|$ &
(b) $\bm u_h^N$ &
(c) $p_h^N$
\\
\includegraphics[width=0.3\textwidth]{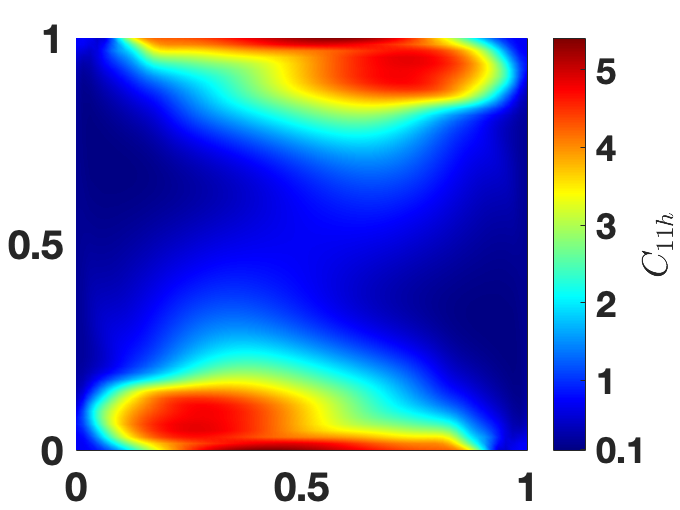} &
\includegraphics[width=0.3\textwidth]{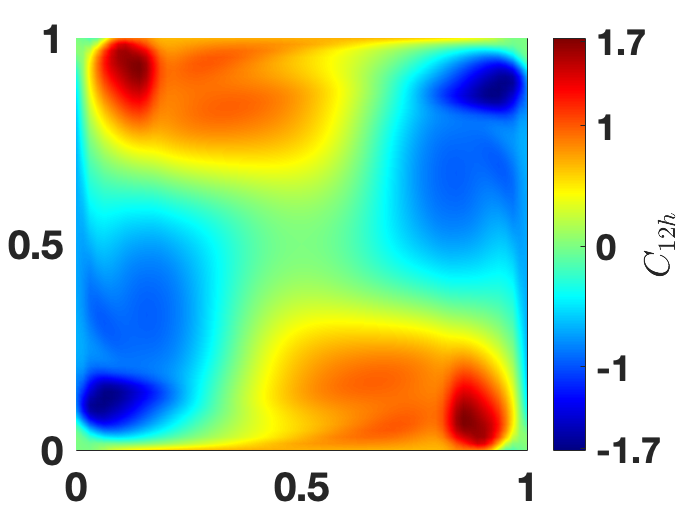} &
\includegraphics[width=0.3\textwidth]{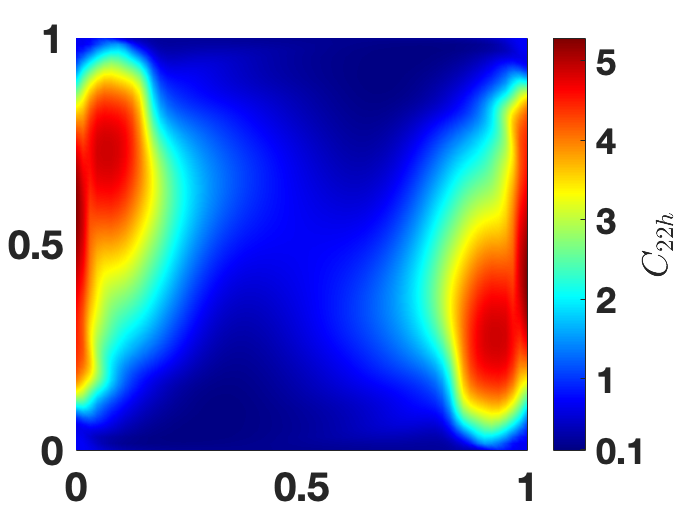} 
\\
(d) $C_{11h}^N$ &
(e) $C_{12h}^N$ &
(f) $C_{22h}^N$
\\
\end{tabular}
\caption{The numerical solution $(\bm{u}_h,p_h,\bm{C}_h)$ of the HDG scheme at $t=1$ with $(h=2^{-6},\tau=\frac{1}{100})$}\label{fig-final}
\end{figure}

First we apply the HDG scheme with $k=1$ and the penalty coefficients $\alpha=600$ and $\beta=600$ for numerical simulation.
The pictures of the discrete solution $(\bm{u}_h,p_h,\bm{C}_h)$ are displayed in Figure~\ref{fig-final} at $t=1$.

\begin{figure}[htbp] 
\begin{tabular}{ccc}
\centering
\includegraphics[width=0.3\textwidth]{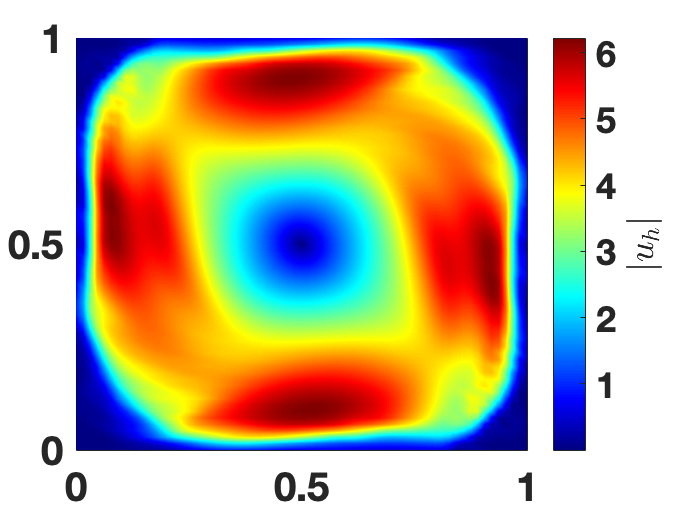} &
\includegraphics[width=0.3\textwidth]{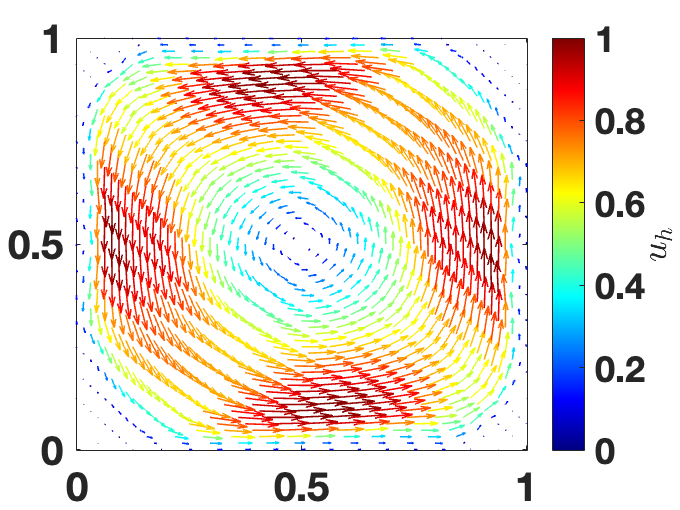} &
\includegraphics[width=0.3\textwidth]{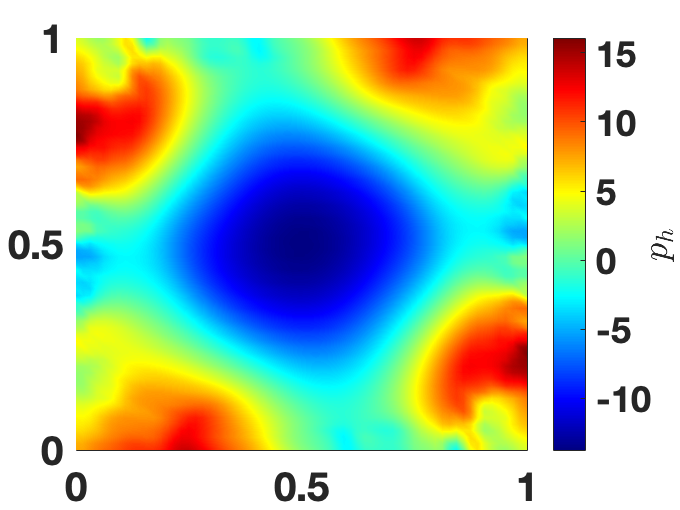} 
\\
(a) $|\bm u_h^N|$ &
(b) $\bm u_h^N$ &
(c) $p_h^N$
\\
\includegraphics[width=0.3\textwidth]{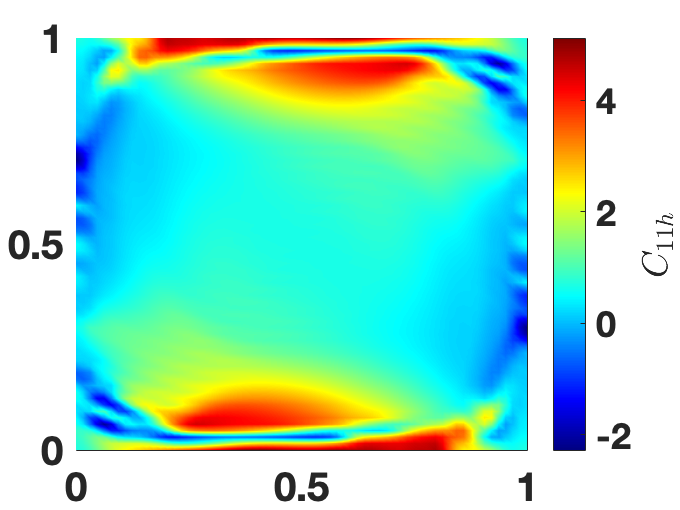} &
\includegraphics[width=0.3\textwidth]{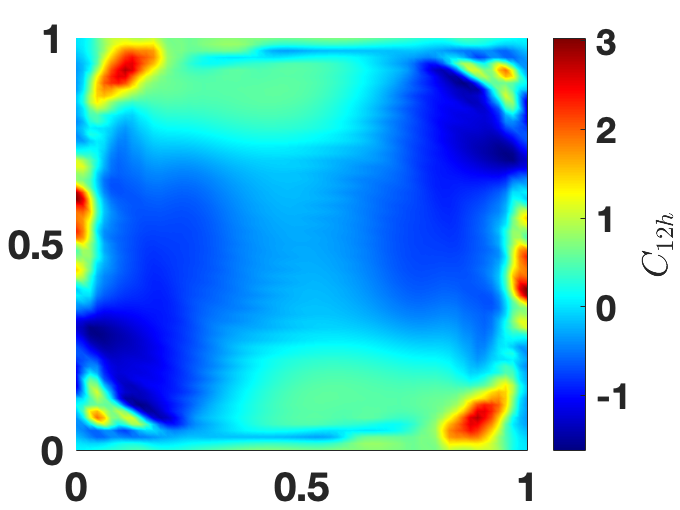} &
\includegraphics[width=0.3\textwidth]{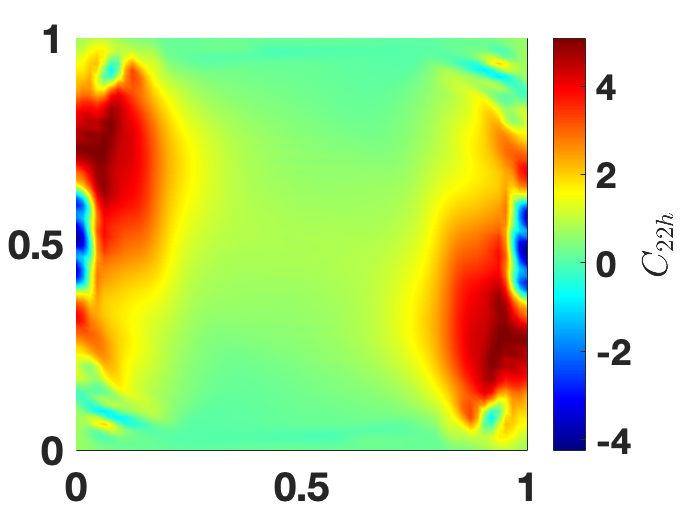} 
\\
(d) $C_{11h}^N$ &
(e) $C_{12h}^N$ &
(f) $C_{22h}^N$
\\
\end{tabular}
\caption{The numerical solution $(\bm{u}_h,p_h,\bm{C}_h)$ of the ordinary FEM scheme at $t=0.55$ with $(h=2^{-6},\tau=\frac{1}{100})$}\label{fig2}
\end{figure}

\begin{figure}[htbp] 
\begin{tabular}{ccc}
\centering
\includegraphics[width=0.3\textwidth]{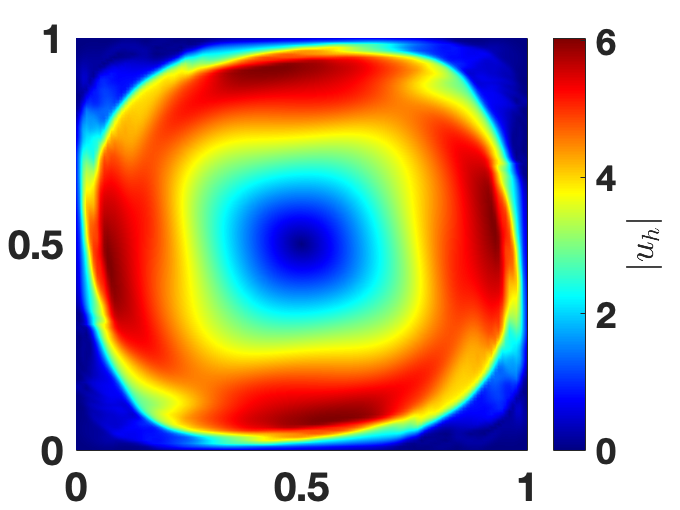} &
\includegraphics[width=0.3\textwidth]{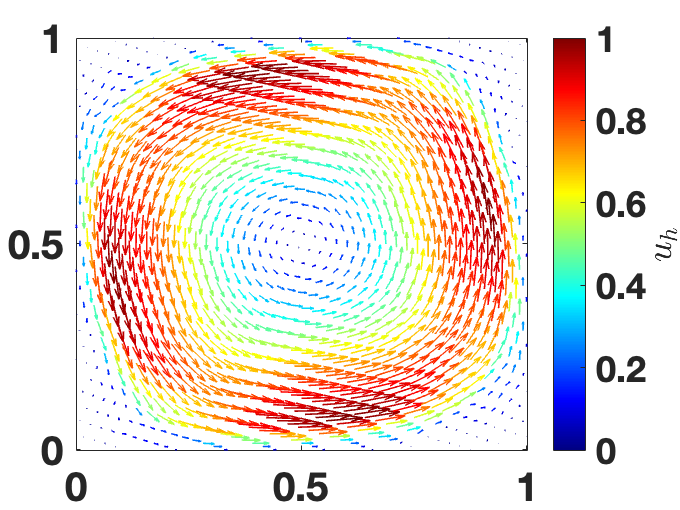} &
\includegraphics[width=0.3\textwidth]{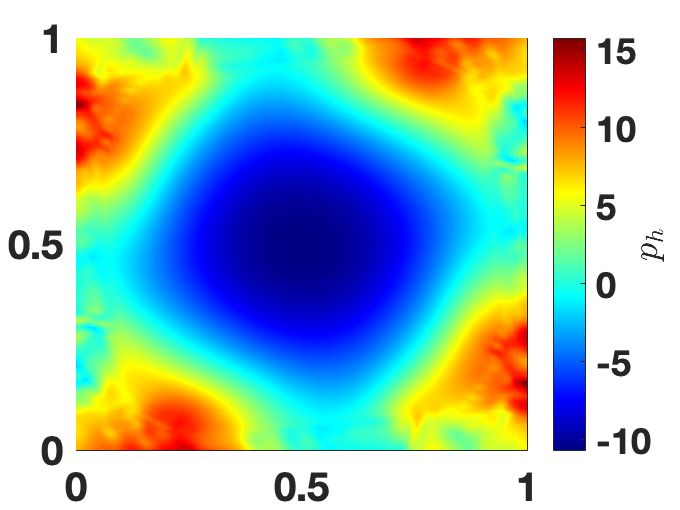} 
\\
(a) $|\bm u_h^N|$ &
(b) $\bm u_h^N$ &
(c) $p_h^N$
\\
\includegraphics[width=0.3\textwidth]{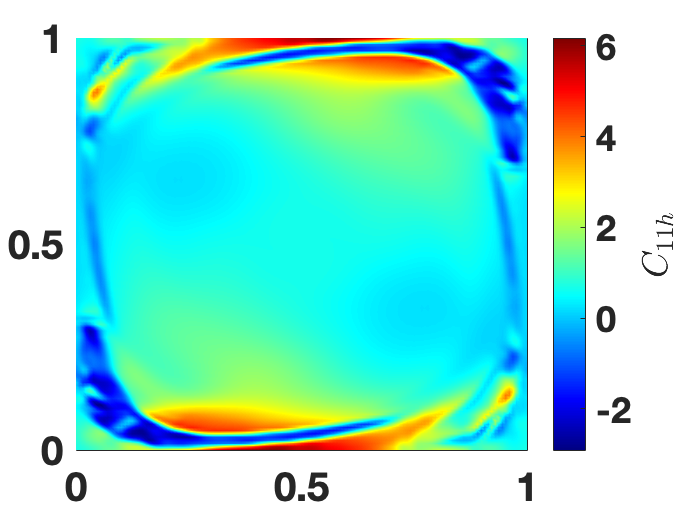} &
\includegraphics[width=0.3\textwidth]{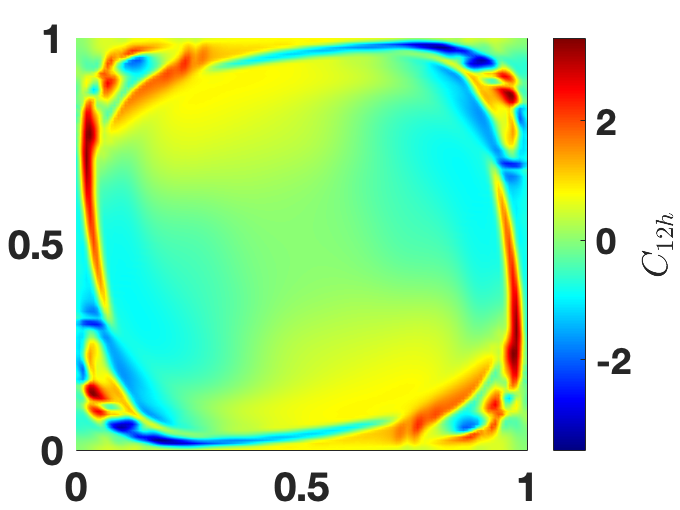} &
\includegraphics[width=0.3\textwidth]{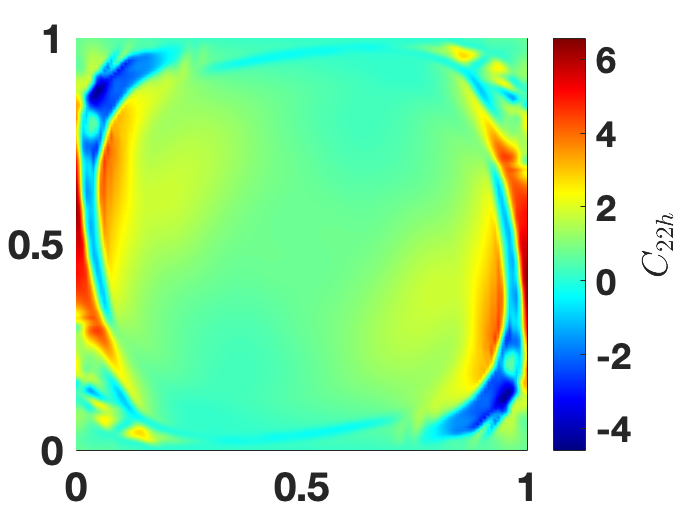} 
\\
(d) $C_{11h}^N$ &
(e) $C_{12h}^N$ &
(f) $C_{22h}^N$
\\
\end{tabular}
\caption{The numerical solution $(\bm{u}_h,p_h,\bm{C}_h)$ of the ordinary FEM scheme at $t=0.45$ with $(h=2^{-7},\tau=\frac{1}{100})$}\label{fig3}
\end{figure}

\begin{figure}[htbp] 
\begin{tabular}{ccc}
\centering
\includegraphics[width=0.3\textwidth]{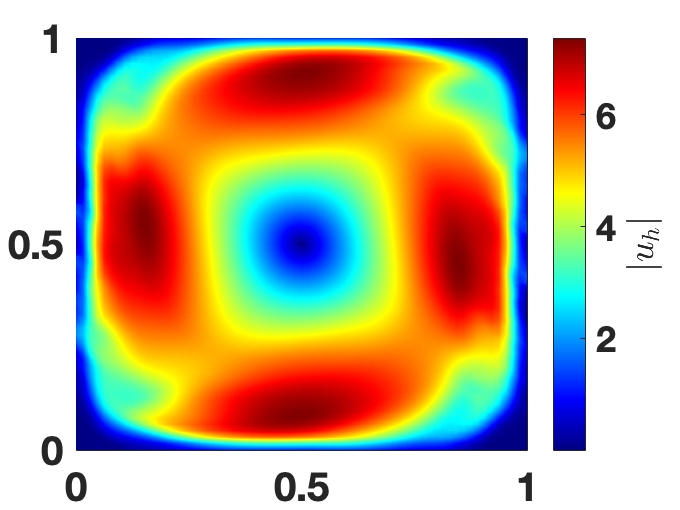} &
\includegraphics[width=0.3\textwidth]{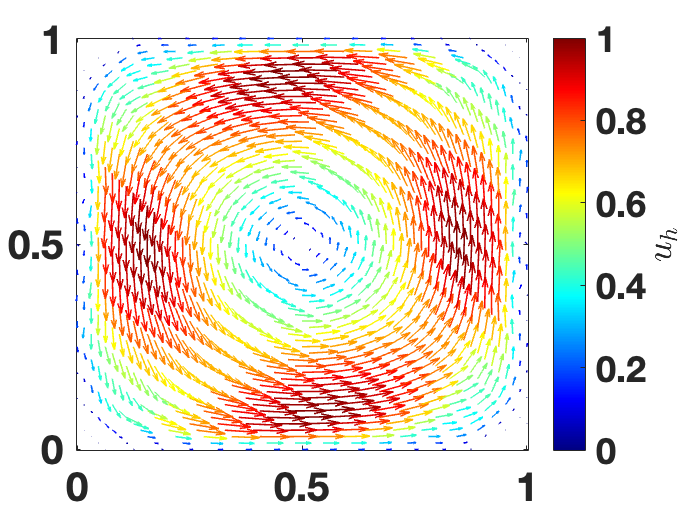} &
\includegraphics[width=0.3\textwidth]{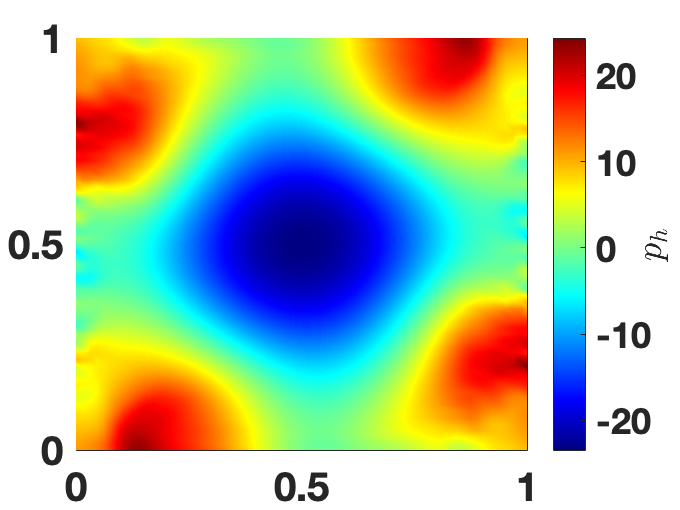} 
\\
(a) $|\bm u_h^N|$ &
(b) $\bm u_h^N$ &
(c) $p_h^N$
\\
\includegraphics[width=0.3\textwidth]{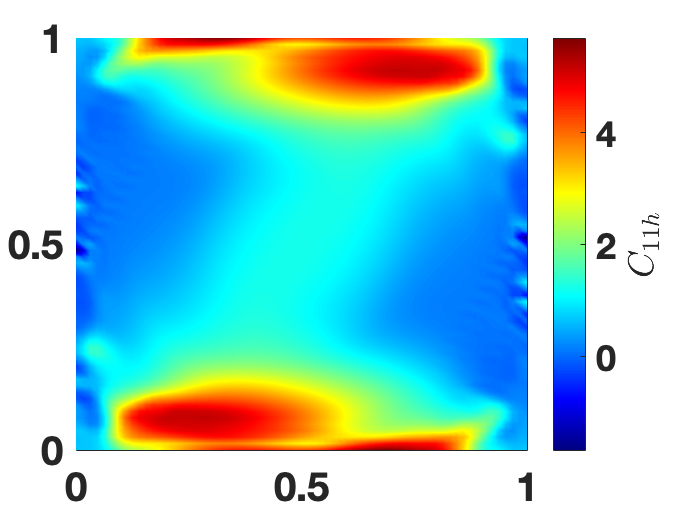} &
\includegraphics[width=0.3\textwidth]{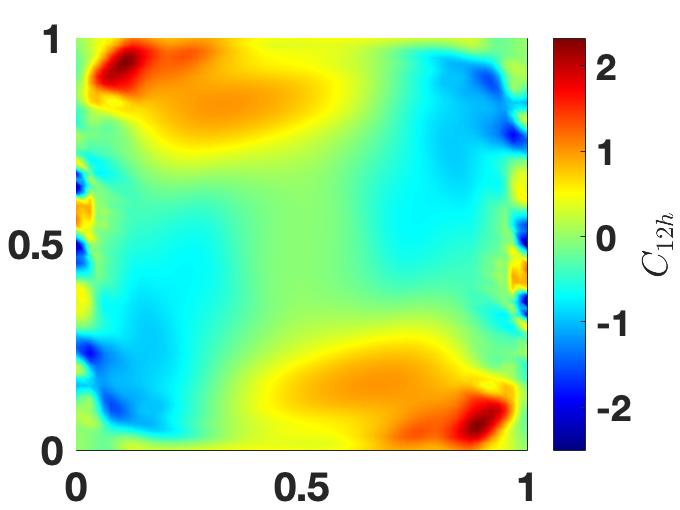} &
\includegraphics[width=0.3\textwidth]{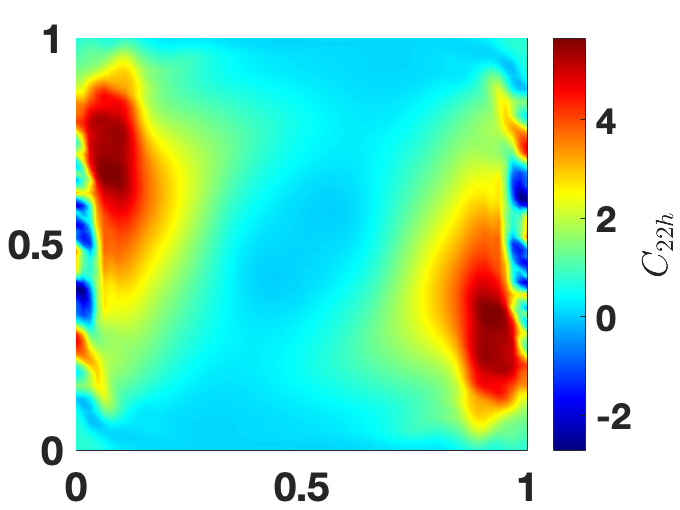} 
\\
(d) $C_{11h}^N$ &
(e) $C_{12h}^N$ &
(f) $C_{22h}^N$
\\
\end{tabular}
\caption{The numerical solution $(\bm{u}_h,p_h,\bm{C}_h)$ of the ordinary FEM scheme at $t=1$ with $(h=2^{-6},\tau=\frac{1}{120})$}\label{fig4}
\end{figure}

As a comparison, the general FEM using P1-bubble/P1/P1-element is applied to solve the same problem. We use the backward Euler method for the time discretization, which is the same as \eqref{eq:Full}, and implicit-explicit time-stepping scheme to deal with the nonlinearity.
We find that the discrete solution $(\bm{u}_h,p_h,\bm{C}_h)$ becomes unstable after $t = 0.55$ as shown in Figure~\ref{fig2}, indicating that our proposed HDG scheme is more stable than the finite element method for small $\ep$ values.
% In particular, the HDG scheme performed well at $t = 0.55$ (see Figure~\ref{fig2}), 
%
Then, we carry out the simulation using a fine mesh $h=2^{-7}$ and discover that the instability appears around $t=0.45$ (see Figure~\ref{fig3}).
On the other hand, we apply the FEM scheme with a smaller time-step size $\tau = \frac{1}{120}$. 
The result is plotted in Figure~\ref{fig4}, which show the instability of the numerical solution.

\begin{figure}[htbp] 
\begin{tabular}{cc}
\centering
\includegraphics[width=0.3\textwidth]{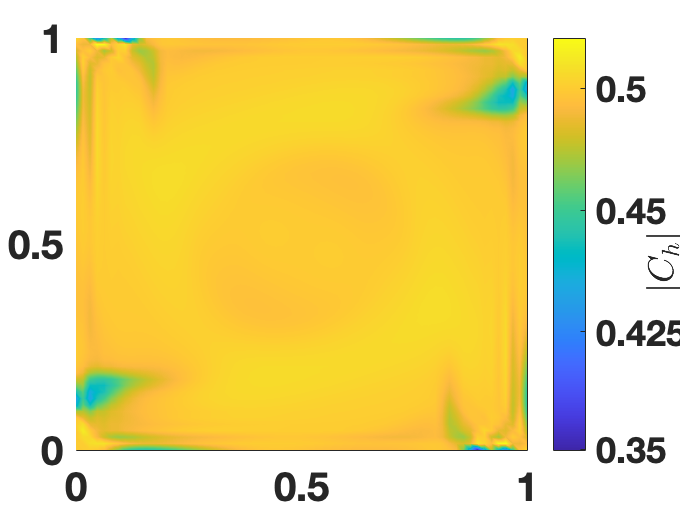} &
\includegraphics[width=0.3\textwidth]{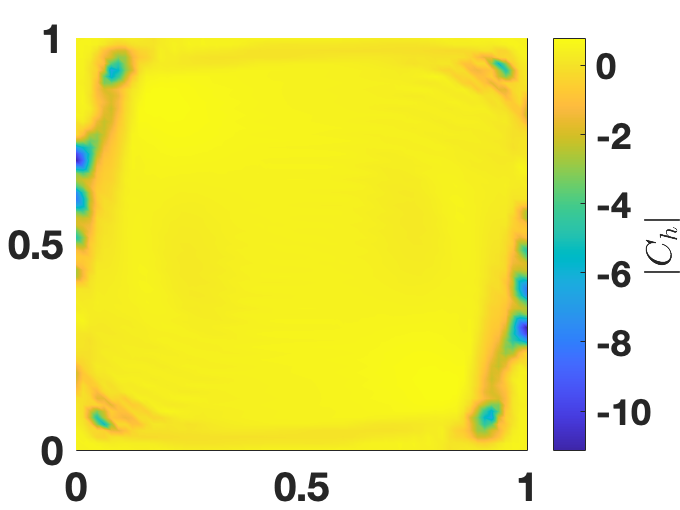} 
\\
(a) HDG scheme with $(h=2^{-6},\tau=\frac{1}{100})$&
(b) FEM scheme with $(h=2^{-6},\tau=\frac{1}{100})$ 
\\
\includegraphics[width=0.3\textwidth]{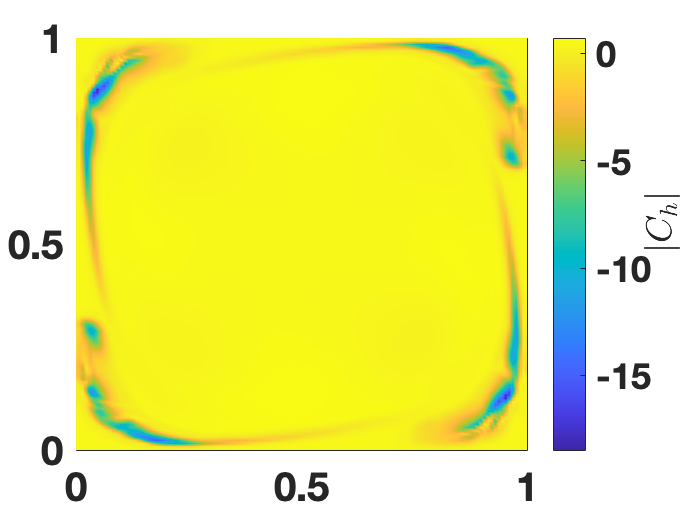} &
\includegraphics[width=0.3\textwidth]{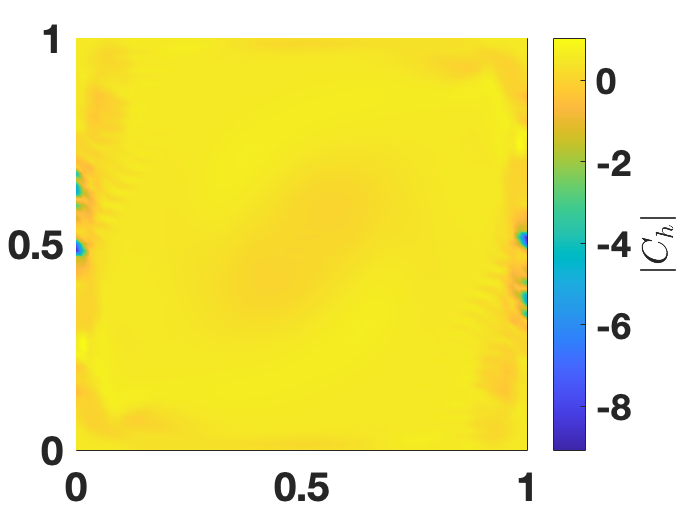} 
\\
(c) FEM scheme with $(h=2^{-7},\tau=\frac{1}{100})$ &
(d) FEM scheme with $(h=2^{-6},\tau=\frac{1}{120})$
\\
\end{tabular}
\caption{The determinant of $\bm{C}_h$ for HDG scheme and the ordinary FEM scheme ((b)--(d))}\label{fig-C}
\end{figure}

Additionally, we plot the determinant of discrete conformation tensor in Figure~\ref{fig-C}. 
For the HDG scheme, we see that both $\bm C_{11h}$ and $\mathrm{det}(\bm C_h)$ are positive functions (see Figure~\ref{fig-final} (d) and Figure~\ref{fig-C} (a)), which indicates that $\bm C_h$ is positive definite. 
However, for ordinary FEM, we observe negative values of $\bm C_{11h}$ (see Figures~\ref{fig2} (d), \ref{fig3} (d) and \ref{fig4} (d)) and the determinant also has negative values (see Figure~\ref{fig-C} (b)--(d)), which implies that $\bm C_h$ of the FEM loses the positive-definiteness.  

%%%%%%%%%%%%%%%%%%%%%%%%%%%%%%%%%%%%
% HDG final time
% %%%%%%%%%%%%%%%%%%%%%%%%%%%%%%%%%%%%
% %%%%%%%%%%%%%%%%%%%%%%%%%%%%%%%%%%%%
%%%%%%%%%%%%%%%%%%%%%%%%%%%%%%%%%%%%
% FEM t = 0.68s
%%%%%%%%%%%%%%%%%%%%%%%%%%%%%%%%%%%%
%%%%%%%%%%%%%%%%%%%%%%%%%%%%%%%%%%%%
%%%%%%%%%%%%%%%%%%%%%%%%%%%%%%%%%%%%
%%%%%%%%%%%%%%%%%%%%%%%%%%%%%%%%%%%%
\bibliographystyle{plain}
\bibliography{HDG}

\begin{thebibliography}{10}

\bibitem{Becker2023}
F.~Becker, K.~Rauthmann, L.~Pauli, and P.~Knechtgesa.
\newblock {An eigenvalue-free implementation of the log-conformation formulation}.
\newblock {\em J. Non-Newton. Fluid Mech.}, 322:105133, 2023.

\bibitem{Boyaval2009}
S.~Boyaval, T.~Leli\`{e}vre, and C.~Mangoubi.
\newblock {Free-energy-dissipative schemes for the Oldroyd-B model}.
\newblock {\em ESAIM: M2AN}, 43(3):523--561, 2009.

\bibitem{FEM2008}
S.~C. Brenner and L.~R. Scott.
\newblock {\em {The Mathematical Theory of Finite Element Methods}}.
\newblock Springer, New York, 2008.

\bibitem{Brunk2022}
A.~Brunk, Y.~Lu, and M.~Luk\'a\v{c}ov\'a-Medvid'ov\'a.
\newblock {Existence, regularity and weak-strong uniqueness for the three-dimensional Peterlin viscoelastic model}.
\newblock {\em Commun. Math. Sci.}, 20(1):201--230, 2022.

\bibitem{Aycil2017}
A.~Cesmelioglu, B.~Cockburn, and W.~Qiu.
\newblock {Analysis of a hybridizable discountinuous Galerkin method for the steady-state incompressible Navier-Stokes equations}.
\newblock {\em Math. Comput.}, 86(306):1643--1670, 2017.

\bibitem{Cockburn20122}
B.~Cockburn and J.~Cui.
\newblock Divergence-free {HDG} methods for the vorticity-velocity formulation of the {Stokes} problem.
\newblock {\em J. Sci. Comput.}, 52(1):256--270, 2012.

\bibitem{Cockburn2008}
B.~Cockburn, B.~Dong, and J.~Guzm\'{a}n.
\newblock {A superconvergent LDG-hybridizable Galerkin method for second-order elliptic problems}.
\newblock {\em Math. Comput.}, 77(264):1887--1916, 2008.

\bibitem{Cockburn20092}
B.~Cockburn, J.~Gopalakrishnan, and R.~Lazarov.
\newblock Unified hybridization of discontinuous {Galerkin}, mixed, and continuous {Galerkin} methods for second order elliptic problems.
\newblock {\em SIAM J. Numer. Anal.}, 47(2):1319--1365, 2009.

\bibitem{Cockburn2011}
B.~Cockburn, J.~Gopalakrishnan, N.~C. Nguyen, J.~Peraire, and F.~J. Sayas.
\newblock {Analysis of HDG methods for Stokes flow}.
\newblock {\em Math. Comput.}, 80(274):723--760, 2011.

\bibitem{HDG2019}
S.~Du and F.~J. Sayas.
\newblock {\em {An Invitation to the Theory of the Hybridizable Discontinuous Galerkin Method, Projections, Estimates, Tools}}.
\newblock Springer, Cham, 2019.

\bibitem{Ervin2004}
V.~J. Ervin and N.~Heuer.
\newblock {Approximation of time-dependent, viscoelastic fluid flow: Crank--Nicolson, finite element approximation}.
\newblock {\em Numer. Methods Partial Differ. Equ.}, 20(2):248--283, 2004.

\bibitem{Fattal2004}
R.~Fattal and R.~Kupferman.
\newblock {Constitutive laws for the matrix-logarithm of the conformation tensor}.
\newblock {\em J. Non-Newton. Fluid Mech.}, 123:281--285, 2004.

\bibitem{Guillope90}
C.~Guillop\'e and J.~C. Saut.
\newblock Global existence and one-dimensional nonlinear stability of shearing motions of viscoelastic fluids of {Oldroyd} type.
\newblock {\em ESAIM: M2AN}, 24:369--401, 1990.

\bibitem{Han23}
W.~Han, Y.~Jiang, and Z.~Miao.
\newblock On a second-order decoupled time-stepping scheme for solving a finite element problem for the approximation of {Peterlin} viscoelastic model.
\newblock {\em Comput. Math. Appl}, 142:48--63, 2023.

\bibitem{Hu2007}
D.~Hu and T.~Leli\`{e}vre.
\newblock {New entropy estimates for Oldroyd-B and related models}.
\newblock {\em Commun. Math. Sci}, 5(4):909--916, 2007.

\bibitem{Jiang2018}
Y.~Jiang and Y.~Yang.
\newblock Semi-discrete {Galerkin} finite element method for the diffusive {Peterlin} viscoelastic model.
\newblock {\em Comput. Methods Appl. Math.}, 18(2):275--296, 2018.

\bibitem{Kirk2019}
K.~L.~A. Kirk and S.~Rhebergen.
\newblock {Analysis of a pressure-robust hybridized discontinuous Galerkin method for the stationary Navier-Stokes equations}.
\newblock {\em J. Sci. Comput.}, 81(2):881--897, 2019.

\bibitem{Labeur2007}
J.~R. Labeur and G.~N. Wells.
\newblock A {Galerkin} interface stabilisation method for the advection-diffusion and incompressible {Navier}-{Stokes} equations.
\newblock {\em Comput. Methods Appl. Mech. Engrg.}, 196(49-52):4985--5000, 2007.

\bibitem{Lions2000}
P.~L. Lions and N.~Masmoudi.
\newblock {Global solutions for some Oldroyd models of non-Newtonian flows}.
\newblock {\em Chin. Annal Math.}, 21(2):131--146, 2000.

\bibitem{Lukacova15}
M.~Luk\'a\v{c}ov\'a-Medvid'ov\'a, H.~Mizerov\'a, and S.~Ne\v{c}asov\'a.
\newblock Global existence and uniqueness result for the diffusive {Peterlin} viscoelastic model.
\newblock {\em NONLINEAR ANAL-THEOR}, pages 154--170, 2015.

\bibitem{Medivdova2017}
M.~Luk\'a\v{c}ov\'a-Medvid'ov\'a, H.~Mizerov\'a, S.~Ne\v{c}asov\'a, and M.~Renardy.
\newblock Global existence result for the generalized {Peterlin} viscoelastic model.
\newblock {\em SIAM J. Math. Anal.}, 49(4):2950--2964, 2017.

\bibitem{M1}
M.~Luk\'a\v{c}ov\'a-Medvid'ov\'a, H.~Mizerov\'a, H.~Notsu, and M.~Tabata.
\newblock {Numerical analysis of the Oseen-type Peterlin viscoelastic model by the stabilized Lagrange--Galerkin method. Part I: A nonlinear scheme}.
\newblock {\em ESAIM: M2AN}, 51(5):1637--1661, 2017.

\bibitem{M2}
M.~Luk\'a\v{c}ov\'a-Medvid'ov\'a, H.~Mizerov\'a, H.~Notsu, and M.~Tabata.
\newblock Numerical analysis of the {Oseen}-type {Peterlin} viscoelastic model by the stabilized {Lagrange}--{Galerkin} method. {Part} {II}: {A} linear scheme.
\newblock {\em ESAIM: M2AN}, 51(5):1663--1689, 2017.

\bibitem{Hana2015}
H.~Mizerov\'a.
\newblock {\em {Analysis and Numerical Solution of the Peterlin Viscoelastic Model}}.
\newblock PhD thesis, Universita\"tsbibliothek, Mainz, 2015.

\bibitem{Najib1995}
K.~Najib and D.~Sandri.
\newblock {On a decoupled algorithm for solving a finite element problem for the approximation of viscoelastic fluid flow}.
\newblock {\em Numer. Math.}, 72(2):223--238, 1995.

\bibitem{Nguyen2011}
N.~C. Nguyen, J.~Peraire, and B.~Cockburn.
\newblock {An implicit high-order hybridizable discontinuous Galerkin method for the incompressible Navier--Stokes equations}.
\newblock {\em J. Comput. Phys.}, 230(4):1147--1170, 2011.

\bibitem{Oikawa2014}
I.~Oikawa.
\newblock {Hybridized discontinuous Galerkin method for convection-diffusion problems}.
\newblock {\em Jpn. J. Ind. Appl. Math.}, 31(2):335--354, 2014.

\bibitem{Oikawa2015}
I.~Oikawa.
\newblock A hybridized discontinuous {Galerkin} method with reduced stabilization.
\newblock {\em J. Sci. Comput.}, 65:327--340, 2015.

\bibitem{Qiuwf2018}
W.~Qiu, J.~Shen, and K.~Shi.
\newblock {An HDG method for linear elasticity with strong symmetric stresses}.
\newblock {\em Math. Comput.}, 87(309):69--93, 2018.

\bibitem{Ravindran2020}
S.~Ravindran.
\newblock {Analysis of a second-order decoupled time-stepping scheme for transient viscoelastic flow}.
\newblock {\em Int. J. Numer. Anal. Model}, 17(1):87--109, 2020.

\bibitem{Sander2017}
S.~Rhebergen and G.~N. Wells.
\newblock Analysis of a hybridized/interface stabilized finite element method for the {Stokes} equations.
\newblock {\em SIAM J. Numer. Anal.}, 55(4):1982--2003, 2017.

\bibitem{Sander2018}
S.~Rhebergen and G.~N. Wells.
\newblock A hybridizable discontinuous {G}alerkin method for the {Navier}--{Stokes} equations with pointwise divergence-free velocity field.
\newblock {\em J. Sci. Comput.}, 76(3):1484--1501, 2018.

\bibitem{Sander20182}
S.~Rhebergen and G.~N. Wells.
\newblock Preconditioning of a hybridized discontinuous {Galerkin} finite element method for the {Stokes} equations.
\newblock {\em J. Sci. Comput.}, 77:1936--1952, 2018.

\bibitem{H2016}
A.~H. Sequeira.
\newblock {\em {Hemorheology: Non-Newtonian Constitutive Models for Blood Flow Simulations, Non-Newtonian Fluid Mechanics and Complex Flows}}.
\newblock Springer, Cham, 2018.

\bibitem{Wittschieber2022}
S.~Wittschieber, L.~Demkowicz, and M.~Behr.
\newblock {Stabilized finite element methods for a fully-implicit logarithmic reformulation of the Oldroyd-B constitutive law}.
\newblock {\em J. Non-Newton. Fluid Mech.}, 306:104838, 2022.

\bibitem{Xia2023}
L.~Xia and G.~Zhou.
\newblock A linearizing-decoupling finite element method with stabilization for the {Peterlin} viscoelastic model.
\newblock {\em Jpn. J. Ind. Appl. Math.}, 2023.

\bibitem{Zhang2013}
M.~Zhang, I.~Lashgari, T.~A. Zaki, and L.~Brandt.
\newblock Linear stability analysis of channel flow of viscoelastic {Oldroyd-B} and {FENE-P} fluids.
\newblock {\em J. Fluid Mech.}, 737:249--279, 2013.

\bibitem{Zhangyz2024}
Y.~Zhang, X.~Yong, and X.~Du.
\newblock Numerical analysis of time filter method for the stabilized incompressible diffusive {P}eterlin viscoelastic fluid model.
\newblock {\em Comput. Math. Appl}, 168:239--253, 2024.

\bibitem{Zheng2017}
H.~Zheng, J.~Yu, and L.~Shan.
\newblock {Unconditional error estimates for time dependent viscoelastic fluid flow}.
\newblock {\em Appl. Numer. Math.}, 119(1):1--17, 2017.

\end{thebibliography}
\end{document}